\documentclass[12pt]{amsart}

\usepackage{amsmath}
\usepackage{amscd}
\usepackage{amssymb}
\usepackage{enumerate}
\usepackage{amsfonts}
\usepackage{graphicx}
\usepackage[all]{xypic}
\usepackage{mathrsfs}

\numberwithin{equation}{section}

\newcommand{\eq}[1]{\eqref{eq:#1}}

\newcommand{\al}{\alpha}  \newcommand{\ga}{\gamma}
 \newcommand{\de}{\delta} \newcommand{\De}{\Delta}
\newcommand{\la}{\lambda} 
  \newcommand{\na}{\nabla}
\newcommand{\eps}{\epsilon} 
 \newcommand{\te}{\theta}

\newcommand{\bN}{{\mathbb N}}

\newcommand{\bR}{{\mathbb R}}

\newcommand{\cE}{{\mathcal E}}

\newcommand{\dd}{\partial}

\newcommand{\abs}[1]{\left\lvert #1\right\rvert}
\newcommand{\norm}[1]{\left\lvert\!\left\lvert #1\right\rvert\!\right\rvert}

\newcommand{\supp}{\operatorname{supp}}

\newcommand{\sign}{\operatorname{sign}}

\newcommand{\mbr}{\medbreak}

\newtheorem{thm}{Theorem}[section]
\newtheorem{cor}[thm]{Corollary}
\newtheorem{lem}[thm]{Lemma}
\newtheorem{prop}[thm]{Proposition}

\theoremstyle{remark}
\newtheorem{rem}[thm]{Remark}

\newtheorem{defin}[thm]{Definition}

\newcommand{\bbe}{\begin{equation}}
\newcommand{\ee}{\end{equation}}
\newcommand{\bbea}{\begin{eqnarray*}}
\newcommand{\eea}{\end{eqnarray*}}

\newcommand{\infint}{(-\infty,+\infty)}

\newcommand{\td}{\tilde{\delta}}
\newcommand{\wi}{with maximal interval of existence $I$}
\newcommand{\ssub}{\subset\subset}
\newcommand{\h}{\dot{H^1}}
\newcommand{\hm}{\dot{H}}
\newcommand{\w}[1]{\overset{.}{W}_x^{\frac{1}{2},\,#1}}
\newcommand{\sd}{\sqrt{-\De}}
\newcommand{\hl}{\h\times L^2}
\newcommand{\rt}{(u_0,u_1)}
\newcommand{\rtt}{(\tilde{u}_0,\tilde{u}_1)}

\newcommand{\sde}{\overline{s}_\de}

\newcommand{\dv}{\operatorname{div}}

\title[energy-critical focusing Wave]{Global well-posedness, scattering
and blow-up for the energy critical focusing non-linear wave
equation}

\author[C.E. Kenig and F. Merle]{Carlos E. Kenig\address{\hskip-1.15em Carlos E.
Kenig: \hfill\newline Department of Mathematics, \hfill\newline
University of Chicago, \hfill\newline Chicago, IL 60637
\hfill\newline USA} \and Frank Merle
\address{\hskip-1.15em Frank Merle: \hfill\newline D\'epartement de Math\'ematiques, \hfill\newline
Universit\'e de Cergy-Pontoise, \hfill\newline Pontoise,
\hfill\newline 95302 Cergy-Pontoise, \hfill\newline FRANCE} }

\thanks{The first author was supported in part by NSF and the second one in part by CNRS.
Part of this research was carried out during visits of the second
author to the University of Chicago and I.H.E.S. and of the first 
author to Paris XIII.}
\begin{document}

\maketitle

\section{Introduction}\label{s:1}

In this paper we consider the energy critical non-linear wave
equation
\[
\begin{cases}
\dd^2_t u - \De u = \pm \abs{u}^{\frac{4}{N-2}} u & (x,t)\in\bR^N\times\bR
\\
u\bigl\lvert_{t=0} = u_0 \in \h(\bR^N) & \\
\dd_t u \bigl\lvert_{t=0} = u_1 \in L^2(\bR^N)
\end{cases}
\]
Here the $-$ sign corresponds to the defocusing problem, while the
$+$ sign corresponds to the focusing problem. The theory of the
local Cauchy problem (CP) for this equation was developed in many
papers, see for instance \cite{23}, \cite{9}, \cite{20}, \cite{25},
\cite{26}, \cite{27}, \cite{K} etc. In particular, one can show that if
$\norm{(u_0,u_1)}_{\h\times L^2}\leq\de$, $\de$ small, there exists
a unique solution with $(u,\dd_t u)\in C\bigl(\bR;\h(\bR^N)\times
L^2(\bR^N)\bigr)$ with the norm
$$\norm{u}_{L_{xt}^{\frac{2(N+1)}{N-2}}}<\infty$$
$\bigl($i.e., the solution
scatters in $\h(\bR^N)\times L^2(\bR^N)${}$\bigr)$. See section 2 of
this paper for a review and an update of the results.

\mbr


In the defocusing case, Struwe \cite{30} in the radial case, when $N=3$, Grillakis
\cite{11} in the general case when $N=3$, and then Grillakis \cite{12}, Shatah-Struwe
\cite{24}, \cite{25}, \cite{26} (and others \cite{K}) in higher dimensions, proved that this
also holds for any $(u_0,u_1)$ with $\norm{(u_0,u_1)}_{\h\times L^2}<\infty$ and that,
(for $3\leq N\leq 5$) for more regular $(u_0,u_1)$ the solution preserves the smoothness for
all time. This topic has been the subject of intense investigation. See the recent
work of Tao \cite{32} for a recent installment in it and further references.

\mbr

In the focusing case, these results do not hold. In fact, the
classical identity
 \bbe
\frac{d^2}{dt^2} \int \abs{u(x,t)}^2 = 2\left[ \int (\dd_t u)^2 -
\abs{\nabla u}^2 - \abs{u(t)}^{\frac{2N}{N-2}} \right]
 \ee
(see the work of H.~Levine \cite{19} and also sections 3 and 5) was
used by Levine \cite{19} to show that if $(u_0,u_1)\in H^1\times
L^2$ is such that
\[
E\bigl((u_0,u_1)\bigr) = \int \frac{1}{2} \abs{\nabla u_0}^2 +
\frac{1}{2} \abs{u_1}^2 - \frac{(N-2)}{2N} \abs{u_0}^{\frac{2N}{N-2}}<0,
\]
the solution must break down in finite time. Moreover,
\[
W(x) = W(x,t) = \frac{1}{\left(
1+\frac{\abs{x}^2}{N(N-2)}\right)^{\frac{(N-2)}{2}}}
\]
is in $\h(\bR^N)$ and solves the elliptic equation
\[
\De W + \abs{W}^{\frac{4}{N-2}} W = 0,
\]
so that scattering cannot always occur even for global (in time)
solutions.

\mbr

In this paper we initiate the detailed study of the focusing case
(see also \cite{18} for an interesting recent work in this 
direction). We show:
\begin{thm}
Let $(u_0,u_1)\in\h\times L^2$, $3\leq N\leq 5$. Assume that
$E\bigl((u_0,u_1)\bigr)<E\bigl((W,0)\bigr)$. Let $u$ be the
corresponding solution of the Cauchy problem, with maximal interval
of existence \hfill\newline $I=\bigl( -T_-(u_0,u_1), T_+(u_0,u_1)
\bigr)$. $($See Definition \ref{d:2.13}.$)$ Then:
\begin{enumerate}[i$)$]
\item If $\int\abs{\na u_0}^2 < \int \abs{\na W}^2$, then

$$I=(-\infty,+\infty) \ \mbox{and} \
\norm{u}_{L^{\frac{2(N+1)}{N-2}}_{xt}}<\infty.$$

\item If $\int\abs{\na u_0}^2 > \int \abs{\na W}^2$, then
$$T_+(u_0,u_1)<+\infty, \ T_-(u_0,u_1)<+\infty.$$
\end{enumerate}
\end{thm}
 \setcounter{section}{1}

\mbr

Our proof follows the new point of view into these problems that we
introduced in \cite{15}, where we obtained the corresponding result
for the energy critical non-linear Schr\"odinger equation for radial
data. In section 3 we prove some elementary variational estimates
which yield the necessary coercivity for our arguments and which
follows from arguments in \cite{15}. In section 4, using the work of
Bahouri-Gerard \cite{4} and the concentration compactness argument
from \cite{15} we produce a ``critical element'' for which
scattering fails and which enjoys a compactness property because of
its criticality. (Propositions \ref{p:4.1} and \ref{p:4.2}.) At this point, we show
a crucial orthogonality property of ``critical elements'' related to a second conservation law in the energy space
(Proposition \ref{p:4.14} and Proposition \ref{p:4.21}) which exploits the finite
speed of propagation for the wave equation and its Lorentz
invariance. This is the extra ingredient that allows us to go beyond
the radial case as in \cite{15}. In sections 5 and 6 we prove a
rigidity theorem (Theorem \ref{t:5.1}), which allows us to conclude the
argument. The first case of the rigidity theorem deals with infinite
time of existence. This uses localized conservations laws of the
type (1.1) and related ones, very much in the spirit of the
corresponding localized virial identity used in \cite{15}. The second
case of the rigidity theorem deals with finite time of existence.
This case is dealt with in \cite{15} through the use of the $L^2$
conservation law, which is absent for the wave equation. We proceed
in two stages. First we show that the solution must have
self-similar behavior (Proposition \ref{p:5.14}). Then, in section 6,
following Merle-Zaag (\cite{22}) and earlier work on non-linear heat
equations by Giga-Kohn (\cite{8}), we introduce self-similar
variables and the new resulting equation, which has a monotonic
energy. We then show that there exists a non-trivial asymptotic
solution $w^*$, which solves a (degenerate) elliptic non-linear
equation. Finally, using the estimates we proved on $w^*$ and the
unique continuation principle, we show that $w^*$ must be zero, a
contradiction which gives our rigidity theorem. In section 7 we
prove our main theorem as a consequence of the rigidity theorem.

\mbr

Finally, we would like to point out that we expect that our
arguments will extend to $N\geq 6$, using arguments similar to those
in the work of Tao-Visan \cite{33}
for 
the local solvability in time of the equation and the corresponding 
extension of the work of Bahouri-Gerard \cite{4} (the rest of argument is independent of the dimension).

\section{A review of linear estimates and the Cauchy
problem}\label{s:2}

In this section we will review the theory of the Cauchy problem
 \bbe
\tag{CP}
\begin{cases}
\dd^2_t u - \De u = \abs{u}^{\frac{4}{N-2}} u & (x,t)\in\bR^N\times\bR
\\
u\bigl\lvert_{t=0} = u_0 \in \h(\bR^N) & \\
\dd_t u \bigl\lvert_{t=0} = u_1 \in L^2(\bR^N)
\end{cases}
 \ee
i.e. the $\h$ critical, focusing Cauchy problem for NLW, and some of
the associated linear theory. We start out with some preliminary
notation and linear estimates. Consider thus
 \bbe \tag{LCP}
\begin{cases}
\dd^2_t w - \De w = h & (x,t)\in\bR^N\times\bR
\\
w\bigl\lvert_{t=0} = w_0 \in \h(\bR^N) & \\
\dd_t w \bigl\lvert_{t=0} = w_1 \in L^2(\bR^N)
\end{cases}
 \ee
the associated linear problem. The solution operator to (LCP) is given by:
 \bbea
w(x,t)&=&\cos(t\sqrt{-\De})w_0+(-\De)^{1/2}\sin(t\sqrt{-\De})w_1\\
&+&\int_{0}^{t}\frac{\sin\left((t-s)\sd\right)}{\sd}h(s)ds\\
&=&S(t)(w_0,w_1)+\int_{0}^{t}\frac{\sin\left((t-s)\sd\right)}{\sd}h(s)ds.
 \eea

\begin{lem}[Strichartz estimates \cite{20}, \cite{10}]\label{l:2.1}
There is a constant $C$, independent of $T$, such that
 \bbea
&& \sup_{0<t<T}(\norm{w(t)}_{\h} + \norm{\dd_tw(t)}_{L^2}) +\norm{w}_{L_t^{\frac{2(N+1)}{(N-1)}}
\w{\frac{2(N+1)}{(N-1)}}}\\
&&
\norm{\dd_tw}_{L_t^{\frac{2(N+1)}{(N-1)}}W_x^{-\frac{1}{2},\frac{2(N+1)}{N-1}}}+\norm{w}_{L_t^{\frac{2(N+1)}{N-2}}L_x^{\frac{2(N+1)}{N-2}}}+\norm{w}_{L_t^{\frac{N+2}{N-2}}L_x^{\frac{2(N+2)}{N-2}}}\\
&& \leq C\left[\norm{w_0}_{\h(\bR^N)}+\norm{w_1}_{L^2(\bR^N)}+
\norm{h}_{L_t^{\frac{2(N+1)}{(N+3)}}\w{\frac{2(N+1)}{(N+3)}}} \right].
 \eea
\end{lem}

\begin{lem}[Trace Theorem]\label{l:2.2}
Let $w_0,w_1,h,w$ be as in Lemma \ref{l:2.1}. Then, for $\abs{d}\leq
1/4$,
 \bbea
&& \sup_t\norm{\na_xw\left(\frac{x_1-dt}{\sqrt{1-d^2}},x',\frac{t
-dx_1}{\sqrt{1-d^2}} \right)}_{L^2(dx_1dx')}\\
&& + \sup_t\norm{\dd_tw\left(\frac{x_1-dt}{\sqrt{1-d^2}},x',\frac{t
-dx_1}{\sqrt{1-d^2}}\right)}_{L^2(dx_1dx')}\\
&&\leq
C\left\{\norm{w_0}_{\h(\bR^N)}+\norm{w_1}_{L^2(\bR^N)}+\norm{h}_{L^1_tL^2_x}\right\}.
 \eea
\end{lem}
\begin{proof} Let $v(x,t)=U(t)f$ be given by
$\hat{v}(\xi,t)=e^{it\abs{\xi}}\hat{f}(\xi)$, with $f\in L^2$. We
will show that
$$
\sup_t\norm{v\left(\frac{x_1-dt}{\sqrt{1-d^2}},x',\frac{t
-dx_1}{\sqrt{1-d^2}} \right)}_{L^2(dx_1dx')}\leq C\norm{f}_{L^2},
$$
which easily implies the desired estimate. But
 \bbea
v(x,t)&=&\int e^{ix.\xi}e^{it\abs{\xi}}\hat{f}(\xi)d\xi=
\int e^{ix_1\xi_1}e^{it\abs{\xi}}e^{ix'.\xi'}\hat{f}(\xi)d\xi_1d\xi'\\
&=& \int
e^{ix_1\xi_1}e^{it\sqrt{\xi_1^2+\abs{\xi'}^2}}e^{ix'.\xi'}\hat{f}(\xi_1,\xi')d\xi_1d\xi',
 \eea
 so that
 \bbea
&& v\left(\frac{x_1-dt}{\sqrt{1-d^2}},x',\frac{t
-dx_1}{\sqrt{1-d^2}} \right)\\&&
=\int
e^{i(x_1-dt)\xi_1/\sqrt{1-d^2}}e^{i(t-dx_1)\sqrt{\xi_1^2+\abs{\xi'}^2}/\sqrt{1-d^2}}e^{ix'.\xi'}\hat{f}(\xi)d\xi_1d\xi'\\
&& =\int
e^{ix_1(\xi_1-d\abs{\xi})/\sqrt{1-d^2}}e^{ix'.\xi'}e^{-idt\xi_1/\sqrt{1-d^2}}e^{it\abs{\xi}/\sqrt{1-d^2}}\hat{f}(\xi)d\xi_1
d\xi'\\ &&=\int
e^{ix_1(\xi_1-d\abs{\xi})/\sqrt{1-d^2}}e^{ix'.\xi'}\hat{g}_t(\xi)d\xi_1d\xi',
 \eea
where
$\hat{g}_t(\xi)=e^{-idt\xi_1/\sqrt{1-d^2}}e^{it\abs{\xi}/\sqrt{1-d^2}}\hat{f}(\xi)$.
We now define $\eta_1=\frac{\xi_1-d\abs{\xi}}{\sqrt{1-d^2}}$,
$\eta'=\xi'$ and compute
 \bbea
\left\lvert \frac{d\eta}{d\xi} \right\rvert &=& \det \left(
\begin{array}{ccccc}
\frac{1-d\xi_1/\abs{\xi}}{\sqrt{1-d^2}} &
\frac{-d\xi_2/\abs{\xi}}{\sqrt{1-d^2}} & \ldots & \ldots &
\frac{-d\xi_N/\abs{\xi}}{\sqrt{1-d^2}} \\
0 & 1 & 0 & \ldots & 0 \\
0 & 0 & 1 & \ldots & 0 \\
\vdots & \vdots & \vdots & \ddots  & \vdots \\
0 & 0 & 0 & \ldots & 1
\end{array}
\right) \\
&=& \left(\frac{1-d\xi_1/\abs{\xi}}{\sqrt{1-d^2}}\right)\approx
1\,\text{ for }\,\abs{d}\leq 1/4.
 \eea
The result now follows from Plancherel's Theorem.
\end{proof}

\begin{rem}\label{r:2.3}
A density argument in fact shows that
$$
t\mapsto w\left(\frac{x_1-dt}{\sqrt{1-d^2}},x',\frac{t-dx_1}{\sqrt{1-d^2}} \right)\in
C\left(\bR;\h\left(\bR^N,dx_1,d\bar{x}'\right) \right),
$$
and similarly for $\dd_tw$.
\end{rem}
\begin{rem}\label{r:2.4}
Let $F(u)=\abs{u}^{\frac{4}{N-2}}u$, then clearly for
$3\leq N\leq 6,$
\[\abs{F(u)}\leq\abs{u}^{\frac{N+2}{N-2}}, \ \abs{(\na F)(u)}\leq
C\abs{u}^{\frac{4}{N-2}},\] \[\abs{(\na F)(u)-(\na F)(v)}\leq C\abs{u-v}
\left\{\abs{u}^{\frac{6-N}{N-2}}+\abs{v}^{\frac{6-N}{N-2}}\right\},\] 
 \bbea
&&
\abs{\na_x\left(F\left(u(x)\right)\right)-\na_x\left(F\left(v(x)\right)\right)}\leq
C\abs{u(x)}^{\frac{4}{N-2}} \abs{\na u(x)-\na v(x)}\\ && +C\abs{\na
v(x)}\left\{\abs{u}^{\frac{6-N}{N-2}}+\abs{v}^{\frac{6-N}{N-2}}\right\}
\abs{u-v}.
 \eea
We will need also a version of the chain rule for fractional
derivatives (see \cite{6}, \cite{16}, \cite{28}, \cite{34}).
\end{rem}

\begin{lem}\label{l:2.5}
Assume $F(0)=F'(0)=0$ and that for all a,b
$$\abs{F'(a+b)}\leq
C\left\{\abs{F'(a)}+\abs{F'(b)} \right\}, \
\abs{F''(a+b)}\leq C\left\{\abs{F''(a)}+\abs{F''(b)} \right\}.$$ We
then have, for $0<\al<1$ $$\norm{D^{\al}F(u)}_{L^p_x}\leq
C\norm{F'(u)}_{L_x^{p_1}}\norm{D^{\al}u}_{L^{p_2}_x},$$
where $\frac{1}{p}=\frac{1}{p_1}+\frac{1}{p_2}$, $1<p_i<\infty$ and
$$\norm{D^{\al}\left( F(u)-F(v)\right)}_{L_x^p}\leq
C \left\{\norm{F'(u)}_{L_x^{p_1}}+\norm{F'(v)}_{L_x^{p_1}}\right\} 
 \norm{D^{\al}(u-v)}_{L_x^{p_2}}$$  
$$+ C\left\{
\norm{F''(u)}_{L_x^{r_1}}+\norm{F''(v)}_{L_x^{r_1}}\right\} 
\left\{ \norm{
D^{\al}u}_{L_x^{r_2}}+\norm{D^{\al}v}_{L_x^{r_2}}\right\} 
\norm{u-v}_{L_x^{r_3}},$$
where $\frac{1}{p}=\frac{1}{r_1}+\frac{1}{r_2}+\frac{1}{r_3}$, $1< r_i <\infty$, $1< p
<\infty$.
\end{lem}

\begin{rem}\label{r:2.6}
In our application of Lemma \ref{l:2.5}, we will have
$$F(u)=\abs{u}^{\frac{4}{N-2}}u,\quad 3\leq N\leq 5,\text{ and }
F'(u)=C_N\abs{u}^{\frac{4}{N-2}},$$
$$F''(u)=\tilde{C}_N\sign(u)\abs{u}^{\frac{4}{N-2}-1}=
\tilde{C}_N\sign(u)\abs{u}^{\frac{6-N}{N-2}}.$$ We will choose $p=\frac{2(N+1)}{(N+3)},$
$p_2=\frac{2(N+1)}{(N-1)}$, so that $\frac{1}{p_1}=\frac{1}{p}-\frac{1}{p_2}=\frac{2}{N+1}$;
$r_3=\frac{2(N+1)}{N-2}$, $r_2=\frac{2(N+1)}{N-1}$, so that
$\frac{1}{r_1}=\frac{1}{p}-\frac{1}{r_2}-\frac{1}{r_3}= \frac{6-N}{2(N+1)}$. Notice
that $p_1\frac{4}{N-2}=\frac{2(N+1)}{N-2}$; $\frac{6-N}{N-2}
r_1=\frac{2(N+1)}{N-2}$. Let us now define the $S(I)$, $W(I)$ norm for an interval $I$
by 
$$\norm{v}_{S(I)}=\norm{v}_{L_I^{\frac{2(N+1)}{N-2}}L_x^{\frac{2(N+1)}{N-2}}}  \ \mbox{and} 
\ \norm{v}_{W(I)}=\norm{v}_{L_I^{\frac{2(N+1)}{N-1}}L_x^{\frac{2(N+1)}{N-1}}}.$$
\end{rem}
\begin{thm}[See \cite{23}, \cite{9}, \cite{25}]\label{t:2.7}
Assume $(u_0,u_1)\in\h\times L^2$, $0\in I$ an interval and
$\norm{(u_0,u_1)}_{\h\times L^2}\leq A$. Then, $($for $3\leq N\leq
5)$ there exists $\de=\de(A)$ such that if
$$\norm{S(t)\left((u_0,u_1)\right)}_{S(I)} <\de,$$ there exists a
unique solution $u$ to $(CP)$ in $\bR^N\times I$, with $(u,\dd_t
u)\in C(I;\h\times L^2)$, $\norm{D_x^{1/2}u}_{W(I)} + \norm{\dd_tD_x^{-1/2}u}_{W(I)}<+\infty$, $\norm{u}_{S(I)}\leq
2\de$. Moreover, if $(u_{0,k},u_{1,k})\to (u_0,u_1)$ as $k \to + \infty$ in $\h\times
L^2$ $($so that, for $k$ large $\norm{S(t)\left((u_0,u_1)
\right)}_{S(I)}$ $<\de)$, the corresponding solutions as $k \to + \infty$
$(u_k,\dd_t(u_k))\to (u,\dd_tu)$ in $C(I;\h\times L^2)$.
\end{thm}
\noindent {\em Sketch of the proof.} (CP) is equivalent to the
integral equation
$$
u(t)=S(t)\left((u_0,u_1) \right)+\int_{0}^t\frac{\sin\left((t-s)\sqrt{-\De} \right)}{\sqrt{-\De}}F(u)(s)ds,
$$
where $F(u)=\abs{u}^{\frac{4}{N-2}}u$. We let $$B_{a,b}=\left\{v\text{ on
}\bR^N\times I:\norm{v}_{S(I)}\leq a \text{ and
}\norm{D_x^{1/2}v}_{W(I)}\leq b\right\},$$
$$\Phi_{(u_0,u_1)}(v)=S(t)(u_0,u_1)+\int_0^t\frac{
\sin\left((t-s)\sqrt{-\De} \right)}{\sqrt{-\De}}F(v)(s)ds.$$ We will
next choose $\de,a,b$ so that $\Phi_{(u_0,u_1)}:B_{a,b}\to B_{a,b}$
and is a contraction there. Note that, by Lemma \ref{l:2.1},
$$\norm{D_x^{1/2}\Phi_{(u_0,u_1)}(v)}_{W(I)}\leq
CA+C\norm{F(v)}_{L_x^{\frac{2(N+1)}{(N+3)}}\w{\frac{2(N+1)}{(N+3)}}}.$$ But, by
Lemma \ref{l:2.5}, $\norm{D_x^{1/2}F(v)}_{L_x^{\frac{2(N+1)}{(N+3)}}}$ is bounded by
$$
C\norm{F'(v)}_{L_x^{\frac{N+1}{2}}}
\norm{D_x^{1/2}v}_{L_x^{\frac{2(N+1)}{N-1}}} 
\leq
C\norm{v}^{\frac{4}{N-2}}_{L_x^{\frac{2(N+1)}{N-2}}}\norm{D_x^{1/2}v}_{L_x^{\frac{2(N+1)}{(N-1)}}},$$
so that
 \bbea
&& \norm{D_x^{1/2}F(v)}_{L_I^{\frac{2(N+1)}{(N+3)}}L_x^{\frac{2(N+1)}{(N+3)}}}\leq C\norm{v}^{\frac{4}{N-2}}_{L_I^{\frac{2(N+1)}{N-2}}
L_x^{\frac{2(N+1)}{N-2}}}
\norm{D_x^{1/2}v}_{L_I^{\frac{2(N+1)}{(N-1)}}L_x^{\frac{2(N+1)}{(N-1)}}} \\
&&\le C\norm{v}^{\frac{4}{N-2}}_{S(I)}\norm{D_x^{1/2}v}_
{W(I)}.
 \eea
Hence, for $v\in B_{a,b}$,
$$
\norm{D_x^{1/2}\Phi_{(u_0,u_1)}(v)}_{W(I)}\leq CA+Ca^{\frac{4}{N-2}}
b.
$$
Similarly, using Lemma \ref{l:2.1} for the second term in $\Phi_{(u_0,u_1)}$, and the
argument above, together with our assumption on $(u_0,u_1)$ for the first term, we obtain:
$$
\norm{\Phi_{(u_0,u_1)}}_{S(I)}\leq\de+Ca^{\frac{4}{N-2}}b.
$$
Next, choose $b=2AC$, $a$ so that $Ca^{\frac{4}{N-2}}\leq \frac 12$. Then,
$$\norm{D^{1/2}_x\Phi_{(u_0,u_1)}(v)}_{W(I)}\leq b.$$ If $\de=a/2$ and
$Ca^{\left(\frac{4}{N-2}-1\right)}b\leq 1/2$ (possible if $N<6$) we
obtain $\norm{\Phi_{(u_0,u_1)}(v)}_{S(I)}\leq a$, so that
$\Phi_{(u_0,u_1)}:B_{a,b}\to B_{a,b}$. Next, for the contraction, we
again use Lemma \ref{l:2.1} and Lemma \ref{l:2.5}, to see that:
 \bbea
&& \norm{D_x^{1/2}\left(\Phi_{(u_0,u_1)}(v)-\Phi_{(u_0,u_1)}(v')
\right)}_{W(I)}+\norm{\Phi_{(u_0,u_1)}(v)-\Phi_{(u_0,u_1)}(v')}_{S(I)} \\
&& \leq C\norm{D_x^{1/2}\left(F(v)-F(v') \right)}_{L_I^{\frac{2(N+1)}{(N+3)}}L_x^{\frac{2(N+1)}{(N+3)}}} \\
&& \leq C\left[\left\{\norm{v}^{\frac{4}{N-2}}_{L_I^{\frac{2(N+1)}{N-2}}L_x^{\frac{2(N+1)}{N-2}}}+
\norm{v'}^{\frac{4}{N-2}}_{L_I^{\frac{2(N+1)}{N-2}}L_x^{\frac{2(N+1)}{N-2}}}  \right\}\right.\\
&&
\norm{D_x^{1/2}(v-v')}_{L_I^{\frac{2(N+1)}{(N-1)}}L_x^{\frac{2(N+1)}{(N-1)}}}+\left\{\norm{v}^{\frac{6-N}{N-2}}_
{L_I^{\frac{2(N+1)}{N-2}}L_x^{\frac{2(N+1)}{N-2}}}\right.\\ &&+\left.
\norm{v'}^{\frac{6-N}{N-2}}_
{L_I^{\frac{2(N+1)}{N-2}}L_x^{\frac{2(N+1)}{N-2}}}\right\}\left\{\norm{D_x^{1/2}v}_{L_I^{\frac{2(N+1)}{(N-1)}}L_x^{\frac{2(N+1)}{(N-1)}}}\right.\\
&&\left.\left.+\norm{D_x^{1/2}v'}_{L_I^{\frac{2(N+1)}{(N-1)}}L_x^{\frac{2(N+1)}{(N-1)}}}\right\}\norm{v-v'}_{L_I^{\frac{2(N+1)}{N-2}}
L_x^{\frac{2(N+1)}{N-2}}}\right]\\ &&\leq
2Ca^{\frac{4}{N-2}}\norm{D_x^{1/2}(v-v')}_{W(I)}+2Ca^{\frac{6-N}{N-2}}
2b\norm{v-v'}_{S(I)}
 \eea
and the contraction property follows for $N<6$. We then find $u\in B_{a,b}$ solving $\Phi_{(u_o,u_1)}(u)=u$.
To show that $(u,\dd_tu)\in C(I;\h\times L^2)$ we use Lemma \ref{l:2.1}, together with the fact that
$D_x^{1/2}F(u)\in L_I^{\frac{2(N+1)}{(N+3)}}L_x^{\frac{2(N+1)}{(N+3)}}$. This also shows that $\dd_tD_x^{-1/2}u\in W(I)$.
The continuity statement at the end is an easy consequence of the fixed point argument, so that the
proof is complete.
\begin{rem}\label{r:2.8}
$u\in L_I^{\frac{N+2}{N-2}}L_x^{2\frac{(N+2)}{N-2}}$, because of Lemma \ref{l:2.1} and the fact
that $D_x^{1/2}F(u)\in L_I^{\frac{2(N+1)}{(N+3)}}L_x^{\frac{2(N+1)}{(N+3)}}$. Note that because of
this and the integral equation, the conclusion of Lemma \ref{l:2.2} holds for $u$,
provided the integrations on the left hand side are restricted to $(x,x',t)\in\bR^N$
so that $\left(\frac{x_1-dt}{\sqrt{1-d^2}},x',\frac{t-dx_1}{\sqrt{1-d^2}}
\right)\in\bR^N\times I$.
\end{rem}

\begin{rem}[Higher regularity of solutions, see for example \cite{9}]\label{r:2.9}
If $(u_0,u_1)\in(\h\bigcap\hm^{1+\mu},H^{\mu})$, $0\leq\mu\leq 1$,
and $(u_0,u_1)$ verifies the conditions in Theorem \ref{t:2.7}, then
$(u,\dd_t u)\in C(I;\h\bigcap\hm^{1+\mu}\times H^{\mu})$ and
$$
\norm{D_x^{1/2+\mu}u}_{W(I)}+\norm{D_x^{1/2}u}_{W(I)} + 
\norm{\dd _t D_x^{\mu-1/2}u}_{W(I)} + \norm{\dd _t
D_x^{-1/2}u}_{W(I)}<\infty,$$
$\norm{u}_{S(I)}\leq 2\de.$
(In this
result we also need to use the assumption $3\leq N\leq 5$).
\end{rem}

\begin{rem}\label{r:2.10}
There exists $\tilde{\de}$ such that if $\norm{(u_0,u_1)}_{\h\times
L^2}\leq \tilde{\de}$, the conclusion of Theorem \ref{t:2.7} applies
to any interval $I$. In fact, by Lemma \ref{l:2.1},
$\norm{S(t)(u_0,u_1)}_{S((-\infty,+\infty))}\leq C\tilde{\de}$ and the
claim follows.
\end{rem}

\begin{rem}\label{r:2.11}
Given $(u_0,u_1)\in \h\times L^2$, there exists $(0\in) I$ such that
the hypothesis of Theorem \ref{t:2.7} is verified on $I$. This is
clear because, by Lemma \ref{l:2.1},
$\norm{S(t)(u_0,u_1)}_{S(I)}<+\infty$.
\end{rem}

\begin{rem}[Finite speed of propagation, see for instance \cite{26}]\label{r:2.12}
Let $R$ denote the fundamental solution of the Cauchy problem, i.e.
$R$ solves
 \bbe
\begin{cases}
(\dd^2_t - \De_x) u = 0 & (x,t)\in\bR^N\times\bR
\\
u\bigl\lvert_{t=0} = 0 & \\
\dd_t u \bigl\lvert_{t=0} = \de(x),
\end{cases}
 \ee
where $\de(x)$ is the Dirac mass at $0$. Then, we can write the
solution  of (LCP) in the form
$$
w(t)=\dd_tR(t)*w_0+R(t)*w_1-\int_o^tR(t-s)*h(s)ds,
$$
where $*$ denotes convolution in the spatial variable. As is well
known, $\supp R(-,t)\subset\overline{B(0,t)}$ and $\supp
\dd_tR(-,t)\subset\overline{B(0,t)}$. Thus, if \\
$\supp
u_0\subset{}^c\overline{B(x_0,a)},\,\, \supp
u_1\subset{}^c\overline{B(x_0,a)}$ and \\
$\supp
h\subset{}^c\left(\bigcup_{0\leq t\leq
a}\left[\overline{B(x_0,a-t)}\times(a-t)\right]\right), $ we have
$$w\equiv 0 \text{ on }\bigcup_{0\leq t\leq
a}\left[B(x_0,a-t)\times(a-t)\right].$$
\end{rem}

These remarks have immediate consequences for the solutions of (CP) given in Theorem
\ref{t:2.7}. In fact, suppose that $(u_0,u_1)$, $(u'_0,u'_1)$ are data verifying the
conditions of Theorem \ref{t:2.7} and such that $(u_0,u_1)=(u'_0,u'_1)$ in
$\bar{B}(x_0,a)$. Then, the corresponding solutions $u,u'$ agree on $\bigcup_{0\leq
t\leq a}\left[\bar{B}(x_0,(a-t))\times(a-t)\right]\bigcap\left\{\bR^N\times I
\right\}$. To see this, for $n \in \bN$, define $u^{(n+1)}(x,t)=S(t)\left((u_0,u_1)
\right)+\int_0^t\frac{\sin\left(
(t-s)\sqrt{-\De}\right)}{\sqrt{-\De}}F\left(u^{(n)}\right)ds$ (for $n=0$, we set
$u^{(0)}(x,t)=S(t)\left((u_0,u_1) \right)$). We define correspondingly
$u'^{(n+1)}(x,t)$. The proof of Theorem \ref{t:2.7} gives us $u=\lim_n u^{(n)}$ and
$u'=\lim_nu'^{(n)}$. The previous remarks allow us to show inductively that
$u^{(n+1)}=u'^{(n+1)}$ on $\bigcup_{0\leq t\leq
a}\left[B(x_0,(a-t))\times(a-t)\right]\bigcap\left\{\bR^N\times I \right\}$, which
establishes the claim. Typical applications of this remark are the following: \newline
\noindent a) If $\supp (u_0)\subset B(0,b)$, $\supp (u_1)\subset
 B(0,b)$ and $(u_0,u_1)$ verifies the hypothesis of Theorem
 \ref{t:2.7}, then $$u(x,t)\equiv 0\text{ on }\left\{(x,t):\abs{x}>b+t,t\geq 0,t\in I
 \right\}.$$ \newline\noindent b) We can approximate solutions $u$
 in $\bR\times I'$, $I'\subset\subset I$ by means of regular,
 compactly supported solutions, combining a), Remark \ref{r:2.9}
 and the last statement in Theorem \ref{t:2.7}.

Similar statements hold for $t<0$, for instance if
$(u_0,u_1)=(u'_0,u'_1)$ in $\bar{B}(x_0,a)$ then $u,u'$ agree on
$\bigcup_{-a\leq t\leq
0}\left[B(x_0,(a+t))\times(a+t)\right]\bigcap\bR^N\times I $.

\begin{defin}\label{d:2.13}
Let $t_0\in I$. We say that $u$ is a solution of $(CP)$ in $I$ if
$(u,\dd_tu)\in C(I;\h(\bR^N)\times L^2)$, $D_x^{1/2}u\in W(I)$,
$u\in S(I)$, $(u,\dd_tu)\big\vert_{t=t_0}=(u_0,u_1)$ and the
integral equation
$$
u(t)=S(t)\left((u_0,u_1) \right)+\int_{t_0}^{t}\frac{\sin\left(
(t-s)\sqrt{-\De}\right)}{\sqrt{-\De}}F\left(u(s)\right)ds
$$
holds, with $F(u)=\abs{u}^{\frac{4}{N-2}} u$, for $x\in \bR^N$, $t\in
I$.
\end{defin}

Note that if $u^{(1)}$, $u^{(2)}$ are solutions of (CP) on $I$, and
\[\left(u^{(1)}(t_0),\dd_tu^{(1)}(t_0)\right)=\left(u^{(2)}(t_0),\dd_tu^{(2)}(t_0)\right),\]
then $u^{(1)}\equiv u^{(2)}$ on $\bR^N\times I$. (See the argument
in \cite{15}, Definition 2.10). This allows one to define a maximal
interval $$I\left((u_0,u_1)
\right)=\left(t_0-T_-((u_0,u_1)),t_0+T_+((u_0,u_1)) \right),$$ with
$T_{\pm}\left( (u_0,u_1)\right)>0$ where the solution is defined. 
If  $T_1>t_0-T_-((u_0,u_1))$ and 
$T_2<t_0+T_+(u_0)$, $t_0\in (T_1,T_2)$,
then $u$ solves (CP) in $\bR^N\times [T_1,T_2]$, so that
\[
(u,\dd_tu)\in C([T_1,T_2];\h\times L^2)), \quad D_x^{1/2}u\in
W([T_1,T_2]), \quad u\in S([T_1,T_2]),
\]
\[
u\in
L^{\frac{(N+2)}{N-2}}\left([T_1,T_2];L_x^{2\frac{(N+2)}{N-2}}\right), \quad
\dd_tD_x^{-1/2}u\in W([T_1,T_2]).
\]

\begin{rem}\label{r:2.14}
If $u$ is such that $(u,\dd_tu)\in C(I;\h\times L^2)$, $\norm{u}_{S(I)}\leq B$ and
there exist $u_j$ with $(u_j,\dd_t(u_j))\in C(I;\h\times L^2)$, $(u_j,\dd_t (u_j))\to
(u,\dd_tu)$ in $C(I;\h\times L^2)$, with $u_j$ a solution of (CP) in $I$ together with
$\norm{u_j}_{S(I)}\leq B$, then $\norm{D_x^{1/2}u}_{W(I)}<+\infty$ and $u$ is a
solution of (CP) in $I$. This follows by showing that $\norm{D_x^{1/2}u_j}_{W(I)}\leq
B'$, where $B'$ is independent of $j$. To show this, first find $A$ so that $\sup_{t\in
I}\norm{(u_j,\dd_t(u_j))}_{\h\times L^2}\leq A$, for all $j$. Next, partition
$I=\bigcup_{k=1}^M I_k$, where $I_k$ is such that $\norm{u_j}_{S(I_k)}\leq\de$, where
$\de=\de(A)$ is to be chosen. Note that $M=M(B,\de)$. We then use the integral
equation for $u_j$, and the estimate

$$
\norm{D_x^{1/2}F(u_j)}_{L_{I_k}^{\frac{2(N+1)}{(N+3)}}L_{x}^{\frac{2(N+1)}{(N+3)}}}\leq
C\de^{\frac{4}{N-2}}\norm{D_x^{1/2}u_j}_{W(I_k)}
$$
(see the proof of Theorem \ref{t:2.7}), so that
$$
\norm{D_x^{1/2}u_j}_{W(I_k)}\leq
CA+C\de^{\frac{4}{N-2}}\norm{D_x^{1/2}u_j}_{W(I_k)}.
$$
Thus, for $\de$ small we obtain $\norm{D_x^{1/2}u_j}_{W(I_k)}\leq
2CA$ and adding in $k$ we obtain the desired bound.
\end{rem}

\begin{lem}[Standard finite blow-up criterion]\label{l:2.15}
If $T_+\left((u_0,u_1) \right)<\infty$, then
$$\norm{u}_{S([t_0,t_0+T_+(u_0,u_1)])}=+\infty.$$ A corresponding
result holds for $T_-\left( (u_0,u_1) \right)$.
\end{lem}

The proof is similar to the one in Lemma 2.11 of \cite{15}.

\begin{rem}[Energy and moment identities]\label{r:2.16}
Let $(u_0,u_1)\in \h\times L^2$ and let $0\in I$ be the maximal
interval of existence. Then, for $t\in I$, with
$\frac{1}{2^*}=\frac{1}{2}-\frac{1}{N}$ ($2^*= \frac{2N}{N-2}$), we have
 \bbea
E
(u(t),\dd_tu(t))&=&\int_{\bR^N}\frac{1}{2}\abs{\dd_tu(x,t)}^2+\frac{1}{2}\abs{\na_xu(x,t)}^2
-\frac{1}{2^*}\abs{u(x,t)}^{2^*}dx\\&=&E\left( (u_0,u_1)\right),
 \eea
 and
 \bbe\label{eq:4.15}
\int\na_xu(x,t)\dd_tu(x,t)dx=\int\na u_0u_1.
 \ee
\end{rem}
\begin{proof}
Let
$e(u)(x,t)=\frac{1}{2}(\dd_tu)^2(x,t)+\frac{1}{2}\abs{\na_xu(x,t)}^2-\frac{1}{2^*}\abs{u(x,t)}^{2^*}.$
Then, for sufficiently smooth solutions of (CP) we have
 \bbe\label{eq:2.17}
\dd_te(u)(x,t)=\sum_{j=1}^N \dd_{x_j}\left(\dd_{x_j}
u(x,t)\dd_tu(x,t) \right),
 \ee
as is readily seen. Now, fix any $I'\subset\subset I$, so that
$\norm{u}_{S(I')}<+\infty$. By dividing $I'=\bigcup_{k=1}^MI_k$,
with $\norm{u}_{S(I_k)}\leq\de(A)$, where $$A=\sup_{t\in
I'}\norm{(u(t),\dd_tu(t))}_{\h\times L^2},$$ we can use Theorem
\ref{t:2.7} to approximate $u$ by compactly supported solutions in
$\bR^N\times I_k$ (see Remarks \ref{r:2.9}, \ref{r:2.12}). We then
apply \eqref{eq:2.17} and integrate by parts, and then pass to the
limit, for $t\in I_k$. The proof of second equality is similar.
\end{proof}

\begin{lem}\label{l:2.18}
Let $(u_0,u_1)\in\h\times L^2$, $\norm{(u_0,u_1)}_{\h\times L^2}\leq
A$ with maximal interval of existence
$I=(-T_-(u_0,u_1),T_+(u_0,u_1))$. There exists
$\epsilon_0>0$ so that, if for some $M>0$ and  $0<\epsilon<\epsilon_0$, we have
$\int_{\abs{x}\geq
M}\abs{\na_xu_0}^2+\frac{\abs{u_0}^2}{\abs{x}^2}+\abs{u_1}^2\leq\epsilon$,
then for $t\in I_+=[0,\infty)\bigcap I$, we have
$$
\int_{\abs{x}\geq \frac32 M + t}\abs{\na_xu(x,t)}^2+\abs{\dd_tu(x,t)}^2dx\leq
C\epsilon.
$$
\end{lem}
\begin{proof}
Choose $\Psi_M\equiv 1$ for $\abs{x}\geq \frac{3}{2}M$,
$\Psi_M\equiv 0$ for $\abs{x}\leq M$, $\abs{\na_x\Psi_M}\leq C/M$.
Define $u_{0,M}=\Psi_{M}u_0$, $u_{1,M}=\Psi_{M}u_1$. Because of our
assumption, we have $\norm{(u_{0,M},u_{1,M})}_{\h\times L^2}\leq
C\epsilon$. Choose now $\epsilon_0$ so small that
$C\eps_0\leq\tilde{\de}$, where $\tilde{\de}$ is as in Remark
\ref{r:2.10}. Then, there exists $u_M$ solving the (CP) in
$I=(-\infty,+\infty)$, with $(u_M(0),\dd_tu_M(0))=(u_{0,M},u_{1,M})$
and such that $\sup_{t\in
(-\infty,+\infty)}\norm{(u_M(t),\dd_tu_M(t))}_{\h\times L^2}\leq
2C\eps$. But, by Remark \ref{r:2.12}, $u_M(x,t)=u(x,t)$ for
$\abs{x}\geq \frac{3M}{2}+t$, $t\in I_+$. The Lemma follows.
\end{proof}

\begin{defin}\label{d:2.19}
Let $(v_0,v_1)\in\h\times L^2$, $v(x,t)=S(t)((v_0,v_1))$ and let $\{
t_n\}$ be a sequence, with
$\lim_{n\to\infty}t_n=\bar{t}\in[-\infty,+\infty]$. We say that
$u(x,t)$ is a non-linear profile associated with $((v_0,v_1),\{t_n
\})$ if there exists an interval $I$, with $\bar{t}\in
\overset{\circ}I$ (if $\bar{t}=\pm \infty, I=[a,+\infty)$ or
$(-\infty,a]$) such that $u$ is a solution of (CP) in $I$ and
$$
\lim_{n\to
\infty}\norm{\left(u(t_n)-v(t_n),\dd_tu(t_n)-\dd_tv(t_n)\right)}_{\h\times
L^2}=0.
$$
\end{defin}

\begin{rem}\label{r:2.20}
There always exists a non-linear profile associated to $((v_0,v_1),\{t_n\})$. The
proof is similar to the one in \cite{15}, Remark 2.13, once we use the proof of
Theorem \ref{t:2.7} and the linear estimates, (with
$w(x,t)=\int_t^{\infty}\frac{\sin\left((t-s)\sqrt{-\De} \right)}{\sqrt{-\De}}h(s)ds$,
$I=(a,+\infty)$, $a>0$)
\bbea
\sup_{t\in I}\norm{(w(t),\dd_tw(t))}_{\h\times
L^2}+&&\norm{D_x^{1/2}w}_{W(I)}+\norm{w}_{S(I)} \\ \leq &&C\norm{h}_{L_I^{\frac{2(N+1)}{(N+3)}}\dot{W}_x^{1/2,\frac{2(N+1)}{(N+3)}}},
\eea
 which follow from \cite{10}, Proposition 3.1, (2) and (3). Also, as
in \cite{15}, Remark 2.13, we have uniqueness of the non-linear
profile and a maximal interval of existence of the non-linear
profile associated to $((v_0,v_1),\{t_n\})$.
\end{rem}

\begin{thm}[Long time perturbation theory, see also \cite{33}, \cite{15}]\label{t:2.21}
Let $I\subset \bR$ be a time interval. Let $t_0\in I$, $(u_0,u_1)\in\h\times
L^2$ and  some constants $M,A,A'>0$. Let $\tilde{u}$ be defined
on $\bR^N\times I$ $(3\leq N\leq 5)$ and satisfy $\sup_{t\in
I}\norm{(\tilde{u}(t),\dd_t\tilde{u}(t))}_{\h\times L^2}\leq A$,
$\norm{\tilde{u}(t)}_{S(I)}\leq M$ and $\norm{D_x^{1/2}\tilde{u(t)}}_{W(I')}<\infty$ for
each $I'\subset\subset I$. Assume that
$$
(\dd_t^2-\De_x)(\tilde{u})-F(\tilde{u})=e,\quad (x,t)\in\bR^N\times
I
$$
$($in the sense of the appropriate integral equation$)$ and that
$$
\norm{(u_0-\tilde{u}(t_0),u_1-\dd_t\tilde{u}(t_0))}_{\h\times
L^2}\leq A',
$$
$$
\norm{D_x^{1/2}e}_{L_I^{\frac{2(N+1)}{(N+3)}}L_x^{\frac{2(N+1)}{(N+3)}}}+\norm{S(t-t_0)\left((u_0-\tilde{u}(t_0),u_1-\dd_t\tilde{u}(t_0))
\right)}_{S(I)}\leq\eps,
$$Then there exists $\eps_0=\eps_0(M,A,A')$ such that there exists a
solution of $(CP)$ in $I$ with $(u(t_0),\dd_tu(t_0))=(u_0,u_1)$, for
$0<\eps<\eps_0$, with $\norm{u}_{S(I)}\leq C(M,A,A')$ and $\forall
t\in I$,
$$
\norm{(u(t),\dd_tu(t))-(\tilde{u}(t),\dd_t\tilde{u}(t))}_{\h\times
L^2}\leq C(A,A',M)(A'+\eps).
$$
\end{thm}
The proof is analogous to the one given in \cite{15}, Theorem 2.14, using the ideas in
the proof of Theorem \ref{t:2.7}.

\begin{rem}\label{r:2.22}
Theorem \ref{t:2.21} yields the following continuity fact, which
will be used later. Let $(\tilde{u}_0,\tilde{u}_1)\in\h\times L^2$,
$\norm{\rtt}_{\hl}\leq A$, let $\tilde{u}$ be the solution of (CP),
with maximal interval of existence $$(-T_-(\rtt),T_+(\rtt)).$$ Let
$\left(u_0^{(n)},u_1^{(n)}\right)\to\rtt$ in $\hl$ and let $u^{(n)}$
be the corresponding solution of (CP), with maximal interval of
existence
$$(-T_-((u_0^{(n)},u_1^{(n)})),T_+((u_0^{(n)},u_1^{(n)}))).$$ Then
$$T_-(\rtt)\leq\underline{\lim}\,T_-((u_0^{(n)},u_1^{(n)})),\quad
T_+(\rtt)\leq\underline{\lim}\,T_+((u_0^{(n)},u_1^{(n)}))$$ and for
each $t\in(-T_-(\rtt),T_+(\rtt))$ we have
$$(u^{(n)}(t),\dd_tu^{(n)}(t))\to (\tilde{u}(t),\dd_t\tilde{u}(t))
\text{ in } \hl.$$ Indeed, let $I\ssub (-T_-(\rtt),T_+(\rtt))$, so that
$$\sup_{t\in
I}\norm{(\tilde{u}(t),\dd_t\tilde{u}(t))}_{\hl}\leq\tilde{A},\quad
\norm{\tilde{u}}_{S(I)}\leq M<+\infty.$$ We will show that, for $n$
large, $u^{(n)}$ exists on $I$, and that
 \bbea
&& \sup_{t\in
I}\norm{(u^{(n)}(t),\dd_tu^{(n)}(t))-(\tilde{u}(t),\dd_t\tilde{u}(t))}_{\hl}\\
&& \leq C(M,\tilde{A})\norm{(u_0^{(n)},u_1^{(n)})-\rtt}_{\hl},
 \eea
and additionally,
$\norm{u^{(n)}}_{S(I)}\leq\tilde{M}(\tilde{A},M)$. To show this,
apply Theorem  \ref{t:2.21}, with $u=u^{(n)}$,
$\rt=(u_0^{(n)},u_1^{(n)})$, $e\equiv 0$. If
$\eps_0=\eps_0(M,\tilde{A},2\tilde{A})$ and $n$ is large enough
that
$$\norm{S(t)((\tilde{u}_0-u_0^{(n)},\tilde{u}_1-u_1^{(n)}))}_{S(I)}\leq\eps, \norm{(\tilde{u}_0-\tilde{u}_0^{(n)},\tilde{u}_1-\tilde{u}_1^{(n)})}_{\hl}\leq2\tilde{A},$$ the
desired conclusions follow from Theorem \ref{t:2.21}. Note also that if we choose
$u_0^{(n)}$, $u_1^{(n)}$ in $C^{\infty}_0(\bR^N)$, the approximating solutions
$u^{(n)}$ will be regular in view of Remark \ref{r:2.9} and for $t\in I$ will have
compact support in $x$, in view of Remark \ref{r:2.12}, and will verify
$\norm{u^{(n)}}_{S(I)}\leq \tilde{M}$.
\end{rem}

\begin{rem}\label{r:2.23}
If $u$ is a solution of (CP) in $\bR^N\times I$, for each $I'\ssub
I$, $I=[a,+\infty)$ (or $I=(-\infty,a]$), such that
$\norm{u}_{S(I)}<\infty$, there exists $(u^+_0,u^+_1)\in\hl$ such
that
$$
\lim_{t\uparrow+\infty}\norm{[(u(t),\dd_tu(t))-(S(t)(u^+_0,u^+_1),
\dd_tS(t)(u^+_0,u^+_1))]}_{\hl}=0.
$$
See Remark 2.15 in \cite{15}, \cite{4} for a similar proof, based on
the fact, in our case, that
$\norm{D_x^{1/2}F(u)}_{L_I^{\frac{2(N+1)}{(N+3)}}L_x^{\frac{2(N+1)}{(N+3)}}}<\infty$
and the inequality used in the proof of Remark \ref{r:2.20}.
\end{rem}
\begin{rem}\label{r:2.24}
We recall that, since we are working in the focusing case, from the work of Levine
(\cite{19}, \cite{29}) we have that if $\rt\in
H^1\times L^2
$ is such that $E(\rt)<0$, then the
maximal interval of existence is finite. 
We will return to the issue of break-down in finite time (blow-up), in
the next section and at the end of the paper.

\end{rem}

\section{Variational estimates}\label{s:3}
Let $W(x) = W(x,t) = \frac{1}{\left(
1+\frac{\abs{x}^2}{N(N-2)}\right)^{\frac{(N-2)}{2}}}$
be a
stationary solution of (CP). That is $W$ solves the non-linear
elliptic equation
 \bbe\label{eq:3.1}
\De W+\abs{W}^{\frac{4}{N-2}}W=0.
 \ee
Moreover, $W\geq 0$ and it is radially symmetric and decreasing.
Note that $W\in\h$, but $W$ need not belong to $L^2$, depending on the dimension. By invariances
of the equation \eqref{eq:3.1}, for $\theta_0 \in[-\pi,\pi]$,
$\la_0>0$, $x_0\in\bR^N$,
$W_{\theta_0,x_0,\la_0}(x)=e^{i\theta_0}\la_0^{\frac{(N-2)}{2}}W(\la_0(x-x_0))$
is still a solution of \eqref{eq:3.1}. By the work of Aubin
\cite{3}, Talenti \cite{31} we have the following characterization
of $W$:
 \bbe\label{eq:3.2}
\forall u\in \h, \quad\norm{u}_{L^ {2^*}}\leq C_N\norm{\na u}_{L^2};
 \ee
moreover,
 \bbe\label{eq:3.3}
\text{If }\norm{u}_{L^{2^*}}=C_N\norm{\na u}_{L^2},\, u\ne 0,\text{
then }\exists(\theta_0,\la_0,x_0):u=W_{\theta_0,x_0,\la_0},
 \ee
where $C_N$ is the best constant of the Sobolev inequality
\eqref{eq:3.2} in dimension $N$.

\mbr

Remark that $\int\abs{\na W}^2=\frac1{C_N^N}$ and $\cE(W)=\frac{1}{N}\frac{1}{C_N^N}$, where $$\cE(u)=\int\frac{1}{2}\abs{\na u}^2-\frac{1}{2^*}\abs{u}^{2^*}.$$
Indeed, the equation \eqref{eq:3.1} gives $\int\abs{\na
W}^2=\int\abs{W}^{2^{*}}$. Also, \eqref{eq:3.3} yields
$C_N^2\int\abs{\na W}^2=\left(\int\abs{W}^{2^*}
\right)^{\frac{N-2}N},$
 so that $C_N^2\int\abs{\na
W}^2=\left(\int\abs{\na W}^{2} \right)^{\frac{(N-2)}{2}}$. Hence,
$\int\abs{\na W}^2=\frac1{C_N^N}$ and  $\cE(W)=(\frac{1}{2}-\frac{1}{2^*})\int\abs{\na
W}^2=\frac{1}{NC_N^N}$.

\begin{lem}\label{l:3.4}
Let $u\in\h(\bR^N)$ be such that for $\de_0>0$,
$$
\norm{\na u}^2_{L^2}<\norm{\na W}^2_{L^2} \ \mbox{and} \ \cE(u)\leq (1-\de_0)\cE(W).
$$
Then there exists $\bar{\de}=\bar{\de}(\de_0)>0$ such that
$$
\norm{\na u}^2_{L^2}\leq (1-\bar{\de})\norm{\na W}^2_{L^2}
\ \mbox{and} \
\cE(u)\geq 0.
$$
\end{lem}
\begin{proof}

It is contained in Lemma 3.4 of \cite{15}.
\end{proof}
\begin{cor}\label{c:3.5}
If $u$ is as in Lemma \ref{l:3.4}, then there exists
$C_{\bar{\de}}>0$ so that
$$
\int\abs{\na u}^2-\abs{u}^{2^*}\geq C_{\bar{\de}}\int\abs{\na u}^2.
$$
\end{cor}
\begin{proof}
Note that \eqref{eq:3.2} implies that

 \bbea
\int\abs{\na u}^2-\abs{u}^{2^*}&&\geq \int\abs{\na
u}^2-C_N^{2^*}\left(\int\abs{u}^2\right)^{\frac{2^*}2}\\
&& \geq\int\abs{\na u}^2
\left[1-C_N^{2^*}\left(\int\abs{\na u}^2\right)^{\frac2{N-2}}\right]\\
&&\geq \int\abs{\na u}^2
\left[1-C_N^{2^*}(1-\bar{\de})^{\frac1{N-2}}\left(\int\abs{\na
W}^2\right)^{\frac2{N-2}}\right],
 \eea
by Lemma \ref{l:3.4}. But
$
\left(\int\abs{\na W}^2\right)^{\frac2{N-2}}=
\frac{1}{C_N^{\frac{2N}{N-2}}}=\frac{1}{C_N^{2^*}},
$
so that the corollary follows with
$C_{\bar{\de}}=[1-(1-\bar{\de})^{\frac1{N-2}}].$
\end{proof}
\begin{cor}\label{c:3.6}
Let $u\in\h$, $\norm{\na u}_{L^2}<\norm{\na W}_{L^2}.$ Then
$\cE(u)\geq 0$.
\end{cor}
\begin{proof}
If $\cE(u)<\cE(W)=\frac{1}{N}\frac{1}{C_N^N}$, the claim follows
from Lemma \ref{l:3.4}. If, on the other hand $\cE(u)\geq\cE(W)$,
the statement is obvious.
\end{proof}
\begin{rem}\label{r:3.7}
Let $u\in\h(\bR^N)$ be such that $\cE(u)\leq (1-\de_0)\cE(W)$.
Assume that $\norm{\na u}^2_{L^2}>\norm{\na W}^2_{L^2}.$ Then there
exists $\bar{\de}=\bar{\de}(\de_0,N)$ such that
$$
\norm{\na u}^2_{L^2}\geq(1+\bar{\de})\norm{\na W}^2_{L^2}.
$$
The proof of this is similar to the one of Lemma \ref{l:3.4}. See
Remark 3.14 in \cite{15}.
\end{rem}
\begin{thm}[Energy trapping]\label{t:3.8}
Let $u$ be a solution of $(CP)$, with
$(u,\dd_tu)\big\vert_{t=0}=\rt\in\hl$ and maximal interval of
existence $I$. Assume that, for $\de_0>0$,
$$E(\rt)\leq(1-\de_0)E((W,0))  \ \mbox{and} \ \norm{\na u_0}^2_{L^2}<\norm{\na
W}^2_{L^2}.$$ Then, there exists $\bar{\de}=\bar{\de}(\de_0)$ such
that, for $t\in I$, we have
 \bbe\label{eq:3.9}
\norm{\na_xu(t)}^2_{L^2}\leq (1-\bar{\de})\norm{\na W}^2_{L^2}
 \ee
 \bbe\label{eq:3.10}
\int\abs{\na _xu(t)}^2-\abs{u(t)}^{2^*}\geq
C_{\bar{\de}}\int\abs{\na_xu(t)}^2
 \ee
 \bbe\label{eq:3.11}
\cE(u(t))\geq 0\quad (\text{and hence }E((u(t),\dd_tu(t)))\geq 0).
 \ee
\end{thm}
\begin{proof}
By Remark \ref{r:2.16}, $E((u(t),\dd_tu(t)))=E(\rt),$ $t\in I$.
Also, $\cE(u(t))\leq E((u(t),\dd_tu(t)))$. Thus, the Theorem follows
from Lemma \ref{l:3.4}, Corollary \ref{c:3.5}, Corollary \ref{c:3.6}
and a continuity argument.
\end{proof}
\begin{cor}\label{c:3.9}
Let $u$ be as in Theorem \ref{t:3.8}. Then for all $t\in I$ we have
$E((u(t),\dd_tu(t)))\simeq\norm{(u(t),\dd_tu(t))}^2_{\hl}\simeq\norm{\rt}^2_{\hl}$,
with comparability constants which depend only on $\de_0$.
\end{cor}

\begin{proof}
For $t\in I$,
$E((u(t),\dd_tu(t)))\leq\norm{(u(t),\dd_tu(t))}^2_{\hl}$. Also,
 \bbea
&&
E((u(t),\dd_tu(t)))=\frac{1}{2}\int\left(\dd_tu(t)\right)^2+\cE(u(t))\\
&&
=\frac{1}{2}\int\left(\dd_tu(t)\right)^2+\frac{1}{2}\left[\int\abs{\na_x
u(t)}^2-\abs{u(t)}^{2^*} \right]
+\left(\frac{1}{2}-\frac{1}{2^*}\right)\int\abs{u(t)}^{2^*}\\
&&\geq\frac{1}{2}\int(\dd_tu(t))^2
+C_{\bar{\de}}\int\abs{\na_xu(t)}^2.
 \eea
Finally,
$
E((u(t),\dd_tu(t)))=E(\rt)\simeq \norm{\rt}^2_{\hl}.
$
\end{proof}
\begin{thm}[Finite time blow-up, see also Remark \ref{r:2.24}]\label{t:3.10}
Assume that $\rt\in\hl$, $u_0\in L^2$ and that $u$ is the solution
of $(CP)$ with maximal interval of existence $I$. Assume that
$E(\rt)<E((W,0))$ and $\int\abs{\na u_0}^2>\int\abs{\na W}^2.$ Then
$I$ must be a finite interval.
\end{thm}

\begin{proof}
Fix $\de_0$ positive so that $E(\rt)\leq (1-\de_0)E((W,0))$. Define
$$y(t)=\int\abs{u(x,t)}^2dx.$$ We then have $$y'(t)=2 \int u\dd_tu \ \mbox{and} \
y''(t)=2\left[\int(\dd_tu)^2-\abs{\na_xu}^2+\abs{u}^{2^*} \right].$$
(To check these identities, we proceed as in Remark \ref{r:2.16},
starting with data in $C_0^{\infty}$ and using a  limiting
argument.) Let $\tilde{\de}_0=\de_0E((W,0))$, so that $E((W,0))\geq
E((u(t),\dd_tu(t)))+\tilde{\de}_0$ and hence
$\frac{1}{2^*}\int\abs{u(t)}^{2^*}\geq\frac{1}{2}\int\left((\dd_tu(t))^2+\abs{\na_xu(t)}^2
\right)-E((W,0))+\tilde{\de}_0$ so that
$$
\int\abs{u(t)}^{2^*}\geq
\frac{N}{N-2}\int\left((\dd_tu(t))^2+\abs{\na_xu(t)}^2\right)-2^*E((W,0))+2^*\tilde{\de}_0.
$$
But then, (with $\tilde{\tilde{\de}}_0=22^*\tilde{\de}_0$) we have
 \bbea
&&
y''(t)\geq2\int(\dd_tu(t))^2+\frac{2N}{N-2}\int(\dd_tu(t))^2-2
2^*E((W,0))\\
&&
+\frac{2N}{N-2}\int\abs{\na_xu(t)}^2-2\int\abs{\na_xu(t)}^2+\tilde{\tilde{\de}}_0\\
&&=4\frac{(N-1)}{N-2}\int(\dd_tu(t))^2+\frac{4}{N-2}\int\abs{\na_xu}^2-\frac{4}{N-2}\int\abs{\na
W}^2+\tilde{\tilde{\de}}_0\\ && \geq 4\frac{(N-1)}{N-2}
\int(\dd_tu(t))^2+\tilde{\tilde{\de}}_0
 \eea
(by Remark \ref{r:3.7} and a continuity argument.) \\
Assume now that
$I\bigcap[0,\infty)=[0,\infty)$. Then, by our lower bound on
$y''(t)$, there exists $t_0>0$ such that $y'(t_0)>0$, and hence
$y'(t)> 0$ for $t>t_0$. Hence, for $t>t_0$,
$$
y''(t)y(t)\geq4\left[\frac{N-1}{N-2}\right]\left(\int(\dd_tu)^2(t)\right)\left(\int
u(t)^2\right)\geq \frac{N-1}{N-2} y'(t)^2,
$$
so that, for $t>t_0$,
$$\frac{y''(t)}{y'(t)}\geq\frac{N-1}{N-2}\frac{y'(t)}{y(t)}\quad\text{
or } \quad\left(\log y'(t)\right)'\geq \frac{N-1}{N-2}\left(\log
y(t)\right)'.$$ Hence for $t>t_0$, $$\log y'\geq \frac{N-1}{N-2}\log
y-C_0\quad\text{ or }\quad y'(t)\geq \tilde{C_0}y^{\frac{N-1}{N-2}},$$
which leads to finite time blow-up of $y$, because
$\frac{N-1}{N-2}>1$. This is a contradiction which gives the result.
\end{proof}
An extension of Theorem \ref{t:3.10} will be given in Section
\ref{s:7}.

\section{Existence and compactness of a critical element; further properties of critical
elements}\label{s:4}

Let us consider the statement:

(SC) For all $\rt\in\hl$, with $\int\abs{\na u_0}^2<\int\abs{\na
W}^2$ and $E(\rt)<E((W,0))$, if $u$ is the corresponding solution of
(CP) with maximal interval of existence $I$ (see Definition
\ref{d:2.13}) then $I=(-\infty,+\infty)$ and
$\norm{u}_{S((-\infty,+\infty))}<\infty$.

In addition, for a fixed $\rt\in\hl$, with $\int\abs{\na u_0}^2<\int\abs{\na
W}^2$ and $E(\rt)<E((W,0))$, we say that $(SC)(\rt)$ holds if, for
$u$ the corresponding solution of (CP), \wi,
we have $I=(-\infty,+\infty)$ and
$\norm{u}_{S(-\infty,+\infty)}<\infty$.

Note that, because of Remark \ref{r:2.10}, if
$\norm{\rt}_{\hl}\leq\td$, then $SC(\rt)$ holds. Thus, in light of
Corollary \ref{c:3.9}, there exists $\eta_0>0$ such that if $\rt$ is
as in (SC), and $E(\rt)\leq\eta_0$, then $SC(\rt)$ holds. Moreover,
for any $\rt$ as in (SC), \eqref{eq:3.11} shows that $$E(\rt)\geq0.$$
Thus, there exists a number $E_C$, $\eta_0\leq E_C\leq E((W,0))$
such that, if $\rt$ is as in (SC) and $E(\rt)<E_C$, then $(SC)(\rt)$
holds and $E_C$ is optimal with this property. For the rest of this
section we will assume that $E_C<E((W,0))$. Using concentration
compactness ideas, following the argument in \cite{15}, Section 4,
we prove that there exists a critical element $(u_{0,C},u_{1,C})$ at
the critical level of energy $E_C$, so that $SC((u_{0,C},u_{1,C}))$
does not hold and from the minimality, this element has a
compactness property up to the symmetries of the equation (which will give rigidity in the problem). We then
use the finite speed of propagation and Lorentz transformations to
establish support and orthogonality properties of critical elements, which are essential to treat the nonradial case.

\begin{prop}\label{p:4.1}
There exists $(u_{0,C},u_{1,C})$ in $\hl$, with
$$E((u_{0,C},u_{1,C}))=E_C<E((W,0)),\quad \int\abs{\na
u_{0,C}}^2<\int\abs{\na W}^2$$ such that if $u_C$ is the solution of
$(CP)$ with data $(u_{0,C},u_{1,C})$ and \wi,
$0\in\overset{\circ}I$, then $\norm{u_C}_{S(I)}=+\infty$.
\end{prop}

\begin{prop}\label{p:4.2}
Assume that $u_C$ is as in Proposition \ref{p:4.1} and that
$($say$)$ $\norm{u_C}_{S(I_+)}=+\infty$, where
$I_+=[0,\infty)\bigcap I$. Then there exists $x(t)\in\bR^N$,
$\la(t)\in \bR^+$, for $t\in I_+$, such that $K=
\left\{\vec{v}(x,t),t\in I_+ \right\}$ has the property that $\overline{K}$ is compact in $\hl$, where
$$
\vec{v}(x,t)=\left(\frac{1}{\la(t)^{\frac{(N-2)}{2}}}u_C\left(\frac{x-x(t)}{\la(t)},t
\right)\right.,
\left.\left.\frac{1}{\la(t)^{\frac N 2}}\dd_tu_C\left(\frac{x-x(t)}{\la(t)},t
\right)\right).\right.
$$
A
corresponding conclusion is reached if
$\norm{u_C}_{S(I_-)}=+\infty$, where $I_-=(-\infty,0)\bigcap I$.
\end{prop}

The proofs of Propositions \ref{p:4.1} and \ref{p:4.2} are identical
to the corresponding ones in \cite{15}, using Lemma \ref{l:4.3}
below and the results of Section \ref{s:2}, especially Theorem
\ref{t:2.21}. We will therefore omit them.

\begin{lem}[Concentration compactness]\label{l:4.3}
Let $\{(v_{0,n},v_{1,n})\}\in\hl$,
$\norm{(v_{0,n},v_{1,n})}_{\hl}\leq A$. Assume that
$$\norm{S(t)((v_{0,n},v_{1,n}))}_{S(-\infty,+\infty)}\geq\de> 0,$$
where $\de=\de(A)$ is as in Theorem \ref{t:2.7}. Then there exists
a sequence $\{(V_{0,j},V_{1,j}) \}$ in $\hl$, a subsequence of
$\{(v_{0,n},v_{1,n})\}$ $($which we still call
$\{(v_{0,n},v_{1,n})\})$ and a triple
$(\la_{j,n};x_{j,n};t_{j,n})\in\bR^+\times\bR^N\times\bR$ with
$$
\frac{\la_{j,n}}{\la_{j',n}}+\frac{\la_{j',n}}{\la_{j,n}}+\frac{\abs{t_{j,n}-t_{j',n}}}{\la_{j,n}}+
\frac{\abs{x_{j,n}-x_{j',n}}}{\la_{j,n}}\to\infty
$$
as $n\to\infty$, for $j\ne j'$ $($we say that
$(\la_{j,n};x_{j,n};t_{j,n})$ is orthogonal if this property is
verified$)$ such that
 \bbe\label{eq:4.4}
\norm{\left(V_{0,1},V_{1,1}\right)}_{\hl}>\al_0(A)>0.
 \ee
If $V_j^l(x,t)=S(t)\left((V_{0,j},V_{1,j}) \right)$, then given
$\eps_0>0$, there exists $J=J(\eps_0)$ and $\left\{(w_{0,n},w_{1,n}) \right\}_{n=1}^{\infty}\in\hl,$ so
that
 \begin{eqnarray}\label{eq:4.5}
&&
v_{0,n}=\sum_{j=1}^J\frac{1}{\la_{j,n}^{\frac{(N-2)}{2}}}V^l_j\left(\frac{x-x_{j,n}}{\la_{j,n}},-\frac{t_{j,n}}{\la_{j,n}}
\right)+w_{0,n},\\
&&
v_{1,n}=\sum_{j=1}^J\frac{1}{\la_{j,n}^{\frac{N}{2}}}\dd_tV^l_j\left(\frac{x-x_{j,n}}{\la_{j,n}},-\frac{t_{j,n}}{\la_{j,n}}
\right)+w_{1,n},\nonumber\\
&& \text{ with
}\norm{S(t)((w_{0,n},w_{1,n}))}_{S(\infint)}\leq\eps_0,\text{ for
}n\text{ large}\nonumber
 \end{eqnarray}
\begin{eqnarray}\label{eq:4.6}
&&\int\frac{1}{2}\abs{\na_xv_{0,n}}^2+\frac{1}{2}\abs{v_{1,n}}^2=\sum_{j=1}^J\int\frac{1}{2}\abs{\na_xV_{0,j}}^2+\frac{1}{2}\abs{V_{1,j}}^2 \\
&&+\int\frac{1}{2}\abs{\na_xw_{0,n}}^2+\frac{1}{2}\abs{w_{1,n}}^2+o(1)\text{
as }n\to\infty\nonumber
\end{eqnarray}
\begin{eqnarray}\label{eq:4.7}
&& E((v_{0,n}v_{1,n}))=\sum_{j=1}^JE\left(V_j^l\left(
-\frac{t_{j,n}}{\la_{j,n}}\right),\dd_tV_j^l\left(
-\frac{t_{j,n}}{\la_{j,n}}\right)\right)+\\
&& +E((w_{0,n},w_{1,n}))+o(1)\text{ as }n\to \infty.\nonumber
\end{eqnarray}

\end{lem}
\begin{rem}\label{r:4.8}
Lemma \ref{l:4.3} is due to Bahouri-G\'erard \cite{4}. There it is
proved for $N=3$, but the proof extends to all $N\geq 3$. Also, the
norm $\norm{\,\,}_{S\infint}$ is replaced by
$\norm{\,\,}_{L_t^{\frac{(N+2)}{N-2}}L_x^{2\frac{(N+2)}{N-2}}}$ in \cite{4}, but
as is mentioned in page 136 of \cite{4}, it works equally well for
$\norm{\,\,}_{S\infint}$. See the Remark on page 159 of \cite{4} to
eliminate their condition (1.6). (See also the work of Keraani
\cite{17}, where the corresponding result is proved for NLS and
where the analogue of \eqref{eq:4.4} is shown.) See also Remark 4.8
in \cite{15}.
\end{rem}

\begin{cor}\label{c:4.9}
There exists a decreasing function $g:(0,E_C]\to[0,\infty)$ such
that for every $\rt$ as in $(SC)$, with $E(\rt)=E_C-\eta$, we have
$$
\norm{u}_{S(\infint)}\leq g(\eta).
$$
\end{cor}
For a proof of Corollary \ref{c:4.9}, see Corollary 2 in\cite{4} and
Corollary 1.14 in \cite{17}.

We next turn our attention to further properties of critical
elements.

\begin{lem}\label{l:4.10}
Assume that $u$ is a solution of $(CP)$, \wi. Assume that for $t\in
I^+=I\cap[0,\infty)$, there exist $x(t)\in\bR^N$, $\la(t)\in \bR^+$
so that $K=\left\{\vec{v}(x,t),t\in I_+ \right\}$ has the property that $\overline{K}$ is compact in $\hl$, where 
$$\vec{v}(x,t)=\left(\frac{1}{\la(t)^{\frac{(N-2)}{2}}}u\left(\frac{x-x(t)}{\la(t)},t
\right)\right.,
\left.\left.\frac{1}{\la(t)^{\frac{N}{2}}}\dd_tu\left(\frac{x-x(t)}{\la(t)},t
\right)\right)\right..$$
 Then we
can choose $\tilde{\la}(t),\tilde{x}(t)$, continuous in $I_+$, so
that the corresponding $\tilde{K}$ has compact closure in $\hl$.

\end{lem}

\begin{proof}
The proof given in Remark 5.4 of \cite{15} applies verbatim.
\end{proof}

From now on, we always use the $\tilde{\la}(t)$, $\tilde{x}(t)$
provided by Lemma \ref{l:4.10}.

\begin{lem}\label{l:4.11}
Let $u$ be as in Lemma \ref{l:4.10} and assume that $I_+$ is a
finite interval. After scaling, we can assume then that $I_+=[0,1)$.
Then, $$0<\frac {C_0(K)}{1-t}\leq \la(t).$$
\end{lem}
\begin{proof}

Consider $0<t_j\to 1$. (Because of Lemma \ref{l:4.10}, this
suffices.) Let $$(v_{0,j},
v_{1,j})=\left(\frac{1}{\la(t_j)^{\frac{N-2}{2}}}u\left(\frac{x-x(t_j)}{\la(t_j)}
,t_j\right),\frac{1}{\la(t_j)^{\frac{N}{2}}}\dd_tu\left(\frac{x-x(t_j)}{\la(t_j)},t_j \right)\right).$$ Since $(v_{0,j},v_{1,j})\in K$,
$\overline{K}$ is compact in $\hl$, there exists $C_0=C_0(K)>0$
independent of $j$, so that $T_+\left((v_{0,j},v_{1,j})\right)\geq
C_0$. (Here we are using the notation in Definition \ref{d:2.13}.)
(This is an easy consequence of Theorem \ref{t:2.7}.) Let $v_j(t)$
be the corresponding solution of (CP). Note that $\la(t_j)^{\frac{(N-2)}{2}}
v_{0,j}(\la(t_j)y+x(t_j))=u(y,t_j), \la(t_j)^{\frac{N}{2}}
v_{1,j}(\la(t_j)y+x(t_j))=\dd_tu(y,t_j).$ Hence, by uniqueness in
(CP) (see the argument in Definition \ref{d:2.13}) for $t$ such that
$t_j+t\leq T_+\left(\rt \right)=1$ we
have
$$\la(t_j)^{\frac{(N-2)}{2}}v_j(\la(t_j)y+x(t_j),\la(t_j)t)=
u(y,t_j+t).$$ Thus, we have
$t_j+t\leq 1$, for all $0<\la(t_j)t\leq C_0$. But then, choose
$t=C_0/\la(t_j)$, so that $\la(t_j)\geq C_0/(1-t_j)$, as desired.
\end{proof}

\begin{lem}\label{l:4.12}
Let $u$ be as in Lemma \ref{l:4.11}. Then $\exists \bar{x}\in\bR^N$
such that $$\supp u\subset B(\bar{x},(1-t)), \supp \dd_tu\subset
B(\bar{x},(1-t)).$$
\end{lem}

\begin{proof}
We first start by showing that for $t\in[0,1)$, there is a ball
$B_{(1-t)}$ of radius $(1-t)$ so that $\supp\na u$, $\supp\dd_t u$
are contained in $B_{(1-t)}$. If not, for a fixed $t$, there exist
$\eps_0>0$, $\eta_0>0$ such that, for all $x_0\in\bR^N$ we have
$$
\int_{\abs{x-x_0}\geq(1+\eta_0)(1-t)}\abs{\na_x
u(t)}^2+\abs{\dd_tu(t)}^2\geq\eps_0.
$$

Choose a sequence $t_n\uparrow 1$. Recall from Lemma \ref{l:4.11}
that $\la(t_n)\geq \frac{C_0}{1-t_n}$. We claim that, given $R_0>0$,
$M>0$, for $n$ large we have
$$
\int_{\abs{x+\frac{x(t_n)}{\la(t_n)}}\geq R_0}\abs{\na_x
u(x,t_n)}^2+\abs{\dd_tu(x,t_n)}^2+\frac{\abs{u(x,t_n)}^2}{\abs{x}^2}\leq\frac{\eps_0}M.
$$
Indeed, if
$
\vec{v}(x,t)=\frac{1}{\la(t)^{\frac{N}{2}}}\left(\na
u\left(\frac{x-x(t)}{\la(t)},t
\right),\dd_tu\left(\frac{x-x(t)}{\la(t)},t \right)\right),
$
$$
\int_{\abs{x+\frac{x(t_n)}{\la(t_n)}}\geq R_0}\abs{\na_x
u(x,t)}^2+\abs{\dd_tu(x,t)}^2dx=\int
_{\abs{y}\geq\la(t)R_0}\abs{\vec{v}(y,t)}^2dy
$$
and our claim follows from the compactness of $\overline{K}$ and the
fact that $\la(t_n)\uparrow +\infty$. Using this estimate, we apply
Lemma \ref{l:2.18} backward in time, to conclude that for $n$ large, 
$$
\forall t\in[0,t_n], \ \ \int_{\abs{x+\frac{x(t_n)}{\la(t_n)}}\geq 2R_0+(t_n-t)}\abs{\na_x
u(x,t)}^2+\abs{\dd_tu(x,t)}^2\leq\eps_0.
$$
But, if $(1+\eta_0)(1-t)\geq 2R_0+(t_n-t)$, we reach a
contradiction. But, for $0\leq t<1$, fixed, we can always choose $n$
large and $R_0$ small so that this is the case.

The next step is to show that $\abs{\frac{x(t)}{\la(t)}}\leq M$, for
$0\leq t<1$. Assume not, so that we can find (in light of Lemma
\ref{l:4.10} we can assume $\abs{\frac{x(t)}{\la(t)}}\leq M_T$,
$0\leq t\leq T<1$) $t_n\uparrow 1$ so that
$\abs{\frac{x(t_n)}{\la(t_n)}}\to +\infty$. Fix a ball $B=B(x_0,1)$,
such that $\supp\na u_0,\supp u_1\subset B$. But, for fixed $R_0>0$,
$\eps_0>0$ given, our previous argument shows that for $n$ large,
$$
\int_{\abs{x+\frac{x(t_n)}{\la(t_n)}}\geq 2R_0+t_n}\abs{\na_x
u_0}^2+\abs{u_1}^2\leq\eps_0.
$$
But, if $\abs{\frac{x(t_n)}{\la(t_n)}}\to +\infty$,
$B(x_0,1)\subset\left\{\abs{x+\frac{x(t_n)}{\la(t_n)}}\geq 2R_0+t_n
\right\}$, so that $\na u_0$, $u_1$ are identically $0$,
contradicting $I_+=[0,1)$. Let now $t_n\uparrow 1$, and choose a
subsequence so that $-\frac{x(t_n)}{\la(t_n)}\to\bar{x}$. Arguing as
before, for $0\leq t< t_n$, we see that, for $n$ large,
$$
\int_{\abs{x+\frac{x(t_n)}{\la(t_n)}}\geq R_0}\abs{\na_x
u(t_n)}^2+\abs{\dd_tu(t_n)}^2\leq\eps_0/M,
$$
for $R_0,M,\eps_0$ given and hence, by our previous argument,
$$
\int_{\abs{x+\frac{x(t_n)}{\la(t_n)}}\geq 2R_0+(t_n-t)}\abs{\na
u(x,t)}^2+\abs{\dd_tu(x,t)}^2dx\leq\eps_0,
$$
for $n$ large. Letting $n\to\infty$, we obtain for all $R_0>0$ and $\eps_0>0$ small,
$$
\int_{\abs{x-\bar{x}}\geq 2R_0+(1-t)}\abs{\na
u(x,t)}^2+\abs{\dd_tu(x,t)}^2dx\leq\eps_0,
$$
so that $\supp\na u(-,t)$, $\supp\dd_tu(-,t)\subset B(\bar{x},1-t)$.
Assume now that $-\frac{x(t_n)}{\la(t_n)}\to\bar{x}$,
$-\frac{x(t_n')}{\la(t_n')}\to\bar{x}'$ for two different sequences
$t_n,t_n'\to 1$. If $\bar{x}\ne \bar{x}'$ and $(1-t)$ is so small
that $1-t<\abs{\bar{x}- \bar{x}'}$, we must have $\na
u(-,t),\dd_tu(-,t)\equiv 0$, a contradiction to $I_+=[0,1)$.
\end{proof}

\begin{rem}\label{r:4.13}
After a translation we can assume $\bar{x}=0$. Also, since
$u(-,t)\in L^{2^*}$ for each $t$, the conditions $\supp u\subset
B(0,1-t)$ and $\supp\na_x u\subset B(0,1-t)$ are equivalent.

\end{rem}

We turn now to the next important property of $u_C$ (at least in the nonradial situation): the second invariant of the equation for $u_C$ is zero. We consider the cases $I_+$ is a
finite interval and then an infinite interval.
\begin{prop}\label{p:4.14}
Assume that $u_C$ is as in Proposition \ref{p:4.2} and $I_+$ is a
finite interval. Then,
$$
\int\na u_{0,C}.u_{1,C}=0.
$$
\end{prop}

\begin{proof}
By scaling, we can assume that $I_+=[0,1)$. By Lemma \ref{l:4.12},
$\supp u_C,\supp\dd_t u_C\subset B(0,1-t)$. Note also that for any
$u$ a solution of (CP) in $I$, the maximal interval of existence,
and $t\in I$, we have from (\ref{eq:4.15}),
$
\int\na_xu(t)\dd_tu(t)dx=\int\na u_0.u_1.
$ Assume now that, (without
loss of generality)
$$
\ga=\int\dd_{x_1}(u_{0,C}).u_{1,C}>0.
$$
We will reach a contradiction, by considering (for convenience)
$u(x,t)=u_C(x,1+t)$, $-1\leq t< 0$. Clearly, for $-1\leq t< 0,$
$$E((u(t),\dd_tu(t)))=E_C, \ \int\abs{\na u(t)}^2\leq
(1-\bar{\de})\norm {\na W}^2_{L^2}, \ga=\int\dd_{x_1}u(t).\dd_tu(t),$$  by Theorem \ref{t:3.8} and our assumption above. We will consider the
action of Lorentz transformations on $u$. (Now, $\supp u(-,t)\subset
B(0,-t)$, $\supp\dd_t u(-,t)\subset B(0,-t),$ $-1\leq t< 0$.) Thus,
for $0<d<1/4$, consider
 \bbe\label{eq:4.16}
z_d(x_1,\bar{x},t)=u\left(\frac{x_1-dt}{\sqrt{1-d^2}},\bar{x},
\frac{t-dx_1}{\sqrt{1-d^2}} \right),
 \ee
where $x=(x_1,\bar{x})\in\bR^N$, $t\in\bR$ and
$s=\frac{t-dx_1}{\sqrt{1-d^2}}$ is such that $-1\leq s< 0$. 

Note
that, for this range of $s$ and $y=(y_1,\bar{y})$ such that
$(y,s)\in\supp u$, we have $\abs{y}\leq\abs{s}$. Thus, if
$y_1=\frac{x_1-dt}{\sqrt{1-d^2}}$, $\bar{y}=\bar{x}$, we obtain
$x_1^2+\abs{\bar{x}}^2\leq t^2$ in support of $z_d,\dd_t z_d$. Fix
now $-\frac{1}{2}\leq t<0$ and $x_1^2+\abs{\bar{x}}^2\leq t^2$.
Then, $\frac{t-dx_1}{\sqrt{1-d^2}}\geq
\frac{(1+d)t}{\sqrt{1-d^2}}\geq-\frac{1}{2}\frac{1+d^2}{\sqrt{1-d^2}}\geq
-1$, while $\frac{t-dx_1}{\sqrt{1-d^2}}\leq
\frac{(1-d)t}{\sqrt{1-d^2}}<0$. Thus, for such $(x,t)$, $z_d$ is
defined and $z_d(x,t)=0$, $\na_xz_d(x,t)=0$, $\dd_tz_d(x,t)=0$ for
$x_1^2+\abs{\bar{x}}^2=t^2$. We extend $z_d(-,t)$ to be zero for
$\abs{x}\geq\abs{t}$, $-\frac{1}{2}\leq t< 0$. An elementary
calculation shows that if $u$ is a regular solution (by regular
solution we will mean one as in Remark \ref{r:2.9}, with $\mu=1$) of
$$
\dd_t^2u-\De u=\abs{u}^{\frac{4}{N-2}}u \quad\text{in }\bR^N\times[-1,0)
$$
the resulting $z_d$ is a solution of the (CP) for this equation in
$-\frac{1}{2}\leq t<0$, $x\in\bR^N$.

We will now show that the $z_d$ we defined in \eqref{eq:4.16} is a
solution of (CP) in $\bR^N\times[-\frac12,0)$. To this end, fix
$\eps_0>0$ and consider $-\frac{1}{2}\leq t\leq -\eps_0$,
$x\in\bR^N$. Note  that in this range we have, on $\supp z_d$, that
$-1\leq s\leq-\frac{3}{\sqrt{15}}\eps_0$. Note that since $S\left([-1,(\frac{-3}{\sqrt{15}})\eps_0]\right)$ norm of $u$ is finite, and $u\in
L^{\frac{(N+2)}{N-2}}_{[-1,-(\frac{3}{\sqrt{15}})\eps_0]}L_x^{2\frac{(N+2)}{N-2}}$ (see
Definition \ref{d:2.13}), in light of Remark \ref{r:2.8}, we have
that $(z_d,\dd_t(z_d))\in C([-1/2,-\eps_0];\hl)$. Also, if we let
$J=\abs{\det\frac{\dd (y,s)}{\dd(x,t)}}$, then $J\equiv 1$ and
hence, if $D_{\eps_0}=\bR^N\times [-1/2,-\eps_0]$,
$\tilde{D}_{\eps_0}=\Phi(D_{\eps_0})$, where $\Phi(x,t)=(y,s)$, then
 \bbea
&&
\int_{D_{\eps_0}}\abs{z_d(x,t)}^{\frac{2(N+1)}{(N-1)}}dxdt=\int_{\tilde{D}_{\eps_0}}\abs{u(y,s)}
^{\frac{2(N+1)}{(N-1)}}dyds \\
&& \leq \int_{-1\leq s\leq
-\frac{3}{\sqrt{15}}\eps_0}\abs{u(y,s)}^{\frac{2(N+1)}{(N-1)}}dyds\leq C_{\eps_0}.
 \eea

Moreover, pick $u_{0,j}\in C^{\infty}_0(B(0,\frac{3}{\sqrt{15}}\eps_0)), \quad
u_{1,j}\in C^{\infty}_0(B(0,\frac{3}{\sqrt{15}}\eps_0))$ with
$(u_{0,j},u_{1,j})\to
(u(-(\frac{3}{\sqrt{15}})\eps_0),\dd_su(-(\frac{3}{\sqrt{15}})\eps_0))$ in $\hl$.
Let $u_j$ be the solution of (CP), defined for
$s<-\left(\frac{3}{\sqrt{15}}\right)\eps_0$. Note that, because of
Remark \ref{r:2.22}, we know that, for $j$ large, $u_j$ is a
solution of (CP) for $-1\leq
s<-\left(\frac{3}{\sqrt{15}}\right)\eps_0$,
$$\norm{(u_j,\dd_s(u_j))}_{C\left((-1,-\frac{3}{\sqrt{15}}\eps_0),\hl
\right)}\leq C$$ and
$$\norm{u_j}_{S\left(-1,-\frac{3}{\sqrt{15}}\eps_0\right)}
 \ + \ \norm{u_j}_{L^{\frac{(N+2)}{N-2}}_{\left[-1,-\frac{3}{\sqrt{15}}\eps_0
\right]}L_x^{\frac{2(N+1)}{N-2}}}\leq \tilde{C}_{\eps_0}.$$ 
Also, by
virtue of Remark \ref{r:2.9}, $u_j$ is regular for
$s\in\left[-1,-\frac{3}{\sqrt{15}}\eps_0 \right]$ and for $-1\leq
s\leq-\frac{3}{\sqrt{15}}\eps_0$, we have $\supp u_j(-,s)\subset
B(0,\abs{s})$, by Remark \ref{r:2.12}. If we now consider $z_{j,d}$
given by \eqref{eq:4.16} with $u$ replaced by $u_j$, the $z_{j,d}$
are solutions of (CP) in $-\frac{1}{2}\leq t\leq -\eps_0$. Moreover,
from the proof of Remark \ref{r:2.22} and the proof that
$(z_d,\dd_t(z_d))\in C\left([-1/2,-\eps_0];\hl \right)$ we can
conclude that $(z_{j,d},\dd_t(z_{j,d}))\to (z_d,\dd_tz_d)$ in
$C([-1/2,-\eps_0];\hl)$ and similarly that $\norm{z_{j,d}}
_{S([-1/2,-\eps_0])}\leq C_{\eps_0}$. From Remark \ref{r:2.14} it
now follows that $z_d$ is a solution of (CP) for
$t\in[-1/2,-\eps_0]$. Since $\eps_0>0$ is arbitrary, we conclude
that $T_+((z_d(-1/2),\dd_tz_d(-1/2)))\geq 0$. But, since for each
$t\in[-1/2,0)$, $\supp z_d,\dd_tz_d\subset\{\abs{x}\leq\abs{t}\}$,
either $T_+((z_d(-1/2),\dd_tz_d(-1/2)))=0$, or $z_d\equiv 0$. We
will soon see that $z_d\not\equiv 0$.

We have, by Remark \ref{r:2.16}, that
 \bbe\label{eq:4.17}
\int_{-1/2}^{-1/4}E((z_d(t),\dd_tz_d(t)))dt=\frac{1}{4}E((z_d(-1/2),\dd_tz_d(-1/2))).
 \ee
We are now going to estimate the left-hand side. Note that
 \bbea
\dd_{x_1}z_d &=&-\frac{d}{\sqrt{1-d^2}}\dd_s
u+\frac{1}{\sqrt{1-d^2}}\dd_{y_1}u, \ \ \
\dd_{\bar{x}}z_d =\dd_{\bar{y}}u,\\
\dd_{t}z_d &=&\frac{1}{\sqrt{1-d^2}}\dd_s
u-\frac{d}{\sqrt{1-d^2}}\dd_{y_1} u.\\
 \eea
Thus, the left-hand side of \eqref{eq:4.17} equals $I_1+I_2$, where
 \bbea
&&
I_1=\int_{-1/2}^{-1/4}\int \frac{1}{2}\left\{\left(\frac{1+d^2}{1-d^2}\right)((\dd_su)^2+
(\dd_{y_1}u)^2)
+\abs{\na_{\bar{y}}u}^2\right\}
-\left.\frac{1}{2^*}\abs{u}^{2^*}\right.dx_1d\bar{x}dt \\
&&
I_2=-\frac{2d}{(1-d^2)}\int_{-1/2}^{-1/4}\int\dd_{y_1}u\dd_sudx_1d\bar{x}dt.
 \eea
 We now have
 \bbe\label{eq:4.18}
\lim_{d\downarrow
0}\frac{I_2}{d}=-2\left(\frac{1}{4}\right)\int\dd_{y_1}u(-\frac12).\dd_tu(-\frac12)=-\frac{1}{2}\ga.
 \ee
To see this, we consider the change of variables $\Phi(x,t)=(y,s)$,
introduced before. Let $D_{-1/4}=\bR^N\times [-1/2,-1/4]$,
$\tilde{D}_{-1/4}=\Phi(D_{-1/4})$. Since
$\abs{\det\frac{\dd(y,s)}{\dd(x,t)}}=1$, we have,
 \bbea
\lim_{d\downarrow 0}\frac{I_2}{d}&=&\lim_{d\downarrow
0}-2\int\int_{D_{-1/4}}\dd_{y_1}u.\dd_sudx_1d\bar{x}dt\\
&=&\lim_{d\downarrow
0}-2\int\int_{\tilde{D}_{1/4}}\dd_{y_1}u.\dd_sudy_1d\bar{y}ds.
 \eea

Since $y_1=\frac{x_1-dt}{\sqrt{1-d^2}}$,
$s=\frac{t-dx_1}{\sqrt{1-d^2}}$, we have that
$t\sqrt{1-d^2}=dy_1+s$, so that
$
\tilde{D}_{1/4}=\left\{(y,s):\abs{y}\leq\abs{s}\text{ and
}-\frac{1}{2}\sqrt{1-d^2}\leq dy_1+s\leq -\frac{1}{4}\sqrt{1-d^2}
\right\}.
$
Note that the restriction $\abs{y}\leq \abs{s}$ comes from the
support of $\dd_{y_1}u$, $\dd_su$. Thus,
 \bbea
&&
\int\int_{\tilde{D}_{1/4}}\dd_{y_1}u.\dd_sudyds=\int_{\abs{y}\leq
2}\int\int_{-\frac{1}{2}\sqrt{1-d^2}-dy_1}^{-\frac{1}{4}\sqrt{1-d^2}-dy_1}\dd_{y_1}u\dd_sudsdy\\
&& =\int_{\abs{y}\leq
2}\int\int_{-\frac{1}{2}}^{-\frac{1}{4}}\dd_{y_1}u.\dd_sudsdy-
\int_{\abs{y}\leq
2}\int\int_{-\frac{1}{2}}^{-\frac{1}{2}\sqrt{1-d^2}-dy_1}\dd_{y_1}u\dd_sudsdy\\
&& -\int_{\abs{y}\leq
2}\int\int_{-\frac{1}{4}\sqrt{1-d^2}-dy_1}^{-\frac{1}{4}}\dd_{y_1}u\dd_sudsdy.
 \eea
Consider, for instance, the second term. There, $-\frac{1}{2}\leq
s\leq -\frac{1}{2}\sqrt{1-d^2}+2d$, so that, it is bounded by
$$
\int_{-\frac{1}{2}}^{-\frac{1}{2}\sqrt{1-d^2}+2d}\int_{\abs{y}\leq
2}\abs{\dd_{y_1}u}\abs{\dd_su}dyds\leq Cd\sup_{-\frac{1}{2}\leq
s\leq -\frac{1}{2}\sqrt{1-d^2}+2d}\norm{(u,\dd_su)}^2_{\hl}.
$$
Thus, the second term goes to $0$ as $d\to 0$, and the third one can
be treated similarly. (Recall that, by compactness, we have
\hfill\newline $\sup_{-1\leq s\leq
0}\norm{(u(s),\dd_su(s))}_{\h}\leq A$.) But, using now
\eqref{eq:4.15}, we obtain \eqref{eq:4.18}.

\medbreak

To study $I_1$, we introduce
$$
I_3=\int_{-\frac{1}{2}\sqrt{1-d^2}}^{-\frac{1}{4}\sqrt{1-d^2}}\int_{\abs{y}\leq
2}\frac{1}{2}\left\{\frac{1+d^2}{1-d^2}((\dd_su)^2+(\dd_{y_1}u)^2)+\abs{\na_{y_1}u}^2
\right\}-\frac{\abs{u}^{2^*}}{2^*}dyds.$$
%
Then, in light of the fact that $\sup_{-1\leq s\leq
0}\norm{(u(s),\dd_su(s))}_{\hl}\leq A$, the identity in Remark
\ref{r:2.16} and the support properties of $u$, $\dd_su$, we have
$$
I_3=\frac{1}{4}E(u(-1/2),\dd_su(-1/2))+O(d^2)=\frac{1}{4}E_C+O(d^2).
$$
We next claim that
 \bbe\label{eq:4.19}
\lim_{d\downarrow 0}\frac{I_1-I_3}{d}=-\int y_1e(u)(-1/4)+\int
y_1e(u)(-1/2).
 \ee
(Recall the definition of $e(u)$ from the proof of Remark
\ref{r:2.16}.) Let us assume \eqref{eq:4.19} temporarily. Recall
that\eqref{eq:2.17}, the support properties of $u$ and integration
by parts, yield
$$
\dd_s\int y_1e(u(s))dy=-\int\dd_{y_1}u(s)\dd_su(s)dy,
$$
so that in light of \eqref{eq:4.15}, $-\int y_1e(u)(-1/4)+\int
y_1e(u)(-1/2)=\frac{1}{4}\ga$ and hence $\lim_{d\downarrow
0}\frac{I_1-\frac{1}{4}E_C}{d}=\frac{1}{4}\ga$ and so, using
\eqref{eq:4.18}, we obtain from \eqref{eq:4.17},
 \bbe\label{eq:4.20}
\lim_{d\downarrow
0}\frac{\frac{1}{4}E\left((z_d(-1/2),\dd_tz_d(-1/2))-\frac{1}{4}E_C
\right)}{d}=\frac{1}{4}\ga-\frac{1}{2}\ga=-\frac{1}{4}\ga.
 \ee
(Note that, for $d$ small this already implies that $z_d$ cannot by
identically $0$.) \eqref{eq:4.20} implies that, for $d>0$ small,
$E\left(z_d(-1/2),\dd_tz_d(-1/2) \right)<E_C$, since $\ga>0$.

We now turn to the verification of \eqref{eq:4.19}.  Note that
$$
I_1=\int_{-1/2}^{-1/4}\int_{\abs{x}\leq
2}\frac{1}{2}\left\{(\dd_su)^2+\abs{\na_yu}^2
\right\}-\frac{1}{2^*}\abs{u}^{2^*}dxdt+ O(d^2),
$$
$$
I_3=\int_{-1/2\sqrt{1-d^2}}^{-1/4\sqrt{1-d^2}}\int_{\abs{y}\leq
2}\frac{1}{2}\left\{(\dd_su)^2+\abs{\na_yu}^2
\right\}-\frac{1}{2^*}\abs{u}^{2^*}dyds+ O(d^2).
$$
By using the change of variables used in the proof of
\eqref{eq:4.18} we see that
 \bbea
\lim_{d\downarrow 0}\frac{I_1-I_3}{d}&=&\lim_{d\downarrow
0}\left\{\frac{1}{d}\int_{\abs{y}\leq
2}\int_{-\frac{1}{4}\sqrt{1-d^2}}^{-\frac{1}{4}\sqrt{1-d^2}-dy_1}e(u)(y_1,\bar{y},s)dsdy\right.\\
&& \left.-\frac{1}{d}\int_{\abs{y}\leq
2}\int_{-\frac{1}{2}\sqrt{1-d^2}}^{-\frac{1}{2}\sqrt{1-d^2}-dy_1}e(u)(y_1,\bar{y},s)dsdy
\right\}\\ & =&\lim_{d\downarrow 0}[A(d)-B(d)].
 \eea
 \bbea
A(d)&=&\frac{1}{d}\int_{\abs{y}\leq
2}\int_{-\frac{1}{4}\sqrt{1-d^2}}^{-\frac{1}{4}\sqrt{1-d^2}-dy_1}e(u)(y_1,\bar{y},s)dsdy\\
&=&-\int_{\abs{y}\leq
2}\int_{0}^{y_1}e(u)(y_1,\bar{y},\sqrt{1-d^2}(-1/4)-hd)dhdy,
 \eea
where we have made the change of variables
$h=\frac{(-1/4)\sqrt{1-d^2}-t}{d}$. Since $(\na u,\dd_tu)\in
C([-1,-\eps_0];L^2(\bR^N))$, for every $\eps_0>0$, we see that, for
$d$ small, we have, for $\abs{h}\leq 2$, that
$$
\int_{\abs{y}\leq
2}[e(u)(y,\sqrt{1-d^2}(-1/4)-hd)-e(u)(y,-1/4)]dy=o(1)
$$
as $d\to 0$, uniformly in $h$. Hence,
\bbea
A(d)&=&-\int_{\abs{y}\leq 2}\int_0^{y_1}e(u)(y,-1/4)dhdy+o(1)\\
&=&-\int y_1e(u)(y,-1/4)dy+o(1).
\eea
Similarly, $B(d)=-\int y_1e(u)(y,-1/2)dy+o(1)$ and hence
\eqref{eq:4.19} follows.

Finally, since $E_C<E((W,0))$, $\norm{\na u(-1)}^2_{L^2}<\norm{\na
W}^2_{L^2}$, because of Theorem \ref{t:3.8}, we have that, for
$-1\leq s\leq 0$, $\norm{\na_xu(-,s)}^2_{L^2}\leq
(1-\bar{\de})\norm{\na W}^2_{L^2}$, $\bar{\de}>0$. We now consider
$\int_{-1/2}^{-1/4}\int\abs{\na_xz_d(x,t)}^2dxdt$. The argument that
we used in the estimate for $\frac{I_2}{d}$ above, (together with
the calculation of $\dd_{x_1}z_d$, $\dd_{\bar{x}}z_d$) show that
 \bbea
 \lim_{d\downarrow
0}\int_{-1/2}^{-1/4}\int\abs{\na_xz_d(x,t)}^2dxdt&=&
\int_{-1/2}^{-1/4}\int\abs{\na_yu(y,s)}^2dyds\\
& \leq &\frac{1}{4}(1-\bar{\de})\norm{\na W}^2_{L^2}.
 \eea
But then, for $d$ small,
$\int_{-1/2}^{-1/4}\int\abs{\na_xz_du(x,t)}^2dxdt\leq\frac{1}{4}(1-\bar{\de}/2)\norm{\na
W}^2_{L^2}$. Thus, there exists $t_0=t_0(d)\in (-1/2,-1/4)$ such
that, for $d$ small, $$\int\abs{\na_xz_d(x,t_0)}^2dx <\norm{\na
W}^2_{L^2},\quad E(z_d(-1/2),\dd_tz_d(-1/2))<E((W,0)).$$ By Theorem
\ref{t:3.8} we have, for all $d$ small,
$\norm{\na_xz_d(x,-1/2)}^2_{L^2}<\norm{\na W}^2_{L^2}$. Since the
interval of existence of $z_d$ is finite, this contradicts the
definition of $E_C$ taking $d>0$ small, and thus $\ga=0$.
\end{proof}
\begin{prop}\label{p:4.21}
Assume that $u_C$ is as in Proposition \ref{p:4.2} and
$I_+=[0,+\infty)$. Assume in addition that for $t>0$, $\la(t)>A_0>0$. Then,
$$
\int\na u_{0,C}u_{1,C}=0.
$$
\end{prop}
\begin{proof}
Because of Proposition \ref{p:4.14}, we can assume that
$T_-(\rt)=+\infty$. To abbreviate the notation, let us denote
$u(x,t)=u_C(x,t)$. Again, without loss of generality, if the
conclusion does not hold, we can assume that
$
\ga=\int\dd_{y_1}u_0.u_1>0
$
and hence, by \eqref{eq:4.15}, for all $s\in\bR$ we have
$$
\int\dd_{y_1}u(s)\dd_su(s)=\ga>0.
$$
We will see that this assumption leads to a contradiction. We first
start out by showing: given $\eps>0$,
 \begin{eqnarray}\label{eq:4.22}
&&\text{there exists }R_0(\eps)>0\text{ such
that,}\text{for all }s\geq 0,\text{ we have }\\
&& \int_{\abs{y+\frac{y(s)}{\la(s)}}\geq R_0(\eps)}\abs{\dd_su}^2+\abs{\na_y
u}^2+\frac{\abs{u}^2}{\abs{y}^2}+\abs{u}^{2^*}\leq\eps.\nonumber
 \end{eqnarray}
In fact, by compactness of $\overline{K}$, given $\eps>0$,
$\exists\tilde{R}_0=\tilde{R}_0(\eps)>0$ such that, 
$$
\forall s\in[0,\infty), \ \int_{\abs{y+\frac{y(s)}{\la(s)}}\geq
\tilde{R}_0/\la(s)}\abs{\dd_su}^2+\abs{\na_y
u}^2+\frac{\abs{u}^2}{\abs{y}^2}+\abs{u}^{2^*}\leq\eps.
$$
Since $\la(s)\geq A_0$, $R_0(\eps)=\tilde{R}_0(\eps)/A_0$ does the
job.\\ 
Next, we show that, as a consequence of \eqref{eq:4.22} we have
good bounds for $\abs{\frac{y(s)}{\la(s)}}$:
 \begin{eqnarray}\label{eq:4.23}
\ \ \text{For }M>0, \text{we have, for all }s\in[0,\infty),
\abs{\frac{y(s)}{\la(s)}}\leq s+M.
 \end{eqnarray}
To verify \eqref{eq:4.23}, recall that, since $E(\rt)=E_C>0$, $\rt$
is not identically $0$ and we have, because of Corollary
\ref{c:3.9},
 \bbea
\inf_{s\geq 0}\int\abs{\na_yu(y,s)}^2+\abs{\dd_su(y,s)}^2dy
\geq C\norm{\rt}^2_{\hl}=B_0>0.
 \eea
Then, use \eqref{eq:4.22} to choose $M_0>0$ so that
$$
\int_{\abs{y+\frac{y(s)}{\la(s)}}\geq M_0}\abs{\na u}^2+\abs{\dd_su}^2\leq
B_0/2;\quad s\in[0,\infty),
$$
to conclude that
$$
\int_{\abs{y+\frac{y(s)}{\la(s)}}\leq M_0}\abs{\na u}^2+\abs{\dd_su}^2\geq
B_0/2;\quad s\in[0,\infty).
$$
Recall now from Lemma \ref{l:2.18}, that there exists $\eps_0>0$ so
that if for some $M_1>0$, we have
 \bbe\label{eq:4.24}
\int_{\abs{x}>M_1}\abs{\na_yu_0}^2+\frac{\abs{u_0}^2}{\abs{y}^2}+\abs{u_1}^2\leq
\eps
 \ee
then
$$
\int_{\abs{x}\geq 2M_1+s}\abs{\na_yu(y,s)}^2+\abs{\dd_su(y,s)}^2\leq
C\eps,
$$
wherever $0<\eps<\eps_0$ and $s\geq 0$. Since we can assume, without
loss of generality, that $y(0)=0$, $\la(0)=1$, in light of
\eqref{eq:4.22} we can  always achieve \eqref{eq:4.24}. We will show
that we can choose $\eps$ so small that
$\abs{\frac{y(s)}{\la(s)}}\leq s+3\max(M_0,M_1)$. If not,
$\abs{\frac{y(s)}{\la(s)}}\geq s+3\max(M_0,M_1)$, and if
$\abs{y+\frac{y(s)}{\la(s)}}\leq M_0$, $\abs{y}\geq s+3\max(M_0,M_1)-M_0\geq
s+2\max(M_0,M_1)\geq s+2M_1$. But then,
$$
\frac{B_0}{2}\leq\int_{\abs{y+y(t)/\la(s)}\leq M_0}\abs{\na_y
u}^2+\abs{\dd_su}^2\leq\int_{\abs{y}\geq
s+2M_1}\abs{\na_yu}^2+\abs{\dd_su}^2 \leq C\eps,$$ by
\eqref{eq:4.24}. If $C\eps<B_0/2$, we reach a contradiction, which
establishes \eqref{eq:4.23}.

Having \eqref{eq:4.22}, \eqref{eq:4.23} at our disposal, we now
define for $R>0$, $d>0$,
 \bbe\label{eq:4.25}
z_{d,R}(x_1,\bar{x},t)=u_R\left(\frac{x_1-dt}{\sqrt{1-d^2}},\bar{x},\frac{t-dx_1}{\sqrt{1-d^2}}
\right),
 \ee
where
$$
u_R(y_1,\bar{y},s)=R^{\frac{(N-2)}{2}}u\left(Ry_1,R\bar{y},Rs\right).
$$
Note that $u_R$ is a solution of (CP) in $\bR^N\times\bR$, that
$E((u_R(0),\dd_su_R(0)))=E_C$ and that $\exists\bar{\de}>0$ such
that $\int\abs{\na_yu_R(y,s)}^2dy\leq (1-\bar{\de})\int\abs{\na
W}^2$. We also have $\sup_{s\in\bR}\norm{(u_R,\dd_su_R)}_{\hl}\leq
A$ and $\norm{u_R}_{S((0,+\infty))}=+\infty$. Moreover, we will use
the fact that, when $(x,t)$ are in a compact set, the identity
$\dd_te(z_{d,R})(x,t)=\sum_{j=1}^N\dd_{x_j}(\dd_{x_j}z_{d,R}\dd_tz_{d,R})$ holds, which can be shown by approximating
$u_R$ by compactly supported regular solutions and making the
observation that the corresponding $z_{d,R}$ are then solutions of
(CP) on finite time intervals. 

We now have:
\begin{eqnarray}\label{eq:4.26}
&& \text{There exists }d_0>0\text{ so that, for
}0<d<d_0\nonumber\\
&& \int_1^2\int_{3\leq\abs{x}\leq
8}\abs{\na_xz_{d,R}}^2+\abs{\dd_tz_{d,R}}^2+\abs{z_{d,R}}^{2^*}\leq\eta_1(R,d),\\
&&\text{where } \eta_1(R,d)\underset{R\to\infty}{\xrightarrow{\,\,\,\,\,}}0,
\text{ uniformly in }d<d_0.\nonumber
\end{eqnarray}
To establish \eqref{eq:4.26}, we use the change of variables
$\Phi(x,t)=(y,s)$ where $y_1=\frac{x_1-dt}{\sqrt{1-d^2}}$,
$\bar{y}=\bar{x}$, $s=\frac{t-dx_1}{\sqrt{1-d^2}}$. Then, for $d$
small, we have, after changing variables, that the left-hand side of
\eqref{eq:4.26} is bounded by
$$
\int_{1-1/8}^{2+1/8}\int_{3-1/8\leq\abs{y}\leq8+1/8}
\abs{\na_yu_R}^2+\abs{\dd_su_R}^2+\abs{u_R}^{2^*}dyds,$$ which after
rescaling, becomes
$$
\frac{1}{R}\int_{(1-1/8)R}^{(2+1/8)R}\int_{(3-1/8)R\leq\abs{y}\leq(8+1/8)R}
\abs{\na_yu}^2+\abs{\dd_su}^2+\abs{u}^{2^*}dyds.
$$
But, by \eqref{eq:4.23}, $\abs{\frac{y(s)}{\la(s)}}\leq(2+1/8)R+M$, for $0\leq
s\leq (2+1/8)R$, so that, for $R$ large,
$\left\{y:\abs{y}\geq(3-1/8)R \right\}\subset
\left\{y:\abs{y+\frac{y(s)}{\la(s)}}\geq R/2 \right\}$ and our claim then
follows from \eqref{eq:4.22}.

We now pick $\theta_1=\theta_1(\al)\in C_0^{\infty}(\abs{\al}<5)$, $\theta_1\equiv
1$ on $\abs{\al}<4$, $0\leq\theta_1\leq 1$, and define
$\theta(x)=\theta_1(x_1)\theta_1(\abs{\bar{x}})$. Note that
$\theta(x)\equiv 1$ on $\abs{x}\leq 4$ and
$\supp\theta\subset\left\{\abs{x}\leq\sqrt{25+25}=\sqrt{50} \right\}$. Our next
task is to study
$$\int_1^2\int\theta^2e(z_{d,R})(x_1,\bar{x},t)dx_1d\bar{x}dt.$$ For
this, we first change variables $\Phi(x,t)=(y,s)$, as before. Our
integral then becomes, (for $d_0$ small), recalling that
$t\sqrt{1-d^2}=dy_1+s$,
 \bbe\label{eq:4.27}
\int_{\sqrt{1-d^2}\leq dy_1+s\leq
2\sqrt{1-d^2}}\int\theta^2\left(\Phi^{-1}(y,s)\right)e(z_{d,R})\left(\Phi^{-1}(y,s)\right)dy_1d\bar{y}ds.
 \ee
Recall that $\theta^2(x)=\theta_1^2(x_1)\theta^2_1(\abs{\bar{x}})$,
$\bar{y}=\bar{x}$, $y_1=\frac{x_1-dt}{\sqrt{1-d^2}}$,
$s=\frac{t-dx_1}{\sqrt{1-d^2}}$, so that
$x_1=\frac{y_1+ds}{\sqrt{1-d^2}}$ and
$\theta^2\left(\Phi^{-1}(y,s)\right)=\theta_1^2\left(\frac{y_1+ds}{\sqrt{1-d^2}}\right)
\theta_1^2(\bar{y}).$ Note that $$\abs{e(z_{d,R})\left(
\Phi^{-1}(y,s)\right)}\leq
C\left\{\abs{\na_yu_R(y,s)}^2+\abs{\dd_su_R(y,s)}^2+\abs{u_R^{2^*}(y,s)}
\right\},$$ for $0<d\leq d_0$, and that
$$\theta_1^2\left(\frac{y_1+ds}{\sqrt{1-d^2}}
\right)=\theta_1^2\left(\sqrt{1-d^2}y_1+\frac{d}{\sqrt{1-d^2}}(s+dy_1)
\right).$$ Thus, since in our domain of integration we have
$\sqrt{1-d^2}\leq dy_1+s\leq 2\sqrt{1-d^2}$, for $0\leq d<d_0$,
$d_0$ small, we have
 \bbea
&& \theta_1^2\left(\sqrt{1-d^2}y_1+\frac{d}{\sqrt{1-d^2}}(s+dy_1)
\right)-\theta_1^2(\sqrt{1-d^2}y_1)\\
&& =O(d)(\theta_1^2)'\left( \sqrt{1-d^2}y_1+\eta O(d)\right),
 \eea
where $\abs{\eta}\leq 1$. Note that
$\supp(\theta_1^2)'(\al)\subset\left\{4\leq\abs{\al}\leq 5
\right\}$, so that, for $d_0$ small, this can only be nonzero for
$3+\frac{1}{4}\leq\abs{y_1}\leq 5+\frac{1}{4}$. Using a similar
argument for
$\theta_1^2\left(\sqrt{1-d^2}y_1\right)-\theta_1^2(y_1),$ and the
argument used in the proof of \eqref{eq:4.26}, we see that the
integral in \eqref{eq:4.27} equals
 \bbe\label{eq:4.28}
\underset{\sqrt{1-d^2}\leq dy_1+s\leq
2\sqrt{1-d^2}}{\int\int}\theta^2(y)e(z_{d,R})(\Phi^{-1}(y,s))dyds+d\eta_2(R,d),
 \ee
where
$\abs{\eta_2(R,d)}\underset{R\to\infty}{\xrightarrow{\,\,\,\,\,}}0$,
uniformly in $d<d_0$.

Now, using the formulas after \eqref{eq:4.17}, we see that the
integral in \eqref{eq:4.28} equals $J_1+J_2$, where
 \bbea
J_1&=&\underset{1\leq \frac{dy_1+s}{\sqrt{1-d^2}}\leq
2}{\int\int}\theta^2(y)\left\{
\frac{1}{2}\left[\frac{1+d^2}{1-d^2}(\dd_su_R)^2+\frac{1+d^2}{1-d^2}(\dd_{y_1}u_R)^2+
\abs{\na_{\bar{y}}u_R}^2\right]\right.\\
&&-\left.\frac{1}{2^*}\abs{u_R}^{2^*}\right\}dy ds,\\
J_2&=&-\frac{2d}{(1-d^2)}\underset{1\leq \frac{dy_1+s}{\sqrt{1-d^2}}\leq
2}{\int\int}\theta^2(y)\dd_{y_1}u_R\dd_su_Rdyds.
 \eea
Let us first analyze $\frac{J_2}{d}$. We clearly have
$$
\frac{J_2}{d}=\underset{1\leq \frac{dy_1+s}{\sqrt{1-d^2}}\leq
2}{-2\int\int}\theta^2(y)\dd_{y_1}u_R\dd_su_Rdyds+O(d^2),
$$
where $O(d^2)$ is uniform in $R$. Consider
 \bbea
 \tilde{J_2} &=& \underset{1\leq \frac{s}{\sqrt{1-d^2}}\leq
2}{-2\int\int}\theta^2(y)\dd_{y_1}u_R\dd_su_Rdyds \\
&=& \underset{1\leq \frac{s}{\sqrt{1-d^2}}\leq
2}{-2\int\int}\dd_{y_1}u_R.\dd_su_Rdyds
-\underset{1\leq \frac{s}{\sqrt{1-d^2}}\leq
2}{2\int\int}[\theta^2(y)-1]\dd_{y_1}u_R\dd_su_R
\\
&=&\tilde{J}_{21}+\tilde{J}_{22}.
 \eea
Note that $\supp[\theta^2-1]\subset\{\abs{y}\geq 4\}$, so that, with
the argument in the proof of \eqref{eq:4.26} we see that
$\tilde{J}_{22}=\eta_3(R,d)$, with $\eta_3\underset{R\to
\infty}{\xrightarrow{\,\,\,\,\,}}0$, uniformly in $d$. Moreover,
\eqref{eq:4.15} and scaling  show that
$\tilde{J}_{21}=-2\ga\sqrt{1-d^2}$, so that
\bbea
 \frac{J_2}{d}=&&\underset{1\leq \frac{dy_1+s}{\sqrt{1-d^2}}\leq
2}{-2\int\int}\theta^2(y)\dd_{y_1}u_R\dd_su_Rdyds 
 + \underset{1\leq \frac{s}{\sqrt{1-d^2}}\leq
2}{2\int\int}\theta^2(y)\dd_{y_1}u_R\dd_su_Rdyds
\\
&&-2\ga\sqrt{1-d^2}
+O(d^2)+\eta_3(R,d).
\eea
We turn to the difference of the two integrals on the right hand
side. It is dominated by
\bbea
&& 2\int\int_{y_1>0}\int_{2\sqrt{1-d^2}-dy_1}^{2\sqrt{1-d^2}}\theta^2(y)\abs{\dd_{y_1}u_R}\abs{\dd_su_R}dsdy\\
&&
+2\int\int_{y_1<0}\int_{2\sqrt{1-d^2}}^{2\sqrt{1-d^2}-dy_1}\theta^2(y)\abs{\dd_{y_1}u_R}
\abs{\dd_su_R}dsdy\\
&&
+2\int\int_{y_1>0}\int_{\sqrt{1-d^2}-dy_1}^{\sqrt{1-d^2}}\theta^2(y)\abs{\dd_{y_1}u_R}
\abs{\dd_su_R}dsdy\\
&&
+2\int\int_{y_1<0}\int_{\sqrt{1-d^2}}^{\sqrt{1-d^2}-dy_1}\theta^2(y)\abs{\dd_{y_1}u_R}
\abs{\dd_su_R}dsdy
=\tilde{A}+\tilde{B}+\tilde{C}+\tilde{D}.
 \eea
We will estimate $\tilde{A}$, the other terms being similar. In our
region of integration, we have $\abs{y_1}\leq 5$. We make the change
of variables in the $s$ integral $h=\frac{2\sqrt{1-d^2}-s}{d}$. We
then have, in our region of integration, $0\leq h\leq y_1$. Thus,
 \bbea
\tilde{A}&\leq&
2d\int\int_{y_1>0}\int_0^{y_1}\theta^2(y)\abs{\dd_{y_1}u_R(y,2\sqrt{1-d^2}-dh)}
\\
&& \abs{\dd_{s}u_R(y,2\sqrt{1-d^2}-dh)} dhdy \\
&\leq&
2d\int_0^5\int\theta^2(y)\abs{\dd_{y_1}u_R(y,2\sqrt{1-d^2}-dh)}\\
&& \abs{\dd_{s}u_R(y,2\sqrt{1-d^2}-dh)}dydh \\ &\leq& 2Ad.
 \eea
We thus have
 \bbe\label{eq:4.29}
\frac{J_2}{d}=-2\ga+O(d)+\eta_3(R,d),
 \ee
where $O(d)$ is uniform in $R$ and $\eta_3(R,d)\underset{R\to
\infty}{\xrightarrow{\,}}0$  uniformly for $0\leq d<d_0$.\\
Next, $\frac{J_1}{d}=\underset{1\leq \frac{dy_1+s}{\sqrt{1-d^2}}\leq
2}{\frac{1}{d}\int\int}\theta^2(y)e(u_R)(y,s)dyds+O(d)$,
where $O(d)$ is uniform in $R$ large. We now consider
$
\tilde{J}_1= \underset{1\leq \frac{s}{\sqrt{1-d^2}}\leq
2}{\int\int}\theta^2(y)e(u_R)(y,s)ds.
$
Note that, arguing as in the case of $\tilde{J}_2$, it is easy to
see that
 \bbe\label{eq:4.3}
\tilde{J}_1=\sqrt{1-d^2}E_C+\eta_4(R,d),\,\,\text{ with
}\,\eta_4(R,d)\underset{R\to\infty}{\xrightarrow{\,\,\,\,\,}}0,
 \ee
uniformly for $0\leq d\leq d_0$.

We next turn to $\frac{J_1-\tilde{J}_1}{d}$, which equals after rescaling to
$$\frac{1}{dR}\left\{\underset{R\leq
\frac{dy_1+s}{\sqrt{1-d^2}}\leq
2R}{\int\int}\theta^2(y/R)e(u)(y,s)dyds-
\right.
\left.
\underset{R\leq \frac{s}{\sqrt{1-d^2}}\leq
2R}{\int\int}\theta^2(y/R)e(u)(y,s)dyds \right\}.
$$
The difference then equals
 \bbea
&&\frac{1}{dR}\left\{
-\int\int_{y_1>0}\int_{2R\sqrt{1-d^2}-dy_1}^{2R\sqrt{1-d^2}}\theta^2(y/R)e(u)(y,s)dyds\right.\\
&&+\left.\int\int_{y_1<0}\int_{2R\sqrt{1-d^2}}^{2R\sqrt{1-d^2}-dy_1}\theta^2(y/R)e(u)(y,s)dyds\right\}\\
&& +\frac{1}{dR}\left\{
\int\int_{y_1>0}\int_{R\sqrt{1-d^2}-dy_1}^{R\sqrt{1-d^2}}\theta^2(y/R)e(u)(y,s)dyds\right.\\
&&-\left.\int\int_{y_1<0}\int_{R\sqrt{1-d^2}}^{R\sqrt{1-d^2}-dy_1}\theta^2(y/R)e(u)(y,s)dyds\right\}\\
&& =\frac{1}{dR}\left\{
-\int\int_{y_1>0}\int_{2R\sqrt{1-d^2}-dy_1}^{2R\sqrt{1-d^2}}\theta^2(y/R)e(u)(y,s)dyds\right.\\
&&+\left.\int\int_{y_1>0}\int_{R\sqrt{1-d^2}-dy_1}^{R\sqrt{1-d^2}}\theta^2(y/R)e(u)(y,s)dyds\right\}\\
&& +\frac{1}{dR}\left\{
\int\int_{y_1<0}\int_{2R\sqrt{1-d^2}}^{2R\sqrt{1-d^2}-dy_1}\theta^2(y/R)e(u)(y,s)dyds\right.\\
&&-\left.\int\int_{y_1<0}\int_{R\sqrt{1-d^2}}^{R\sqrt{1-d^2}-dy_1}\theta^2(y/R)e(u)(y,s)dyds\right\}\\
&& =\frac{1}{dR}(L_1+L_2).
 \eea
We will first study $\frac{1}{dR}L_1$. In the first of its
integrals, we interchange the order of integration, to obtain, in
$\supp\theta^2$,
$$
-\int_{2R\sqrt{1-d^2}-5dR}^{2R\sqrt{1-d^2}}\int_{\frac{2R\sqrt{1-d^2}-s}{d}}^{5R}\int
\theta^2(y/R)e(u)(y,s)d\bar{y}dy_1ds.
$$
We then perform the change of variables $s=2R\sqrt{1-d^2}-5\al dR$,
so that the integral equals
$$
-5dR\int_0^1\int_{5\al R}^{5R}\int
\theta^2(y/R)e(u)(y,2R\sqrt{1-d^2}-5\al dR)d\bar{y}dy_1d\al.
$$
Similarly, the second integral equals
$$
5dR\int_0^1\int_{5\al R}^{5R} \theta^2(y/R)e(u)(y,R\sqrt{1-d^2}-5\al
dR)d\bar{y}dy_1d\al,
$$
so that
 \bbea
&&\frac{1}{dR}L_1=5\left\{-\int_0^1\int_{5\al
R}^{5R}\int\theta^2(y/R)e(u)(y,2R\sqrt{1-d^2}-5\al
dR)dyd\al\right.\\
&&+\left.\int_0^1\int_{5\al
R}^{5R}\int\theta^2(y/R)e(u)(y,R\sqrt{1-d^2}-5\al
dR)dyd\al\right\}\\
&&=-5\int_0^1\int_{5\al
R}^{5R}\int\int_{R\sqrt{1-d^2}}^{2R\sqrt{1-d^2}}\theta^2(y/R)\dd_se(u)(y,s-5\al
dR)dsdyd\al\\
&&=-5\int_0^1\int_{5\al
R}^{5R}\int\int_{R\sqrt{1-d^2}}^{2R\sqrt{1-d^2}}\theta^2(y/R)\sum_{j=1}^N\dd_{y_j}\left(\dd_{y_j}u
(y,s-5\al dR).\right. \\ &&.\left.\dd_su(y,s-5\al dR)\right)
dsdyd\al,
 \eea
by virtue of \eqref{eq:2.17}. Integrating by parts, we obtain that $\frac{1}{dR}L_1$ equals
 \bbea
&&5\int_{0}^{1}\int\int_{R\sqrt{1-d^2}}^{2R\sqrt{1-d^2}}
\theta_1^2(5\al)\theta_1^2(\abs{\bar{y}}/R)\dd_{y_1}u\left(5\al R,\bar{y},s-5\al dR \right)\\ &&\dd_su\left(5\al R,\bar{y},s-5\al dR \right)dsd\bar{y}d\al
+5\sum_{j=1}^N\int_0^1\int_{5\al
R}^{5R}\int\int_{R\sqrt{1-d^2}}^{2R\sqrt{1-d^2}}\\&&\dd_{y_j}
(\theta^2(y/R))\dd_{y_j}u(y,s-5\al
dR)\dd_su(y,s-5\al dR)dsdyd\al.
 \eea
For the second term, note that 
$$\abs{\dd_{y_j}(\theta^2(y/R))}\leq
C/R \  \mbox{and} \ \supp\dd_{y_j}(\theta^2(y/R))\subset\{ 4\leq\abs{y}\leq 8\}.$$
Also,
$R(\sqrt{1-d^2}-5\al d)\leq s-5\al dR\leq2R(\sqrt{1-d^2}-5\al d)$,
$0\leq\al\leq 1$. Therefore, the argument after \eqref{eq:4.26}
shows that the second term, for $0\leq d\leq d_0$ is of the form $\eta_5(R,d)$,
with $\eta_5(R,d) \underset{R\to\infty}{\xrightarrow{\,\,\,\,\,}}0$,
uniformly in $0\leq d\leq d_0$. A similar argument shows that
 \bbea
&&\frac{1}{dR}L_2=5\int_{-1}^0\int\int_{R\sqrt{1-d^2}}^{2R\sqrt{1-d^2}}
\theta^2(5\al)\theta_1^2(\abs{\bar{y}}/R)\dd_{y_1}u\left(5\al
R,\bar{y},s-5\al dR \right)\\&& \dd_su\left(5\al
R,\bar{y},s-5\al dR \right)dsd\bar{y}d\al+\eta_6(R,d),
 \eea
with $\eta_6$ behaving like $\eta_5$. Hence,
 \bbea
&&\frac{J_1-\tilde{J}_1}{d}=5\int_{-1}^1\int\int_{R\sqrt{1-d^2}}^{2R\sqrt{1-d^2}}
\theta_1^2(5\al)\theta_1^2(\abs{\bar{y}}/R)\dd_{y_1}u\left(5\al R,\bar{y},s-5\al dR \right)\\
&&  \dd_su\left(5\al R,\bar{y},s-5\al dR
\right)dsd\bar{y}d\al+\eta_5(R,d)+\eta_6(R,d).
 \eea

Making the change of variables $y_1=5\al R$, $\zeta=s-5\al dR$, the
integral on the right hand side gets transformed and $\frac{J_1-\tilde{J}_1}{d}$ equals
 \bbea
&&\underset{R\leq\frac{\zeta+y_1d}{\sqrt{1-d^2}}\leq
2R}{\frac{1}{R}\int\int\int}
\theta_1^2(y_1/R)\theta_1^2(\abs{\bar{y}}/R)\dd_{y_1}u\left(y_1,\bar{y},y\right)\dd_su\left(y_1,\bar{y},y\right)d\zeta dy_1d\bar{y}\\
&&  
+\eta_5(R,d)+\eta_6(R,d).
 \eea
The calculation of $J_2$ above now yields that the integral equals
$\ga+O(d)+\eta_3(R,d)$, so that
 \bbe\label{eq:4.31}
\frac{J_1-\tilde{J}_1}{d}=\ga+O(d)+\eta_3(R,d)+\eta_5(R,d)+\eta_6(R,d),
 \ee
where $\eta_i(R,d)\underset{R\to\infty}{\xrightarrow{\,\,\,\,\,}}0$,
uniformly for $0\leq d<d_0$ and $O(d)$ is uniform in $R$ large.

Next, we recall that for fixed $R$, $u_R\in
L^{\frac{N+2}{N-2}}_IL^{\frac{2(N+2)}{N-2}}_x$, for any compact time
interval. From this and Lemma \ref{l:2.2} we can see that
$\theta(x)z_{d,R}(x,t)$ is in $C([1,2];\hl)$. Fix now $t_0\in[1,2]$
and recall, from the beginning of the proof that
$\dd_te(z_{d,R})(x,t)=\sum_{j=1}^N\dd_{x_j}(\dd_{x_j}z_{d,R}
\dd_tz_{d,R})$. Hence,
 \bbea
&& \int\theta^2(x)e(z_{d,R})(x,t_0)dx=\int_{1}^2\int\theta^2(x)e(z_{d,R})(x,t)dxdt\\
&&
+\int_1^2\int\theta^2(x)\int_t^{t_0}\sum_{j=1}^N\dd_{x_j}(\dd_{x_j}(z_{d,R})
\dd_tz_{d,R})d\al dxdt\\
&&=\int_1^2\int\theta^2(x)e(z_{d,R})(x,t)dxdt\\&&-\sum_{j=1}^N\int_1^2\int
\int_t^{t_0}\dd_{x_j}(\theta^2(x))\dd_{x_j}(z_{d,R})(x,\al)\dd_tz_{d,R}(x,\al)d\al
dxdt.
 \eea
Because of \eqref{eq:4.26}, the second term equals
$\eta_7(R,d,t_0)$, with \hfill\newline
$\eta_7(R,d,t_0)\underset{R\to\infty}{\xrightarrow{\,\,\,\,\,}}0$,
uniformly in $t_0\in[1,2]$, $0\leq d\leq d_0$. Thus, if
$$E(t_0,d,R)=\int\theta^2(x)e(z_{d,R})(x,t_0)dx,$$ we have (using our previous estimates):
\begin{eqnarray}\label{eq:4.32}
&& E(t_0,d,R)=\eta_7(R,d,t_0)+d\eta_2(R,d)-2\ga d+O(d^2)+d\eta_3(R,d)\\
&&+E_C+\eta_4(R,d)+\ga d+d\eta_3(R,d)+d\eta_5(R,d)+d\eta_6(R,d)\nonumber\\
&&=E_C-\ga d
+d\left\{\eta_2(R,d)+\eta_3(R,d)+\eta_5(R,d)+\eta_6(R,d) \right\}\nonumber\\
&&+\eta_4(R,d)+\eta_7(R,d,t_0)+O(d^2).\nonumber
\end{eqnarray}
We now need to consider
 \bbea
&& \int_1^2\int\theta^2(x)\abs{\na_xz_{d,R}(x,t)}^2dxdt
=\int_1^2\int\theta^2(x)\left\{
\frac{1}{(1-d^2)}\abs{\dd_{y_1}u_R}^2\right.\\ &&\left.
+\abs{\na_{\bar{y}}u_R}^2-
\frac{2d}{(1-d^2)}\dd_{y_1}u_R.\dd_su_R+
+\frac{d^2}{(1-d^2)}\abs{\dd_su_R}^2\right\}dxdt.
 \eea
The arguments used to establish \eqref{eq:4.32} easily yield that
the right hand side equals
$
\int\int_{\sqrt{1-d^2}}^{2\sqrt{1-d^2}}\theta^2(y)\abs{\na_yu_R}^2+O(d),
$
where $O(d)$ is uniform in $R$, i.e.,
 \bbe\label{eq:4.33}
\int_1^2\int\theta^2(x)\abs{\na_xz_{d,R}(x,t)}^2dxdt=\int_{\sqrt{1-d^2}}^{2\sqrt{1-d^2}}\int
\theta^2(y)\abs{\na_yu_R}^2+O(d).
 \ee
Define now $h_{d,R}(x,t)=\theta(x)z_{d,R}(x,t)$. Then,
$$
\abs{\na_xh_{d,R}(x,t)}^2=\theta^2\abs{\na_xz_{d,R}}^2+\abs{\na\theta}^2\abs{z_{d,R}}^2+
2\theta\na\theta.\na z_{d,R}z_{d,R}
$$
and note that the last two terms are supported in $3\leq\abs{x}\leq
8$. Also,
$\abs{h_{d,R}}^{2^*}=\theta^2(x)\abs{z_{d,R}}^{2^*}+\left(\abs{\theta}^{2^*}-\abs{\theta}^2
\right) \abs{z_{d,R}}^{2^*}$ and the last term is supported in
$3\leq\abs{x}\leq 8$.

We are now able to conclude the proof. Choose $d_0$ so that for $0<d<d_0$, uniformly in $R$, we have
$$
\int_0^1\int\theta^2\abs{\na_xz_{d,R}}^2\leq(1-\bar{\de}/2)\int\abs{\na
W}^2,
$$
which we can do because of \eqref{eq:4.33}. Let
$1+\bar{\bar{\de}}=\frac{1-\bar{\de}/4}{1-\bar{\de}/2}$. Consider
 \bbea
S_1&=&S_1(d,R)=\left\{t\in[1,2]:\int\theta^2(x)\abs{\na_xz_{d,R}}^2(x,t)dx\right.\\
&& \left.\leq(1+\bar{\bar{\de}}) (1-\bar{\de}/2)\int\abs{\na
W}^2=(1-\bar{\de}/4)\int\abs{\na W}^2 \right\}.
 \eea
Then $\abs{S_1}\geq\bar{\bar{\de}}/(1+\bar{\bar{\de}})$, for all
$0<d\leq d_0$, $R>0$. Next, choose $d_1$ so small and $R>R_0(d_1)$
so that, for all $t_0\in[1,2]$, $E(t_0,d,R)\leq E_C-\frac{\ga}{2}d_1$.
In addition, we can choose $d_1\leq d_0$. This is possible in view
of \eqref{eq:4.32}. Now, for an $\eps>0$ to be chosen, find
$R_1(\eps)$ so large that for $R\geq R_1(\eps)$, we have
$\eta_1(R,d_1)\leq\eps$, where $\eta_1$ is as in \eqref{eq:4.26}.\\
Consider next the set
\begin{eqnarray*}
S_2 &=& S_2(R,d_1,\eps,M) \\
&=& \left\{t\in[1,2]:\int_{3\leq\abs{x}\leq 8}\abs{\na_xz_{d,R}}^2+
\abs{\dd_tz_{d,R}}^2+\abs{z_{d,R}}^{2^*}\leq M\eps \right\}.
\end{eqnarray*}
Because of \eqref{eq:4.26}, $\abs{S_2}\geq (1-1/M)$, and if we
choose $M=M_{\bar{\de}}$ so large that
$(1-1/M_{\bar{\de}})+\bar{\bar{\de}}/(1+\bar{\bar{\de}})>1$, we can
find $t_0= t_0(R,\eps)\in S_1\cap S_2$. We then have:
\begin{eqnarray}\label{eq:4.34}
&& \int\abs{\na_xh_{d,R}(t_0)}^2\leq\int\theta^2\abs{\na
z_{d,R}(t_0)}^2+CM\eps \\ && \leq(1-\bar{\de}/4)\int\abs{\na
W}^2+CM\eps\leq(1-\bar{\de}/8)\int\abs{\na W}^2,\nonumber
\end{eqnarray}
if we choose $CM\eps\leq\bar{\de}/8\int\abs{\na W}^2$, $R\geq
R_1(\eps)$.
\begin{eqnarray}\label{eq:4.35}
&& \int e(h_{d,R})(t_0)\leq\int\theta^2e(z_{d,R})(t_0)+C\eps M
 \leq E_C-\frac{\ga d_1}{2}+C\eps M,
\end{eqnarray}
for $R\geq R_0(d_1)$, $R\geq R_1(\eps)$. If we now choose $C\eps
M\leq\frac{\ga}{4}d_1$, we have then
 \bbe\label{eq:4.36}
\int e(h_{d,R})(t_0)\leq E_C-\frac{\ga d_1}{4},
 \ee
for all $R>Max(R_0(d_1),R_1(\eps))$,
$\eps=\eps(\ga,d_1,\bar{\de})>0$. Let us now consider $w_R(x,t)$ to
be the solution of (CP) with data at $t=t_0$,
($h_{d_1,R}(t_0),\dd_th_{d_1,R}(t_0)$). In light of the definition
of $E_C$, $w_R(x,t)$ exists for all time and verifies, in view of
Corollary \ref{c:4.9}, \bbe\label{eq:4.37}
\int\int\abs{w_R(x,t)}^{\frac{2(N+1)}{N-2}}dxdt\leq C_{d_1,\ga}, \ee
uniformly for all $R>Max(R_0(d_1),R_1(\eps))$. \\
Next, observe that,
by finite speed of propagation (Remark \ref{r:2.12}),
$w_R(x,t)=z_{d,R}(x,t)$ on $\cup_{-2\leq t\leq 1}(B(0,2+t)\times
t)$. To justify the application of Remark \ref{r:2.12}, we
approximate $\rt$ and hence $(u_{0,R},u_{1,R})$ by
$(u_{0,R}^{(j)},u_{1,R}^{(j)})$ which are in $C_0^{\infty}\times
C_0^{\infty}$. The resulting $u_R^{(j)}$ exists on any finite time
interval, for $j$ large by Remark \ref{r:2.22}, and the
corresponding $z_{d,R}^{(j)}$ are now solutions of (CP) on each
finite time interval. We then have, for $j$ large,
$w_R^{(j)}=z_{d,R}^{(j)}$ on the required set, and a passage to the
limit, (since $x$ and $t$ are in fixed bounded sets, we can apply
Lemma \ref{l:2.2}), gives the required identity. But then,
$$
\underset{{\underset{-2\leq t\leq 1}{\cup}(B(0,2+t)\times
t)}}{\int\int}\abs{z_{d_1,R}}^{\frac{2(N+1)}{N-2}}dxdt\leq C_{d_1,\ga}.
$$
We now use our change of variables $(y,s)=\Phi(x,t)$, and observe
that, (for $d_1$ small enough), $\Phi(\cup_{-2\leq t\leq
1}(B(0,2+t)\times t))\supset\{(y,s):0\leq s\leq 1/4,\abs{y}\leq 1/4
\}$. But then, we obtain
$
\underset{ \substack{
0\leq s\leq 1/4\\
\abs{y}\leq 1/4 }}{\int\int}\abs{u_R}^{\frac{2(N+1)}{N-2}}dyds\leq
C_{d_1,\ga},
$
for all $R\geq Max(R_0(d_1),R_1(\eps))$. If we now rescale the
above interval, we find that for all $R\geq Max (R_0(d_1),R_1(\eps))$,
$$
\underset{ \substack{
0\leq s\leq R/4\\
\abs{y}\leq R/4 }}{\int\int}\abs{u}^{\frac{2(N+1)}{N-2}}dyds\leq
C_{d_1,\ga}.
$$
But, since we have
$
\underset{s\geq 0}{\int\int}\abs{u}^{\frac{2(N+1)}{N-2}}dyds=+\infty,
$
we reach a contradiction, which establishes the proposition.
\end{proof}

\section{Rigidity Theorem. Part 1: Infinite time interval and self-similarity
for finite time intervals}\label{s:5}

In this and the following section we will prove the following:
\begin{thm}\label{t:5.1}
Assume that $(u_0,u_1)\in\h\times L^2$ is such that 
$$E((u_0,u_1))<E((W,0)), \ \ \int\abs{\na u_0}^2<\int\abs{\na W}^2, \int \na u_0 . u_1=0.$$
Let $u$ be the solution of $($CP$)$ with $(u(0),\dd_t
u(0))=(u_0,u_1)$, with maximal interval of existence
$(-T_-(u_0,u_1),T_+(u_0,u_1))$. 
Assume  that there exist $\la(t)>0$, $x(t)\in\bR^N$,
for $t\in[0,T_+(u_0,u_1))$, with the property that \bbea K&=&\biggl\{
\vec{v}(x,t)= \left( \frac{1}{\la(t)^{\frac{(N-2)}{2}}} u\left(
\frac{x-x(t)}{\la(t)},t\right), \frac{1}{\la(t)^{\frac{N}{2}}} \dd_t
u\left(\frac{x-x(t)}{\la(t)},t\right) \right),  \\
&&  t\in[0,T_+(u_0,u_1)) \biggr\}
\eea
has the property that $\overline{K}$ is compact in $\h\times L^2$.\\

Then, $T_+(u_0,u_1)<\infty$ is impossible. \\Ê

Moreover, if
$T_+(u_0,u_1)=+\infty$ and we assume that $\la(t)\geq A_0>0$,
for $t\in[0,\infty)$, we must have $u\equiv 0$.
\end{thm}

\begin{rem}\label{r:5.2}
This Theorem shows the rigidity of $($CP$)$ for optimal small data (consider the solution of $($CP$)$, $u(x,t)=W(x)$). The momentum condition is the ingredient which allows us to treat the nonradial situation and is always true for a radial solution. Lemma \ref{l:4.10} implies that we can choose $x(t)$, $\la(t)$ continuous
in $[0,T_+(u_0,u_1))$. Its proof also shows that we can preserve the
property $\la(t)\geq A_0>0$.
\end{rem}

\mbr

We next turn to the proof of Theorem \ref{t:5.1} in the case when
$$T_+(u_0,u_1)=+\infty, \ \la(t)\geq A_0.$$
Assume that
$(u_0,u_1)\not\equiv(0,0)$. Because of Corollary \ref{c:3.9}, we
have $E((u_0,u_1))=E>0$ and $\sup\limits_{t>0}
\norm{(\na u,\dd_t u)}_{L^2} \leq CE$ as well as, from Theorem
\ref{t:3.8},
 \bbe\label{eq:5.3}
\int \abs{\na_x u(t)}^2 - \abs{u(t)}^{2^*} \geq C_{\overline{\de}}
\int \abs{\na_x u(t)}^2
 \ee
and
 \bbe\label{eq:5.4}
\al \int (\dd_t u)^2 + (1-\al) \left( \int\abs{\na_x u(t)}^2 -
\abs{u(t)}^{2^*} \right) \geq C_\al E,
 \ee
for $0<\al<1$.

\mbr

We will also be applying \eqref{eq:4.22}, which gives the following:
\bbe\label{eq:5.5}
\begin{split}
& \text{Given } \eps>0, \text{ there exists } R_0(\eps)>0
\text{ such that}, 
\text{for all } t\geq 0 
\\
& \int_{\abs{x+\frac{x(t)}{\la(t)}}\geq R_0(\eps)} \abs{\dd_t u}^2
+ \abs{\na_x u}^2 + \frac{\abs{u}^2}{\abs{x}^2} + \abs{u}^{2^*}
\leq \eps E.
\end{split}
\ee (Here we use the assumptions $\la(t)\geq A_0>0$, $E>0$.)

\mbr

We will next summarize some algebraic properties that will be needed
in the sequel. Let us fix $\phi\in C_0^\infty(\bR^N)$, $\phi\equiv 1$
for $\abs{x}\leq 1$, $\phi\equiv 0$ for $\abs{x}\geq 2$, and also
define, for $R>0$, 
$$\phi_R(x)=\phi(x/R), \ \ \psi_R(x)=x\phi(x/R).$$
We will set
\[
r(R) = \int_{\abs{x}\geq R} \frac{\abs{u}^2}{\abs{x}^2} +
\abs{u}^{2^*} + \abs{\na u}^2 + \abs{\dd_t u}^2 dx.
\]
\begin{lem}\label{l:5.6}
The following identities hold: $\text{for all } t\geq 0 $
\begin{enumerate}[i$)$]
\item $\dd_t \left( \int \frac{1}{2}(\dd_t u)^2 + \frac{1}{2}
\abs{\na_x u}^2 - \frac{1}{2^*} \abs{u}^{2^*} \right) = 0$
\item $\dd_t \int \na u.\dd_t u=0$
\item $\dd_t \left( \int \psi_R(x).\na u\dd_t u \right) = -\frac{N}{2}
\int (\dd_t u)^2 + \frac{(N-2)}{2}\left[ \int \abs{\na_x u}^2 -
\abs{u}^{2^*}\right]
 + O(r(R))$
\item $\dd_t \left( \int \phi_R u u_t \right) = \int (\dd_t u)^2
- \int \abs{\na u}^2 + \int \abs{u}^{2^*} + O(r(R))$
\item $\dd_t \left( \int \psi_R \left\{ \frac{1}{2} (\dd_t u)^2
+ \frac{1}{2} \abs{\na_x u}^2 -\frac{1}{2^*} \abs{u}^{2^*} \right\}
\right) = -\int \na u \dd_t u + O(r(R))$.
\end{enumerate}
\end{lem}

Note that i) is Remark \ref{r:2.16}, ii) is \eqref{eq:4.15}, v)
follows from \eqref{eq:2.17}, iv) follows from the arguments in
the proof of Theorem \ref{t:3.10}, and iii) follows by an
integration by parts (and a limiting argument).

\mbr

We now will prove the  lemmas crucial for our purpose. Recall that
we can assume $x(0)=0$.

\begin{lem}\label{l:5.7}
There exist $\eps_1>0$, $C>0$, such that, if $\eps\in(0,\eps_1)$,
there exists $R_0(\eps)$ so that if $R>2R_0(\eps)$, then there
exists $t_0=t_0(R,\eps)$, $0\leq t_0\leq CR$, with the property that for all
$0<t<t_0$ we have $\abs{\frac{x(t)}{\la(t)}}<R-R_0(\eps)$ and
$\abs{\frac{x(t_0)}{\la(t_0)}}=R-R_0(\eps)$.
\end{lem}
\begin{proof}
Since $x(0)=0$, $\la(t)\geq A_0>0$, if not, for all $0<t<CR$ (where $C$ is large), we
have $\abs{\frac{x(t)}{\la(t)}}<R-R_0(\eps)$. Let
\[
z_R(t) = \int \psi_R(x).\na_x u u_t + \bigl( \frac{N}{2}-\al \bigr)
\int \phi_R u u_t, \quad 0<\al<1.
\]
Then, by Lemma \ref{l:5.6} and by \eqref{eq:5.4},
 \bbea
z_R'(t) &=& -\frac{N}{2} \int (\dd_t u)^2 + \frac{(N-2)}{2}
\left(\int \abs{\na u}^2 - \abs{u}^{2^*}\right) + O(r(R)) + \\
&& + \left( \frac{N}{2}-\al \right) \left[
\int (\dd_t u)^2 - \int \abs{\na u}^2 + \int \abs{u}^{2^*} \right]
+ O(r(R)) \\
&=& -\al \int (\dd_t u)^2 - (1-\al) \left[ \int \abs{\na u}^2
- \abs{u}^{2^*} \right] + O(r(R)) \\
&\leq& -C_\al E + O(r(R)) .
 \eea
But, for $\abs{x}\geq R$, $\abs{x+\frac{x(t)}{\la(t)}}\geq R_0(\eps)$,
by our assumption, so that, by \eqref{eq:5.5},
$\abs{r(R)}\leq\widetilde{C}\eps E$. Now, choose $\eps$, so small that
$z_R'(t)\leq\frac{-C_\al E}{2}$. Note that
$\abs{z_R(t)}\leq\widetilde{C}_1 RE$, so that, integrating in $t$, between $0$ and $CR$,
\[
CR \frac{C_\al}{2} E \leq 2\widetilde{C}_1 R E.
\]
This is a contradiction for $C$ large.
\end{proof}

Note that in the radial case,we have $x(t)=0$ (see \cite{15}) and a contradiction follows from Lemma \ref{l:5.7}. This proof, using the momentum, is the algebraic counterpart of virial identity used in \cite{15} for the NLS equation.

\begin{lem}\label{l:5.8}
There exist $\eps_2>0$, $R_1(\eps)>0$, $C_0>0$ such that
if $R>R_1(\eps)$, $t_0=t_0(R,\eps)$ is as in Lemma \ref{l:5.7},
then for $0<\eps<\eps_2$,
$$t_0(R,\eps)\geq  \frac{C_0R}{\eps}.$$
\end{lem}

\begin{proof}
Let for $t\in[0,t_0]$,
$$y_R(t)=\int\psi_R(x) e(u)(x,t) dx.$$ Since $\int \na u_0 u_1 = 0$,
because of ii) in Lemma \ref{l:5.6} and v) in Lemma \ref{l:5.6}, we
have $\abs{y_R'(t)}=O(r(R))$. Since for $0<t<t_0$, if $\abs{x}\geq R$, then
\[
\abs{x + \frac{x(t)}{\la(t)}} \geq R-(R-R_0(\eps)) = R_0(\eps),
\]
we have, integrating in $t$,
\[
\abs{y_R(t_0)-y_R(0)} \leq \widetilde{C}
\eps E t_0.
\]

On the one hand, by \eq{5.5}
\[
\abs{y_R(0)} \leq \widetilde{C} R_0(\eps) E + O(Rr(R_0(\eps))) \leq \widetilde{C}
E \{ R_0(\eps) + \eps R \}.
\]

On the other hand,
\[
\abs{y_R(t_0)} \geq \abs{ \int_{\abs{x+\frac{x(t_0)}{\la(t_0)}}\leq R_0(\eps)}
\psi_R e(u)(t_0) } - \abs{ \int_{\abs{x+\frac{x(t_0)}{\la(t_0)}}\geq R_0(\eps)}
\psi_R e(u)(t_0) }.
\]
In the first integral, $\abs{x}\leq\abs{x+\frac{x(t_0)}{\la(t_0)}}
+ \abs{\frac{x(t_0)}{\la(t_0)}} \leq R$, so that $\psi_R(x)=x$. Note
also that the second integral is bounded by $MR\eps E$. Hence,
\[
\abs{y_R(t_0)} \geq \abs{ \int_{\abs{x+\frac{x(t_0)}{\la(t_0)}}\leq R_0(\eps)}
x e(u)(t_0) } - MR\eps E.
\]
But, $\int_{\abs{x+\frac{x(t_0)}{\la(t_0)}}\leq R_0(\eps)} xe(u)(t_0) $ equals
 \bbea
&& 
-\frac{x(t_0)}{\la(t_0)} \int_{\abs{x+\frac{x(t_0)}{\la(t_0)}}\leq R_0(\eps)}
e(u)(t_0)  + \int_{\abs{x+\frac{x(t_0)}{\la(t_0)}}\leq R_0(\eps)} \left[ x+
\frac{x(t_0)}{\la(t_0)} \right] e(u)(t_0) \\
&&= -\frac{x(t_0)}{\la(t_0)} \int e(u)(t_0) +
\frac{x(t_0)}{\la(t_0)} \int_{\abs{x+\frac{x(t_0)}{\la(t_0)}}\geq R_0(\eps)}
e(u)(t_0) \\ && + \int_{\abs{x+\frac{x(t_0)}{\la(t_0)}}\leq R_0(\eps)} \left[ x+
\frac{x(t_0)}{\la(t_0)} \right] e(u)(t_0).
 \eea
The first term is, in absolute value $[R-R_0(\eps)]E$, while the last
two are bounded in absolute value by $\widetilde{C}(R-R_0(\eps))\eps E
+ \widetilde{C} R_0(\eps) E$. We then find:
\[
\abs{y_R(t_0)} \geq (R-R_0(\eps)) E (1-\widetilde{C}\eps) - MR\eps E
- \widetilde{C} R_0(\eps) E.
\]
The quantity on the right exceeds $\frac{R}{4} E$, if, for $0<\eps<\eps_2$
we have $(1-\widetilde{C}\eps-M\eps)\geq 1/2$ and for $R>R_1(\eps)$ we
have $R/4\geq (1+\widetilde{C}) R_0(\eps)$.

\mbr
Thus,
\[
\frac{R}{4} E - \widetilde{C} E \{ R_0(\eps)+\eps R\} \leq \widetilde{C}
\eps E t_0,
\]
which yields the result for $0<\eps<\eps_2'$ and $R>R_1'(\eps)$.
\end{proof}

\mbr

\begin{proof}[Proof of Theorem \ref{t:5.1},
in the case when $T_+(u_0,u_1)=+\infty$]
 By Lemma \ref{l:5.7}, for $0<\eps<\eps_1$, $R>2R_0(\eps)$ we have
$t_0(R,\eps)\leq CR$, while by Lemma \ref{l:5.8}, for $0<\eps<\eps_2$,
$R>R_1(\eps)$, $t_0(R,\eps)\geq C_0 R/\eps$. Hence, for $R>\max\bigl(
2 R_0(\eps), R_1(\eps) \bigr)$, $\eps<\min(\eps_1,\eps_2)$,
$C_0 R/\eps\leq CR$, which is a contradiction for $\eps$ small.
\end{proof}

We now turn to the start of the analysis of the case $T_+((u_0,u_1))<+\infty$.
By scaling we can assume, without loss of generality, that $$T_+((u_0,u_1))=1.$$
Recall, from Lemma \ref{l:4.11} that
 \bbe\label{eq:5.9}
\la(t) \geq \frac{C_0(K)}{1-t}
 \ee
and, from Lemma \ref{l:4.12}, that (after translation in $x$),
 \bbe\label{eq:5.10}
\begin{split}
\supp u(-,t)\subset B(0,1-t) \ \mbox{and} \
 \supp \dd_t u(-,t)\subset B(0,1-t).
\end{split}
 \ee

\begin{lem}\label{l:5.11}
Let $u$ be as above. Then, there is $C_1(K)>0$ such that
$$\frac{C_1(K)}{1-t}\geq \la(t).$$
\end{lem}
\begin{proof}
Assume not. In light of Lemma \ref{l:4.10}, there exist $t_n\uparrow 1$,
such that $\la(t_n)(1-t_n)\uparrow +\infty$. Consider now
\[
z(t) = \int x\na u \dd_t u + \left( \frac{N}{2}-\al \right)
\int u \dd_t u, \qquad 0<\al<1,
\]
which is defined for $0\leq t<1$ (recall \eqref{eq:5.10}).\\
In view of Lemma \ref{l:5.6}, iii), iv) we have
\[
z'(t) = -\al \int (\dd_t u)^2 - (1-\al) \left[ \int \abs{\na_x u}^2
- \abs{u}^{2^*} \right].
\]
Because of Corollary \ref{c:3.9}, $\bigl( u\not\equiv 0$ since
$T_+((u_0,u_1))=1\bigr)$ we have $E((u_0,u_1))=E>0$,
$\sup\limits_{0<t<1} \norm{ (\na u,\dd_t u) }_{L^2} \leq CE$ and
\[
\al \int (\dd_t u)^2 + (1-\al) \left[ \int \abs{\na_x u}^2 -
\abs{u}^{2^*} \right] \geq C_\al E.
\]
Then, we have
\[
z'(t)\leq -C_\al E, \qquad 0<t<1.
\]
Moreover, \eqref{eq:5.10} and Hardy's inequality give that
$z(t)\underset{t\to 1}{\xrightarrow{\ \ \ \ \ }}0$. Also, the
assumption $\int\na u_0.u_1=0$ and ii) in Lemma \ref{l:5.6} give
that $\int\na u.\dd_t u=0$, $0\leq t<1$.

\mbr

Note that, integrating in $t$, $z(t)\geq C_\al E(1-t)$. We have:
\[
\frac{z(t_n)}{(1-t_n)} = \frac{\int\bigl( x+\frac{x(t_n)}{\la(t_n)}\bigr)
\na u \dd_t u}{(1-t_n)} + \left(\frac{N}{2}-\al\right)
\frac{\int u\dd_t u}{(1-t_n)} \geq C_\al E .
\]
We will show that
 \bbe\label{eq:5.12}
\frac{z(t_n)}{(1-t_n)}\to 0,
 \ee
yielding a contradiction. In fact, for $\eps>0$ given,
\[
\frac{1}{(1-t_n)} \int_{\abs{x+\frac{x(t_n)}{\la(t_n)}}\leq\eps(1-t_n)}
\abs{x+\frac{x(t_n)}{\la(t_n)}} \abs{\na u(t_n)}
\abs{\dd_t u(t_n)} \leq C\eps E.
\]
Next, note that
 \bbe\label{eq:5.13}
\abs{\frac{x(t_n)}{\la(t_n)}} \leq 2(1-t_n).
 \ee
If not,
$
B\bigl( -\frac{x(t_n)}{\la(t_n)}, (1-t_n) \bigr) \cap B(0,1-t_n) = \varnothing,
$
so that
\[
\int_{B\bigl( -\frac{x(t_n)}{\la(t_n)}, (1-t_n) \bigr)} \abs{\na u(x,t_n)}^2 dx=0,
\]
while
 \bbea
&&\int_{\abs{x+\frac{x(t_n)}{\la(t_n)}}\geq(1-t_n)} \abs{\na u(x,t_n)}^2 =
\int_{\abs{\la(t_n)x+x(t_n)}\geq\la(t_n)(1-t_n)} \abs{\na u(x,t_n)}^2 \\
&&= \frac{1}{\la(t_n)^N} \int_{\abs{y}\geq\la(t_n)(1-t_n)}
\abs{\na u\bigl(\frac{y-x(t_n)}{\la(t_n)},t_n\bigr)}^2
\underset{n\to\infty}{\xrightarrow{\ \ \ \ \ }}0,
 \eea
by compactness of $\overline{K}$, since $\la(t_n)(1-t_n)\to+\infty$.
But then, $$E((u(x,t_n),\dd_t u(x,t_n)))\to 0$$ (arguing for $\dd_t u$ in
a similar way) which is a contradiction to $E>0$, and thus establishing
\eqref{eq:5.13}. But then,
 \bbea
&&\frac{1}{(1-t_n)} \int_{\abs{x+\frac{x(t_n)}{\la(t_n)}}\geq\eps(1-t_n)}
\abs{x+\frac{x(t_n)}{\la(t_n)}} \abs{\na u(x,t_n)}  \abs{\dd_t u(x,t_n)}
dx \\
&\leq& 3 \int_{\abs{x+\frac{x(t_n)}{\la(t_n)}}\geq\eps(1-t_n)}
\abs{\na u(x,t_n)}  \abs{\dd_t u(x,t_n)} dx \\
&\leq& \frac{3}{\la(t_n)^N} \int_{\abs{y}\geq\eps(1-t_n)\la(t_n)}
\abs{\na u\bigl(\frac{y-x(t_n)}{\la(t_n)},t_n\bigr)} 
\abs{\dd_t u\bigl(\frac{y-x(t_n)}{\la(t_n)},t_n\bigr)} dy\\
&\underset{n\to\infty}{\xrightarrow{\ \ \ \ \ }}&0,
 \eea
by compactness of $\overline{K}$, and the assumption that
$\la(t_n)(1-t_n)\uparrow+\infty$. This shows \eqref{eq:5.12}
for the first term in $\frac{z(t_n)}{(1-t_n)}$. The second one
gives the same result, using the same argument and the fact that
 \bbea
\frac{1}{(1-t_n)} \int \abs{u(t_n)}\abs{\dd_t u(t_n)} 
\leq
\frac{1}{(1-t_n)} \int \abs{x+\frac{x(t_n)}{\la(t_n)}}
\frac{\abs{u(x,t_n)}}{\abs{x+\frac{x(t_n)}{\la(t_n)}}}
\abs{\dd_t u(x,t_n)},
 \eea
and Hardy's inequality.
\end{proof}

\begin{prop}\label{p:5.14}
Assume that $(u_0,u_1)$ is as in Theorem \ref{t:5.1}, with
$T_+(u_0,u_1)=1$. Then $\supp\na u,\dd_t u\subset B(0,1-t)$
and
\[
\vec{K} = \Bigl\{ (1-t)^{\frac{N}{2}} \bigl( \na u((1-t)x,t), \,
\dd_t u((1-t)x,t) \bigr) \Bigr\}
\]
has compact closure in $L^2(\bR^N)^N\times L^2(\bR^N)$.
\end{prop}
\begin{proof}
We first claim that
\[
(1-t)^{\frac{N}{2}} \bigl( \na u((1-t)(x-x(t)),t), \,
\dd_t u((1-t)(x-x(t)),t) \bigr)
\]
has compact closure in $L^2(\bR^N)^N\times L^2(\bR^N)$. This is
because $C_0(K)\leq(1-t)\la(t)\leq C_1(K)$ and if $\overline{K}$
is compact,
\[
K_1 = \bigl\{ \la^{\frac{N}{2}} \vec{v}(\la x) \,:\, \vec{v}\in\overline{K},
\, c_0\leq\la\leq c_1 \bigr\}
\]
also has the property that $\overline{K}_1$ is compact. Next, let
\[
\widetilde{v}(x,t)= (1-t)^{\frac{N}{2}} \bigl( \na u((1-t)x,t),\,
\dd_t u((1-t)x,t) \bigr),
\]
so that $\widetilde{v}(x,t)=\vec{v}(x+x(t),t)$, where
\[
\vec{v}(x,t)= (1-t)^{\frac{N}{2}} \bigl( \na u((1-t)(x-x(t)),t),\,
\dd_t u((1-t)(x-x(t)),t) \bigr).
\]
Note that, by \eqref{eq:5.10}, $\supp\vec{v}(-,t)
\subset\{x:\abs{x-x(t)}\leq 1\}$. The fact that $E>0$ and the
compactness of $\vec{v}(x,t)$ and preservation of energy now
imply that $\abs{x(t)}\leq C$. But, if $$K_2=\bigl\{ \vec{v}(x+x_0)
\,:\, \vec{v}\in K_1, \abs{x_0}\leq C\bigr\},$$ $\overline{K}_2$ is
also compact and hence the Proposition follows.
\end{proof}

\section{Rigidity Theorem. Part 2: Self-similar variables and
conclusion of the proof of the rigidity theorem}\label{s:6}

In this section our point of departure is Proposition \ref{p:5.14}.\\
For this case, in \cite{15}, we proved an extra decay estimate which allowed us to use the $L^2$ invariance and get a contradiction.\\
Following Merle and Zaag (\cite{22}, see also \cite{1}) we will
introduce self-similar variables to show that a solution as in
Proposition \ref{p:5.14} cannot exist. Merle and Zaag considered
the case of power non-linearities $\abs{u}^{p-1}u$ which have
$p\leq 1+\frac{4}{N-1}$, while here we consider the energy critical
case $p=1+\frac{4}{N-2}$. Nevertheless, many of the calculations in
\cite{22} also apply to our case and one can use an extra Liapunov function. We remark that a similar structure
exists in the case of nonlinear heat equations, as has been used
by Giga and Kohn \cite{8} and others (\cite{MZ}).\\
Again here, we obtain some extra decay estimates which allow us to reduce to an elliptic problem with no solution.

\mbr

We now set,
\[
y=x/(1-t), \quad s=-\log(1-t), \quad 0\leq t<1
\]
and define
 \bbe\label{eq:6.1}
w(y,s,0) = (1-t)^{\frac{N-2}{2}} u(x,t) = e^{-s\frac{(N-2)}{2}} u(e^{-s}y, 1-e^{-s}).
 \ee
Note that $w(y,s,0)$ is defined for $0\leq s<+\infty$, and that
$\supp w(-,s,0)\subset\{\abs{y}\leq 1\}$. We also consider, for
$\de>0$, small,
\[
y= \frac{x}{1+\de-t}, \qquad s = -\log(1+\de-t),
\]
 \bbe\label{eq:6.2}
w(y,s,\de) = (1+\de-t)^{\frac{N-2}{2}} u(x,t) 
= e^{-s\frac{(N-2)}{2}} u(e^{-s}y,1+\de-e^{-s})
 \ee
Note that $w(y,s,\de)$ is defined for $0\leq s<-\log\de$, and that
\[
\supp w(-,\de) \subset \left\{ \abs{y}\leq\frac{e^{-s}-\de}{e^{-s}}
= \frac{(1-t)}{(1+\de-t)} \leq 1 - \de\right\}.
\]
The $w$ solve, in their domain of definition, the equation (see \cite{22}):
 \bbe\label{eq:6.3}
\begin{split}
\dd_s^2 w =& \frac{1}{\rho} \dv \bigl( \rho\na w - \rho (y.\na w)y\bigr)
-\frac{N(N-2)}{4} w \\
+& \abs{w}^{\frac{4}{N-2}} w - 2y \na\dd_s w - (N-1)\dd_s w,
\end{split}
 \ee
where $\rho = (1-\abs{y}^2)^{-\frac12}$.
\begin{lem}\label{l:6.4}
For $\de>0$ fixed, the following hold:For $s\in[0,-\log\de)$,
\begin{enumerate}[i$)$]
\item $\supp w(-,s,\de) \subset \left\{ \abs{y}\leq
\frac{e^{-s}-\de}{e^{-s}} \leq 1 - \de \right\}$
\item[${ }$] $\supp \dd_s w(-,s,\de) \subset \left\{ \abs{y}\leq
\frac{e^{-s}-\de}{e^{-s}} \leq 1 - \de \right\}$
\item  $w(-,s,\de)\in H^1_0(B_1)$ and
\[
\int \abs{w}^{2^*} dy \leq C, \quad \int \abs{\na_y w}^2 < \int
\abs{\na W}^2,
\]
\[
\int \abs{w}^2 + \frac{\abs{w}^2}{(1-\abs{y}^2)^2} \leq C, \quad \int\abs{\dd_s
w}^2\leq C.
\]
\item \[
\int \bigl( \abs{\na w}^2 + \abs{\dd_s w}^2 + \abs{w}^2 + \abs{w}^{2^*}
\bigr) \log\bigl(1/(1-\abs{y}^2)\bigr) dy \leq C \log(\frac{1}{\de}).
\]
\item 
\[
\int \bigl( \abs{\na w}^2 + \abs{\dd_s w}^2 + \abs{w}^2 + \abs{w}^{2^*}
\bigr) (1-\abs{y}^2)^{-1/2} dy \leq \frac{C}{\de^{1/2}}.
\]
\end{enumerate}
\end{lem}

\begin{proof}
The first part of i) was pointed out after \eqref{eq:6.2}. For the second
part, we have, using the notation in \eqref{eq:6.2},
 \bbe\label{eq:6.5}
\begin{split}
\dd_s w(y,s,\de) = & -\frac{(N-2)}{2} e^{-s\frac{(N-2)}{2}} u(e^{-s} y,
1+\de-e^{-s}) \\
&+ e^{-s} e^{-s\frac{(N-2)}{2}} \dd_t u (e^{-s} y, 1+\de-e^{-s})  \\
&- e^{-s} e^{-s\frac{(N-2)}{2}} y.\na u(e^{-s}y, 1+\de-e^{-s})
\end{split}
 \ee
and i) follows from \eqref{eq:5.10}.

\mbr

ii) follows from the support property of $w$, which gives
$w(-,s,\de)\in H^{1,2}_0(B_1)$, a change of variables in $y$ and
\eqref{eq:3.9}, Sobolev embedding and Corollary \ref{c:3.9}, the
Hardy inequality (\cite{5}, for example) and \eqref{eq:6.5}.

\mbr

For iii), iv), note that on $\supp w$, $\supp\dd_s w$, we have
$(1-\abs{y}^2)\geq 1-(1-\de e^s)^2=2\de e^s - \de^2 e^{2s}\geq\de$,
for $\de$ small, $0\leq s<-\log\de$.
\end{proof}

For $w(y,s,\de)$, $\de>0$ as above, we now define (see \cite{22})
 \bbe\label{eq:6.6}
\begin{split}
\tilde{E}(w(s)) = & \int_{B_1} \Bigl\{ \frac{1}{2} \bigl[
(\dd_s w)^2 + \abs{\na w}^2 - (y.\na w)^2 \bigr]  \\
& + \frac{N(N-2)}{8} w^2 - \frac{(N-2)}{2N} \abs{w}^{2^*} \Bigr\}
\frac{dy}{(1-\abs{y}^2)^{1/2}}.
\end{split}
 \ee

\begin{prop}\label{p:6.7}
Let $w=w(y,s,\de)$, $\de>0$ be as above. Then, for $0<s_1<s_2<\log(\frac{1}{\de})$,
the following identities hold:
\begin{enumerate}[i$)$]
\item $\tilde{E}(w(s_2))-\tilde{E}(w(s_1)) = \int_{s_1}^{s_2} \int_{B_1}
\frac{(\dd_s w)^2}{(1-\abs{y}^2)^{3/2}} dy\, ds$
\item $\frac{1}{2} \Bigl[ \int_{B_1} \left( \dd_s w w - \frac{(1+N)}{2}
w^2 \right) \frac{dy}{(1-\abs{y}^2)^{1/2}} \Bigr] \biggl\lvert_{s_1}^{s_2} $
 \bbea
 &=& -\int_{s_1}^{s_2} \tilde{E}(w(s)) ds + \frac{1}{N} \int_{s_1}^{s_2} \int_{B_1} \frac{\abs{w}^{2^*}}{(1-\abs{y}^2)^{1/2}}\\
&& + \int_{s_1}^{s_2} \int_{B_1}
\left\{ (\dd_s w)^2 + \dd_s w y.\na w +
\frac{\dd_s w w\abs{y}^2}{(1-\abs{y}^2)} \right\} \frac{dy}{(1-\abs{y}^2)^{1/2}}.
 \eea
\item $\lim\limits_{s\to\log(\frac{1}{\de})} \tilde{E}(w(s))\leq E=E(u_0,u_1)$.
\end{enumerate}
\end{prop}
\begin{proof}
For i) see the proof of Lemma 2.1 in \cite{22}. For ii), see the proof of (11)
in \cite{22}. We turn to the proof of iii). We analyze term by term, using the
notation in \eqref{eq:6.2}.
 \bbea
&& \int_{B_1} \frac{w^2}{(1-\abs{y}^2)^{1/2}}   \\
 &=&
\int_{\abs{y}<(1-t)/(1+\de-t)} (1+\de-t)^{N-2} \abs{u((1+\de-t)y,t)}^2
\frac{dy}{(1-\abs{y}^2)^{1/2}}  \\
&\leq& C \int_{\abs{x}<1-t} (1+\de-t)^{-2} \abs{u(x,t)}^2
\frac{dx}{\de^{1/2}}  \\
&\leq& \frac{C}{\de^{1/2}(1+\de-t)^2} \left( \int_{\abs{x}<(1-t)}
\abs{u(x,t)}^{2^*} dx \right)^{2/2^*}  (1-t)^{2/N}
\underset{t\to 1}{\xrightarrow{\ \ \ \ \ }}0.
 \eea

 \bbea
&& \int_{B_1} \frac{\abs{w}^{2^*}}{(1-\abs{y}^2)^{1/2}} dy \\
&=& \int_{\abs{y}<(1-t)/(1+\de-t)} (1+\de-t)^N \abs{u((1+\de-t)y,t)}^{2^*}
\frac{dy}{(1-\abs{y}^2)^{1/2}} \\
&=& \int_{\abs{x}<(1-t)} \abs{u(x,t)}^{2^*} \frac{dx}{(1-\abs{y}^2)^{1/2}}.
 \eea
Recall that $\abs{y}^2=\frac{\abs{x}^2}{(1+\de-t)^2}$, and assume that
$1-\eps\de\leq t\leq 1$. Then, we have $\frac{1}{(\eps+1)}\leq
(1-\abs{y}^2)^{1/2}\leq 1$, since $\abs{x}\leq (1-t)$. Thus,
\[
\int_{B_1} \frac{\abs{w}^{2^*}}{(1-\abs{y}^2)^{1/2}} dy \geq
\int_{\abs{x}<1-t} \abs{u(x,t)}^{2^*} dx,
\]
and a similar computation gives that
\[
\int_{B_1} \frac{\abs{\na w}^{2}}{(1-\abs{y}^2)^{1/2}} dy \leq
\frac{1}{(1+\eps)^{1/2}} \int_{\abs{x}\leq 1-t} \abs{\na u}^{2^*} dx.
\]
Also,
 \bbea
 && \int_{B_1} (y.\na w)^2 \frac{dy}{(1-\abs{y}^2)^{1/2}} =
\int_{\abs{x}\leq(1-t)} \frac{\abs{x.\na_x u(x,t)}^2}{(1+\de-t)^2}
\cdot \frac{dx}{(1-\abs{y}^2)^{1/2}}  \\
&& \leq \frac{1}{(1+\eps)} \int_{\abs{x}\leq (1-t)} \abs{\na_x
u(x,t)}^2 dx  \frac{\abs{1-t}^2}{(1+\de-t)^2} \underset{t\to
1}{\xrightarrow{\ \ \ \ \ }}0.
 \eea
With these computations and \eqref{eq:6.5} we see that
\[
\lim\limits_{t\to 1} \frac{1}{2} \int \abs{\dd_s w}^2
\frac{dy}{(1-\abs{y}^2)^{1/2}} = \frac{1}{2} \int \abs{\dd_t u}^2
dx,
\]
which combined with the previous calculations yields iii).
\end{proof}

\begin{cor}\label{c:6.8}
For $s\in[0,\log(\frac{1}{\de}))$, we have 
$$-C/\de^{1/2} \leq \tilde{E}(w(s))\leq E.$$
\end{cor}
\begin{proof}
The first statement follows from Proposition \ref{p:6.7} i), iii), while the
second one follows from Lemma \ref{l:6.4}, iv) and \eqref{eq:6.6}.
\end{proof}
Using space-time estimates, we now obtain our first improvement of the space decay of $w$.
\begin{lem}\label{l:6.9}
For $\de>0$, we have
\[
\int_0^1 \int \frac{\abs{\dd_s w}^2}{(1-\abs{y}^2)} dy \, ds
\leq C \log(\frac{1}{\de}).
\]
\end{lem}

\begin{proof}
We start out with the readily verified identity
 \bbea
 && \frac{d}{ds} \biggl\{ \int \Bigl[ \frac{1}{2} (\dd_s w)^2 + \frac{1}{2}
\bigl(\abs{\na w}^2 - (y.\na w)^2 \bigr) + \frac{(N-2)N}{8} w^2 \\
&& - \frac{(N-2)}{2N} \abs{w}^{2^*} \Bigr]  \bigl( -\log(1-\abs{y}^2)
\bigr) dy \biggr\} \\
&& + \int \bigl[ \log(1-\abs{y}^2) +2 \bigr] y.\na w \dd_s w -
\int \log(1-\abs{y}^2) (\dd_s w)^2  - 2 \int (\dd_s w)^2 \\
&& = - 2 \int \frac{(\dd_s w)^2}{(1-\abs{y}^2)}.
 \eea
We now integrate between $0$ and $1$, change signs. In the estimate
of the left hand side, we can drop the term $\int\log(1-\abs{y}^2)(\dd_s w)^2$
since it is negative. The $\frac{d}{ds}$ term, and the $\int_0^1\int(\dd_s w)^2$ term
are controlled by Lemma \ref{l:6.4} (using $-\log(1-\abs{y}^2) \leq C\log(\frac{1}{\de}))$. It remains to bound
 \bbea
&& \abs{ \int_0^1 \int \Bigl[ \log\bigl( 1-\abs{y}^2 \bigr) +2 \Bigr]
y.\na w \dd_s w dy\,ds } \leq \left( \int_0^1 \int
\frac{\abs{\dd_s w}^2}{(1-\abs{y}^2)} \right)^{1/2} 
\\ &&
\times \left( \int_0^1 \int (1-\abs{y}^2) \abs{\log(1-\abs{y}^2)+2}^2
\abs{\na w}^2 dy\,ds \right)^{1/2}.
 \eea
The second factor is bounded because of Lemma \ref{l:6.4} ii). The proof
is concluded by using $ab\leq \eps a^2 + \frac{1}{\eps} b^2$.
\end{proof}

\begin{lem}\label{l:6.10}
For $\de>0$, we have
\begin{enumerate}[i$)$]
\item $\int_0^1\int_{B_1} \frac{\abs{w}^{2^*}}{(1-\abs{y}^2)^{1/2}}\leq C
\bigl( \log(\frac{1}{\de}) \bigr)^{1/2},$
\item $\tilde{E}(w(1))\geq -C \abs{\log(\frac{1}{\de})}^{1/2}$
\end{enumerate}
\end{lem}
\begin{proof}
We will use Proposition \ref{p:6.7}, ii) to handle i). We have
 \bbea
&& \frac{1}{N} \int_0^1 \int_{B_1} \frac{\abs{w}^{2^*}}{(1-\abs{y}^2)^{1/2}}
dy\,ds \\ && = \frac{1}{2} \left[ \int_{B_1} \left( \dd_s w w -\frac{(1+N)}{2} w^2
\right) \frac{dy}{(1-\abs{y}^2)^{1/2}} \right] \biggl\lvert_0^1  + \int_0^1 \tilde{E}(w(s))ds 
\\ && -
\int_0^1 \int_{B_1} \left\{ (\dd_s w)^2
+ \dd_s w y.\na w + \dd_s w \frac{w\abs{y}^2}{(1-\abs{y}^2)}\right\}
\frac{dy\,ds}{(1-\abs{y}^2)^{1/2}}.
 \eea
By Proposition \ref{p:6.7} i) and iii), the second term on the right
hand side is bounded by $E$. The first term on the right hand side is
bounded using Lemma \ref{l:6.4} ii) and Cauchy-Schwarz. For the third
term, because of the sign, we only need to consider the last two summands,
which are bounded in absolute value by
 \bbea
&& \abs{ \int_0^1 \int_{B_1} \frac{\abs{\dd_s w}}{(1-\abs{y}^2)^{1/2}} \left(\frac{\abs{w}}{(1-\abs{y}^2)} + \abs{\na w} \right) dy\,ds}
\\&&\leq 
2\left( \int_0^1 \int_{B_1} \frac{\abs{\dd_s w}^2}{(1-\abs{y}^2)} dy\,ds
\right)^{\frac12} 
\left( \int_0^1 \int_{B_1} \frac{w^2}{(1-\abs{y}^2)^2}+ \abs{\na w}^2 dy\,ds \right)^{\frac12}
 \\
&&\leq C \bigl(\log(\frac{1}{\de})\bigr)^{1/2},
 \eea
because of Lemma \ref{l:6.4} ii) and Lemma \ref{l:6.9}.
This establishes i). 

To prove
ii), we first consider $\int_0^1 \tilde{E}(w(s))ds$, which is bounded from
below by $-C \bigl(\log(\frac{1}{\de})\bigr)^{1/2}$, by i). The monotonicity
of $\tilde{E}$ $\bigl($(i) in Proposition \ref{p:6.7}$\bigr)$ concludes the
proof of ii).
\end{proof}
We now obtain our second improvement of decay on $w$.
\begin{lem}\label{l:6.11}
For $\de>0$, we have
\[
\int_1^{\bigl(\log(\frac{1}{\de})\bigr)^{3/4}} \int_{B_1}
\frac{(\dd_s w)^2}{(1-\abs{y}^2)^{3/2}} dy\,ds \leq C \bigl(\log(\frac{1}{\de})\bigr)^{1/2}.
\]
\end{lem}
\begin{proof}
Because of i) in Proposition \ref{p:6.7}, we have:
\bbea
&& \int_1^{\bigl(\log(\frac{1}{\de})\bigr)^{3/4}} \int_{B_1}
\frac{(\dd_s w)^2}{(1-\abs{y}^2)^{3/2}} dy\,ds = \\
&& = \tilde{E}\Bigl(w\bigl(\log(\frac{1}{\de})^{3/4}\bigr)\Bigr) - \tilde{E}(w(1)) \leq E +
C \bigl(\log(\frac{1}{\de})\bigr)^{1/2},
 \eea
where we have used Corollary \ref{c:6.8} and Lemma \ref{l:6.10} ii).
\end{proof}

\begin{cor}\label{c:6.12}
For each $\de>0$, there exists $\overline{s}_\de\in \bigl(
1, (\log(\frac{1}{\de}))^{3/4} \bigr)$ such that
\[
\int_{\sde}^{\sde+(\log(\frac{1}{\de}))^{1/8}} \int_{B_1} \frac{(\dd_s w)^2}{(1-\abs{y}^2)^{3/2}}
dy\,ds \leq \frac{2C}{\bigl( \log(\frac{1}{\de}) \bigr)^{1/8}}.
\]
\end{cor}

\begin{proof}
Split the interval $\bigl(1, (\log(\frac{1}{\de}))^{3/4} \bigr)$ into disjoint
intervals of length $\bigl(\log(\frac{1}{\de})\bigr)^{1/8}$. The number of such
intervals is of the order of $\bigl(\log(\frac{1}{\de})\bigr)^{5/8}$. For at
least one such interval $\bigl(\sde,\sde+\bigl(\log(\frac{1}{\de})\bigr)^{1/8}\bigr)$,
\hfill\newline with $\sde\in\bigl(1, (\log(\frac{1}{\de}))^{3/4} \bigr)$, we must have
$$\int_{\sde}^{\sde+\bigl(\log(\frac{1}{\de})\bigr)^{1/8}} \int_{B_1}
\frac{(\dd_s w)^2}{(1-\abs{y}^2)^{3/2}} dy\,ds  \leq
\frac{2C\bigl(\log(\frac{1}{\de})\bigr)^{1/2}}{\bigl( \log(\frac{1}{\de})\bigr)^{5/8}}
 = \frac{2C}{\bigl(\log(\frac{1}{\de}))^{1/8}},
$$
where $C$ is the constant in Lemma \ref{l:6.11}, which proves the
Corollary.
\end{proof}

\begin{rem}\label{r:6.13}
Let $\sde=-\log(1+\de-\overline{t}_\de)$. Note that
\[
\abs{ \frac{(1-\overline{t}_\de)}{1+\de-\overline{t}_\de} -1 } =
\frac{\de}{(1+\de-\overline{t}_\de)} = \frac{\de}{e^{-\sde}} \leq
\de^{1/4} \underset{\de\to 0}{\xrightarrow{\ \ \ \ \ }}0.
\]
\end{rem}
\mbr
Let us now reduce the time evolution problem to a stationary problem in the $w$ variable (i.e. self-similar solutions).
Pick $\de_j\downarrow 0$, so that
\[
\Bigl( (1-\overline{t}_{\de_j})^{\frac{N}{2}} \na u\bigl( (1-\overline{t}_{\de_j})y,
\overline{t}_{\de_j}\bigr), \, (1-\overline{t}_{\de_j})^{\frac{N}{2}} \dd_t u
\bigl( (1-\overline{t}_{\de_j})y, \overline{t}_{\de_j}\bigr) \Bigr)
\longrightarrow \bigl( \na u_0^*, u_1^* \bigr)
\]
in $L^2$. This is possible by Proposition \ref{p:5.14}. Note that,
because of Remark \ref{r:6.13} and the compact closure of
$\overrightarrow{K}$ in Proposition \ref{p:5.14}, we also have that
 \bbea
&& \Bigl( (1+\de_j-\overline{t}_{\de_j})^{\frac{N}{2}} \na u\bigl( (1+\de_j-\overline{t}_{\de_j})y,
\overline{t}_{\de_j}\bigr), \\
&& (1+\de_j-\overline{t}_{\de_j})^{\frac{N}{2}} \dd_t u
\bigl( (1+\de_j-\overline{t}_{\de_j})y, \overline{t}_{\de_j}\bigr) \Bigr)
\longrightarrow \bigl( \na u_0^*, u_1^* \bigr) \text{ in } L^2.
 \eea
Let now $u_j^*$, $u^*$ be solutions of (CP) with data
\[
\Bigl( (1+\de_j-\overline{t}_{\de_j})^{\frac{(N-2)}{2}} u\bigl(
(1+\de_j-\overline{t}_{\de_j})y, \overline{t}_{\de_j}\bigr), \,
(1+\de_j-\overline{t}_{\de_j})^{\frac{N}{2}} \dd_t u \bigl(
(1+\de_j-\overline{t}_{\de_j})y, \overline{t}_{\de_j}\bigr) \Bigr)
\]
and $(u_0^*,u_1^*)$ respectively, in a time interval $[0,T^*]$, independent
of $j$, which we take to have $T^*<1$. By uniqueness in the (CP), we have
 \bbe\label{eq:6.14}
\begin{split}
u_j^*(y,\tau) = & (1+\de_j-\overline{t}_{\de_j})^{\frac{(N-2)}{2}}  \\
 &  u\bigl( (1+\de_j-\overline{t}_{\de_j})y, \overline{t}_{\de_j}
+ (1+\de_j-\overline{t}_{\de_j})\tau \bigr).
\end{split}
 \ee
Note that,
$
\supp u_j^*(-,\tau) \subset \bigl\{
\abs{(1+\de_j-\overline{t}_{\de_j})y} \leq 1-\overline{t}_{\de_j}
- (1+\de_j-\overline{t}_{\de_j})\tau \bigr\}
$
and hence $\abs{y}\leq\frac{1-\overline{t}_{\de_j}}{(1+\de_j-\overline{t}_{\de_j})}
-\tau < 1-\tau$ on the support of $u_j^*(-,\tau)$. Similarly, 
$$
\supp \dd_\tau u_j^*(-,\tau) \subset \left\{
y\,:\, \abs{y}\leq \frac{(1-\overline{t}_{\de_j})}{(1+\de_j-\overline{t}_{\de_j})}
-\tau < 1-\tau \right\}.
$$
Let us compare the solutions in the self-similar 
variables. Recall from \eqref{eq:6.2}, that if $s=-\log(1+\de_j-t)$, then
\[
w(y,s,\de_j) = (1+\de_j-t)^{\frac{(N-2)}{2}} u\bigl( (1+\de_j-t)y, t\bigr).
\]
Define now $\tau$ by $t=\overline{t}_{\de_j}+(1+\de_j-\overline{t}_{\de_j})\tau$,
so that $(1+\de_j-t)=(1+\de_j-\overline{t}_{\de_j})(1-\tau)$. Define also
$s=-\log\bigl( (1+\de_j-\overline{t}_{\de_j})(1-\tau) \bigr)$. We then have
 \bbe\label{eq:6.15}
\begin{split}
w(y,s,\de_j) = & \bigl[ (1+\de_j-\overline{t}_{\de_j})(1-\tau) \bigr]^{\frac{(N-2)}{2}}  \\
 &  u\bigl( (1+\de_j-\overline{t}_{\de_j})(1-\tau)y, \overline{t}_{\de_j}
+ (1+\de_j-\overline{t}_{\de_j})\tau \bigr).
\end{split}
 \ee
If we now set 
$$s'=-\log(1-\tau), y'=\frac{y}{(1-\tau)} \ \mbox{and} \
w_j^*(y',s')=(1-\tau)^{\frac{(N-2)}{2}} u_j^*(y,\tau),$$
then $w_j^*$ is a
solution of \eqref{eq:6.3}, for $0<\tau<T^*$. But, because of
\eqref{eq:6.14}, \eqref{eq:6.15},
 \bbea
w_j^*(y',s') &=& (1-\tau)^{\frac{(N-2)}{2}} (1+\de_j-\overline{t}_{\de_j})^{\frac{(N-2)}{2}}
 \\ &&  u\bigl( (1+\de_j-\overline{t}_{\de_j})y, \overline{t}_{\de_j}
+ (1+\de_j-\overline{t}_{\de_j})\tau \bigr)  \\
&=& w(y',s,\de_j),
 \eea
where
$
s = -\log(1+\de_j-t) = -\log\bigl( (1+\de_j-\overline{t}_{\de_j})
(1-\tau) \bigr) 
= -\log(1+\de_j-\overline{t}_{\de_j}) - \log(1-\tau) = \overline{s}_{\de_j}
+ s',
$
i.e.,
 \bbe\label{eq:6.16}
w_j^*(y',s') = w(y',\overline{s}_{\de_j}+s',\de_j).
 \ee
Consider also, $$w^*(y',s')=(1-\tau)^{\frac{(N-2)}{2}} u^*(y,\tau).$$ We clearly
have $\supp u^*(-,\tau)\subset\{ \abs{y}\leq(1-\tau) \}$ and $w^*$ solves
\eqref{eq:6.3} for $0<\tau<T^*$. Also, recall that $\bigl( u_j^*(-,\tau),
\dd_\tau u_j^*(-,\tau) \bigr) \to \bigl( u^*(-,\tau), \dd_\tau u^*(-,\tau) \bigr)$
in $\h\times L^2$, uniformly for $\tau\in[0,T^*]$, by continuity in (CP).
But then if $0\leq\tau\leq T^*/2=\widetilde{T}$ and $0\leq s'\leq
-\log(1-\widetilde{T})$, we have that
\[
\bigl( w_j^*(-,s'), \dd_{s'} w_j^*(-,s') \bigr)
\underset{j\to\infty}{\xrightarrow{\ \ \ \ \ }}\bigl( w^*(-,s'),
\dd_{s'}w^*(-,s')\bigr) \ \mbox{in} \
\dot{H}^1_0\times L^2
\]
uniformly for $0\leq s'\leq-\log(1-\widetilde{T})$. But, by \eqref{eq:6.16}, we have:
 \bbe\label{eq:6.17}
\bigl( w(y',\overline{s}_{\de_j}+s',\de_j), \dd_{s'}
w(y',\overline{s}_{\de_j}+s',\de_j) \bigr)
\underset{j\to\infty}{\xrightarrow{\ \ \ \ \ }}
\bigl( w^*(-,s'), \dd_{s'}w^*(-,s')\bigr),
 \ee
in $\dot{H}^1_0\times L^2$, uniformly in $0\leq
s'\leq-\log(1-\widetilde{T})$ and $w^*$ is a solution of
\eqref{eq:6.3} and $\supp\bigl( w^*(-,s'), \dd_{s'}w^*(-,s')\bigr)
\subset\bigl\{ \abs{y}\leq 1\bigr\}$.

\begin{lem}\label{l:6.18}
Let $w^*$ be as above. Then, $$w^*(y',s')=w^*(y')  \ \mbox{and} \ w^*\not\equiv 0.$$
\end{lem}

\begin{proof}
Let $S=-\log(1-\widetilde{T})$ and choose $j$ large. Then
\[
\int_0^S \int_{B_1} \frac{(\dd_{s'}w^*(y',s'))^2}{(1-\abs{y'}^2)^{3/2}}
dy'ds' \leq \varliminf\limits_{j\to\infty} \int_0^S \int_{B_1}
\frac{(\dd_{s'}w(y',\overline{s}_{\de_j}+s',\de_j))^2}{(1-\abs{y'}^2)^{3/2}}
dy'ds'
\]
by \eqref{eq:6.17}. The right hand side is bounded by
\[
\varliminf\limits_{j\to\infty} \int_{\overline{s}_{\de_j}}^{S+\overline{s}_{\de_j}}
\int_{B_1} \frac{(\dd_{s'}w(y',s',\de_j))^2}{(1-\abs{y'}^2)^{3/2}}
dy'ds' \leq 2C \lim\limits_{j\to\infty} 1\Bigl/ \bigl(\log(\frac{1}{\de_j})\bigr)^{1/8}=0,
\]
by Corollary \ref{c:6.12}. This shows that $w^*(y',s')=w^*(y')$. \\
To show that
$w^*\not\equiv 0$, assume $w^*\equiv 0$. Then, by \eqref{eq:6.16} and \eqref{eq:6.17},
we would have $\na_{y'}w(y',\overline{s}_{\de_j},\de_j)\to 0$ in $L^2(\bR^N)$,
so that $(1+\de_j-\overline{t}_{\de_j})^{\frac{N}{2}}\na_y u\bigl(
(1+\de_j-\overline{t}_{\de_j})y, \overline{t}_{\de_j} \bigr)\to 0$ in $L^2(\bR^N)$.
Because of Corollary \ref{c:3.9}, we have, for $0<t<1$,
$
\int_{B_1} \abs{\na u(x,t)}^2 + \abs{\dd_t u(x,t)}^2 dx \geq CE>0.
$
But,
\[
\int_{B_1} \abs{\na u(x,\overline{t}_{\de_j})}^2 dx =
\int \abs{(1+\de_j-\overline{t}_{\de_j})^{\frac{N}{2}} \na_y
u\bigl( (1+\de_j-\overline{t}_{\de_j})y, \overline{t}_{\de_j}\bigr) }^2
dy \longrightarrow 0,
\]
so for $j$ large we obtain
 \bbe\label{eq:6.19}
\int_{B_1} \abs{\dd_t u(x,\overline{t}_{\de_j})}^2 dx \geq CE/2.
 \ee
But, by \eqref{eq:6.17} and the fact that $\dd_{s'}w^*(-,s')=0$, we
see that \hfill\newline $\dd_s w(y',\overline{s}_{\de_j},\de_j)\to 0$ in $L^2(\bR^N)$.
We now use the formula \eqref{eq:6.5}, which gives
 \bbea
\dd_s w(y',\overline{s}_{\de_j},\de_j) &=& -\frac{(N-2)}{2}
(1+\de_j-\overline{t}_{\de_j})^{\frac{(N-2)}{2}}
u\bigl( (1+\de_j-\overline{t}_{\de_j})y', \overline{t}_{\de_j}\bigr) + \\
&& + (1+\de_j-\overline{t}_{\de_j})^{\frac{N}{2}}
\dd_t u\bigl( (1+\de_j-\overline{t}_{\de_j})y', \overline{t}_{\de_j}\bigr) - \\
&& - (1+\de_j-\overline{t}_{\de_j})^{\frac{N}{2}}
y'\na u\bigl( (1+\de_j-\overline{t}_{\de_j})y', \overline{t}_{\de_j}\bigr).
 \eea
From our assumption, we see that, since $\abs{y'}\leq 1$, the $L^2$ norm of
the last term goes to $0$. The same can be said for the $L^2$ norm of the
first term, by Sobolev embedding. But this contradicts \eqref{eq:6.19},
so that $w^*\not\equiv 0$.
\end{proof}

\begin{prop}\label{p:6.20}
Let $w^*$ be as above. Then, $w^*\in H^1_0(B_1)$, \hfill\newline
$\int_{B_1} \frac{\abs{w^*(y)}^2}{(1-\abs{y}^2)^2}<\infty$ and $w^*$
solves the $($degenerate$)$ elliptic equation
 \bbe\label{eq:6.21}
\frac{1}{\rho} \dv \bigl( \rho\na w^* - \rho(y.\na w^*)y \bigr)
 - \frac{N(N-2)}{4} w^* + \abs{w^*}^{\frac{4}{N-2}} w^* = 0,
 \ee
where $\rho(y)=(1-\abs{y}^2)^{-1/2}$. \\ Moreover, $w^*\not\equiv 0$ and
 \bbe\label{eq:6.22}
\int \frac{\abs{w^*(y)}^{2^*}}{(1-\abs{y}^2)^{1/2}} dy + \int \frac{\bigl[\abs{\na w^*(y)}^2-(y.\na w^*(y))^2\bigr]}{(1-\abs{y}^2)^{1/2}} dy<+\infty,
 \ee
\end{prop}
\begin{rem} We will see that \eq{6.22} are the critical estimates which allow us to conclude the proof.
\end{rem}
\begin{proof}
It only remains to prove \eq{6.22}. Because of \eq{6.17}
and Lemma \ref{l:6.18}, to bound the first term in  \eq{6.22} it suffices to show that
\[
\int_{\overline{s}_{\de_j}}^{\overline{s}_{\de_j}+S} \int_{B_1}
\frac{\abs{w(y',s',\de_j)}^{2^*}}{(1-\abs{y'}^2)^{1/2}} dy'ds' \leq C,
\]
where $C$ is independent of $j$. In order to show this, we use ii) in
Proposition \ref{p:6.7}, so that
 \bbea
&& \frac{1}{N} \int_{\overline{s}_{\de_j}}^{\overline{s}_{\de_j}+S} \int_{B_1}
\frac{\abs{w(y',s',\de_j)}^{2^*}}{(1-\abs{y'}^2)^{1/2}} dy'ds' = \int_{\overline{s}_{\de_j}}^{\overline{s}_{\de_j}+S} \tilde{E}(w(s')) ds' 
 \\
&& + \frac{1}{2} \left[ \int_{B_1} \left( \dd_s w w - \frac{(1+N)}{2} w^2
\right) \frac{dy'}{(1-\abs{y'}^2)^{1/2}} \right]
\biggl\lvert_{\overline{s}_{\de_j}}^{\overline{s}_{\de_j}+S} \\
&& - \int_{\overline{s}_{\de_j}}^{\overline{s}_{\de_j}+S} \int_{B_1}
\left\{ (\dd_s w)^2 + \dd_s w y'.\na w + \frac{w\dd_s w \abs{y'}^2}{(1-\abs{y'}^2)}
\right\} \frac{dy'}{(1-\abs{y'}^2)^{1/2}}.
 \eea
The first term of the right hand side is bounded by Corollary \ref{c:6.8}, the second one by
Lemma \ref{l:6.4} ii). To bound the last one we only need to estimate
the last two summands. To bound the last summand, we use Cauchy-Schwarz
to bound it by
\[
\left(\int_{\overline{s}_{\de_j}}^{\overline{s}_{\de_j}+S} \int_{B_1}
\frac{w^2}{(1-\abs{y'}^2)^2} dy'ds' \right)^{1/2} 
\left(\int_{\overline{s}_{\de_j}}^{\overline{s}_{\de_j}+S}
\int_{B_1} \frac{\abs{\dd_s w}^2}{1-\abs{y'}^2} dy'ds'
\right)^{1/2} 
\]
\[
\leq C \bigl(\log(1/\de_j)\bigr)^{-1/16},
\]
by Lemma \ref{l:6.4} ii) and Corollary \ref{c:6.12}.

\mbr

To bound the second summand we use Cauchy-Schwarz to bound it by
\[
\left(\int_{\overline{s}_{\de_j}}^{\overline{s}_{\de_j}+S}
\int_{B_1} \abs{\na w}^2 dy'ds' \right)^{1/2} 
\left(\int_{\overline{s}_{\de_j}}^{\overline{s}_{\de_j}+S}
\int_{B_1} \frac{\abs{\dd_s w}^2}{1-\abs{y'}^2} dy'ds'
\right)^{1/2},
\]
which can be estimated similarly. This proves the first estimate. To prove
the gradient estimate, use Corollary \ref{c:6.8} and the previous proof to
conclude that
\[
\int_{\overline{s}_{\de_j}}^{\overline{s}_{\de_j}+S} \int_{B_1}
\Bigl\{ \abs{\na w(y',s',\de_j)}^2 - (y'.\na w(y',s',\de_j))^2 \Bigr\}
\frac{dy'}{(1-\abs{y'}^2)^{1/2}} \leq C.
\]
Using \eq{6.17} and Lemma \ref{l:6.18}, this leads to
\[
\int_{B_1} \bigl\{ \abs{\na w^*}^2 - (y'.\na w^*)^2 \bigr\}
\frac{dy'}{(1-\abs{y'}^2)^{1/2}} \leq C,
\]
which concludes the proof.
\end{proof}

\mbr

The contradiction which finishes the proof of Theorem \ref{t:5.1} is
then provided by the following elliptic result:
\begin{prop}\label{p:6.24}
Assume that $w\in H^1_0(B_1)$, is such that
\begin{enumerate}[i$)$]
\item $\int\frac{\abs{w(y)}^2}{(1-\abs{y}^2)^2}dy<\infty$ $\bigl($a consequence
of $w\in H^1_0(B_1)\bigr)$
\item $\int\frac{\abs{w(y)}^{2^*}}{(1-\abs{y}^2)^{1/2}}dy +\int\frac{\abs{\na w(y)}^2 - (y.\na w(y))^2}{(1-\abs{y}^2)^{1/2}}
dy <\infty$
\item $w$ verifies the (degenerate) elliptic equation \eq{6.21}.
\end{enumerate}
Then, $w\equiv 0$.
\end{prop}

\begin{proof}
We write again the equation \eq{6.21}, with $\rho=(1-\abs{y}^2)^{-1/2}$:
 \bbe\label{eq:6.25}
\frac{1}{\rho} \dv\bigl(\rho\na w - \rho(y.\na w)y\bigr) -
\frac{N(N-2)}{4} w + \abs{w}^{\frac{4}{N-2}} w = 0.
 \ee
Consider first the linear part
\[
Lw = \frac{1}{\rho} \dv\bigl( \rho\na w - \rho (y.\na w) y \bigr)
 = \frac{1}{\rho} \dv\bigl( \rho (I-y\otimes y) \na w \bigr).
\]
For $\abs{y}<1-\de$, $\de>0$, $L$ is a second order elliptic operator
with smooth coefficients. Thus, the well-known argument of Trudinger
\cite{35} shows that $w\in L^\infty(B_{1-\de})$ and hence $w\in C^2(B_{1-\de})$,
where $B_{1-\de}=\{y:\abs{y}<1-\de\}$, for each $\de>0$. From this and
the classical unique continuation theorem of Aronszajn, Krzywicki and Szarski
(see \cite{2} and \cite{13}, Section 17.2) we see that if $w\equiv 0$ on
$1-\de<\abs{y}<1$, we will have $w\equiv 0$.

\mbr

In order to establish this for $1-\de<\abs{y}<1$, it is convenient to
write our equation in polar coordinates $(r,\te)$, $0<r<\infty$,
$\te\in S^{N-1}$. In these coordinates, \eq{6.25} becomes: ($y=r\te$)
 \bbe\label{eq:6.26}
\begin{split}
& (1-r^2)^{1/2} \frac{\dd}{\dd r} (1-r^2)^{1/2} \frac{\dd w}{\dd r}
 + \frac{1}{r^2} \De_\te w  + \frac{(N-1)}{r} (1-r^2) \frac{\dd w}{\dd r}  \\
&= \frac{N(N-2)}{4} w
- \abs{w}^{\frac{4}{N-2}} w,
\end{split}
 \ee
where $\De_\te$ denotes the spherical Laplacian on $S^{N-1}$.

\mbr

For $1-\de<r<1$, we perform the change of variables $v(s,\te)=w(r(s),\te)$,
with $r(s)=1-\frac{(1-s)^2}{4}$. For suitable $\widetilde{\de}$, we have
$1-\widetilde{\de}\leq s\leq 1$, when $1-\de\leq r\leq 1$. Also,
$r'(s)=\frac{1-s}{2}$, $\frac{r'(s)}{(1-r(s))^{1/2}}=1$. Since
\[
(1+r(s))^{1/2} \frac{\dd}{\dd s} v(s,\te) =
(1-r^2(s))^{1/2} \frac{\dd w}{\dd r}(r(s),\te),
\]
$v$ verifies the equation
 \bbe\label{eq:6.27}
\begin{split}
& \frac{\dd}{\dd s} (1+r(s))^{1/2} \frac{\dd v}{\dd s}
 + \frac{1}{(1+r(s))^{1/2}} \frac{1}{r(s)^2} \De_\te v  \\
& + \frac{(N-1)}{r(s)} (1-r(s)^2)^{1/2} \frac{\dd v}{\dd s} =
\frac{N(N-2)}{4(1+r(s))^{1/2}} v - \frac{\abs{v}^{\frac{4}{N-2}}v}{(1+r(s))^{1/2}}.
\end{split}
 \ee
The advantage of \eq{6.27} is that it is elliptic, not degenerate
elliptic, near $s=1$ (or $r=1$). Moreover, since $(1+r(s))$ is bounded
above and below and smooth, the coefficients in \eq{6.27} are smooth.
We now turn to some estimates for $v$, for $1-\widetilde{\de}\leq s\leq 1$,
$\te\in S^{N-1}$.

\mbr

We first claim that
 \bbe\label{eq:6.28}
\int_{1-\widetilde{\de}<s<1} \int_{S^{N-1}} \abs{v(s,\te)}^{2^*}
ds\,d\te < \infty.
 \ee
In fact, the integral in \eq{6.28} equals
$
\int_{1-\de<r<1} \int_{S^{N-1}} \abs{w(s,\te)}^{2^*}
\frac{dr\,d\te}{(1-r)^{1/2}},
$
which is finite by virtue of ii).

\mbr

Next, we notice that, for $1-\de\leq\abs{y}\leq 1$,
\[
\abs{\na_\te w(y)} \simeq \abs{\na w - \left(\frac{y}{\abs{y}}.\na w\right)
\frac{y}{\abs{y}}}
\]
and
 \bbea
 \abs{\na w}^2 - (y.\na w)^2 
&&= \left(\frac{1}{\abs{y}^2}-1\right)(y.\na w)^2 +
\abs{ \na w -\frac{y}{\abs{y}}.\na w\frac{y}{\abs{y}}}^2 \\
&&= (1-\abs{y}^2)\left(\frac{y}{\abs{y}}.\na w\right) ^2+
\abs{ \na w -\frac{y}{\abs{y}}.\na w \frac{y}{\abs{y}}}^2.
 \eea
Thus, since $w\in H^1_0(B_1)$, ii) holds, we see that
\[
\int_{1-\de<r<1} \int_{S^{N-1}} \abs{\na_\te w(s,\te)}^{2}
\frac{dr\,d\te}{(1-r)^{1/2}}<\infty
\]
and hence
 \bbe\label{eq:6.29}
\int_{1-\widetilde{\de}<s<1} \int_{S^{N-1}} \abs{\na_\te v(s,\te)}^{2}
ds\,d\te < \infty.
 \ee

\mbr

Next, we show that
 \bbe\label{eq:6.30}
\int_{1-\widetilde{\de}<s<1} \int_{S^{N-1}} \abs{\frac{\dd v}{\dd s}(s,\te)}^{2}
\frac{ds\,d\te}{(1-s)} < \infty.
 \ee
This estimate, combined with $v\in H_0^{1,2}(B_1)$ is the one that
forces $v$ to vanish, since it means that the Cauchy data for the
solution $v$ of \eq{6.27} vanishes. This is a
consequence of the fact that $w\in H^1_0(B_1)$ and the degeneracy of
\eq{6.25}. On the other hand, \eq{6.28} and \eq{6.29} show that we
are dealing with a ``standard solution'' to \eq{6.27}. To obtain
\eq{6.30}, change variables. The integral equals
 \bbea
&& \int_{1-\widetilde{\de}<s<1}\int \abs{\frac{\dd w}{\dd r}(r(s),\te)}^2
 \frac{\abs{r'(s)}^2}{(1-s)} ds\,d\te \\
&& = \int_{1-\widetilde{\de}<s<1}\int \abs{\frac{\dd w}{\dd r}(r(s),\te)}^2
 \frac{\abs{r'(s)}}{2} ds\,d\te \\
&& = \int_{1-\de<r<1}\int \abs{\frac{\dd w}{\dd r}(r,\te)}^2
\frac{dr}{2}\,d\te  \leq C \int_{1-\de<r<1}\int \abs{\frac{\dd w}{\dd r}(r,\te)}^2 dr\,d\te.
 \eea
Finally, a similar argument, using i) shows that
 \bbe\label{eq:6.31}
\int_{1-\widetilde{\de}<s<1}\int \frac{\abs{v(s,\te)}^2}{(1-s)^3} ds\,d\te <\infty.
 \ee

\mbr

Once we have the estimates \eq{6.28}, \eq{6.29}, \eq{6.30} and \eq{6.31},
we define
 \bbe\label{eq:6.32}
\widetilde{v}(s,\te) = \begin{cases}
v(s,\te) & 1-\widetilde{\de}<s<1, \\
0 & 1<s<2.
\end{cases}
 \ee
Since $v(s,\te)\in H^1_0(ds\,d\te)$, for $1-\widetilde{\de}<s<1$, in
light of \eq{6.29}, \eq{6.30} and \eqref{eq:6.31}, $\widetilde{v}\in
H^1(ds\,d\te)$, $1-\widetilde{\de}<s<2$, $\te\in S^{N-1}$. We claim
that $\widetilde{v}$ solves \eq{6.27} for $1-\widetilde{\de}<s<2$:
to show this, let $\eta(s,\te)$ be a test function. Let
$\mu_\eps(s)$ be a smooth approximation of the characteristic
function of $s<1$. We have to show that,
\[
\int\int (1+r(s))^{1/2} \frac{\dd\widetilde{v}}{\dd s}
\frac{\dd\eta}{\dd s} = \lim\limits_{\eps\downarrow 0}
\int\int (1+r(s))^{1/2} \frac{\dd v}{\dd s}
\frac{\dd}{\dd s}(\eta\mu_\eps).
\]
But, this reduces to showing that
 \bbea
&& \lim\limits_{\eps\downarrow 0} \abs{\int\int \eta(1+r(s))^{1/2}
\frac{\dd v}{\dd s}
\frac{\dd}{\dd s}\mu_\eps} \\
&& \leq C \frac{1}{\eps} \int_{1-2\eps<s<1-\eps} \int \abs{\eta} (1+r(s))^{1/2}
 \abs{\frac{\dd v}{\dd s}} \\
&& \leq C \int_{1-2\eps<s<1-\eps} \int \abs{\frac{\dd v}{\dd s}}
\frac{ds\,d\te}{(1-s)} \underset{\eps\to 0}{\xrightarrow{\ \ \ \ \ }}0,
 \eea
because of \eq{6.30}. We can now apply Trudinger's argument in the critical case \cite{35}
to $\widetilde{v}$, to show that $\widetilde{v}\in C^2\bigl(
(1-\widetilde{\de}<s<2)\times S^{N-1}\bigr)$. Once we have this,
$\widetilde{v}\equiv 0$ on $1-\widetilde{\de}<s<2$, because of
the fact that $\widetilde{v}\equiv 0$ for $1<s<2$ and the unique
continuation theorem of \cite{1}. (See also \cite{13}, Section 17.2.)
From this, we conclude that $w\equiv 0$, as desired.
\end{proof}

\begin{rem}\label{r:6.32}
One can skip the use of Trudinger's argument in \cite{35} and use directly
the more delicate unique continuation theorem of \cite{14}, or rather its
variable coefficient version, due to T.~Wolff (\cite{36}). \end{rem}
\begin{rem}\label{r:6.323}
For this part of the argument, no size or energy conditions are needed. In addition, in the radial case, Lemma \ref{l:6.4} and one dimensional Sobolev inequalities give that $\tilde{E}(w(0))$ is bounded in absolute value, which allows us to reduce directly to the elliptic problem.\end{rem}

\mbr

The results in this section yield the contradiction which completes
the proof of Theorem \ref{t:5.1}.

\section{Main Theorem}\label{s:7}

In this section we establish our main result (see \cite{PS} and \cite{S}  for the subcritical case, where energy controls yield the result).
\begin{thm}
Let $(u_0,u_1)\in\h\times L^2$, $3\leq N\leq 5$. Assume that
$E\bigl((u_0,u_1)\bigr)<E\bigl((W,0)\bigr)$. Let $u$ be the
corresponding solution of the Cauchy problem, with maximal interval
of existence \hfill\newline $I=\bigl( -T_-(u_0,u_1), T_+(u_0,u_1)
\bigr)$. $($See Definition \ref{d:2.13}.$)$ Then:
\begin{enumerate}[i$)$]
\item If $\int\abs{\na u_0}^2 < \int \abs{\na W}^2$, then

$$I=(-\infty,+\infty) \ \mbox{and} \
\norm{u}_{L^{\frac{2(N+1)}{N-2}}_{xt}}<\infty.$$

\item If $\int\abs{\na u_0}^2 > \int \abs{\na W}^2$, then
$$T_+(u_0,u_1)<+\infty, \ T_-(u_0,u_1)<+\infty.$$
\end{enumerate}
\end{thm}
\begin{rem} $\int\abs{\na u_0}^2 = \int \abs{\na W}^2$ is incompatible with the energy condition from \eq{3.2}. (Indeed, in this case $E\bigl((u_0,u_1)\bigr)\geq E\bigl((W,0)\bigr))$.
\end{rem}
\begin{proof}
To establish i) we argue by contradiction. If not, $E_C$, defined in
Section 4, must satisfy $\eta_0\leq E_C<E((W,0))$. Let $u_C$ be as
in Proposition \ref{p:4.2} and assume that $I_+$ is finite. Then, by
Proposition \ref{p:4.14}, $\int \na u_{0,C}.u_{1,C}=0$. But then we
reach a contradiction from Theorem \ref{t:5.1}. If $I_+$ is
infinite, and $\la(t)\geq A_0>0$, Proposition \ref{p:4.21} shows
that $\int \na u_{0,C}u_{1,C}=0$ and Theorem \ref{t:5.1} gives
$u_C\equiv 0$, a contradiction because $E\bigl((u_C,\dd_t
u_C)\bigr)=E_C\geq\eta_0$. \\
To conclude the proof, we need to reduce to case $\la(t)>A_0>0,$ for $t\geq0$, using the argument
in the proof of Theorem 5.1 of \cite{15} (also see \cite{21} for a
similar proof). Recall that $E\bigl((u_C,\dd_t
u_C)\bigr)=E_C\geq\eta_0>0$. Because of Lemma \ref{l:4.10}, we can
assume that there exist $t_n\uparrow+\infty$ so that $\la(t_n)\to
0$. After possibly redefining $\{t_n\}_{n=1}^\infty$, we can assume
that
\[
\la(t_n) \leq \inf\limits_{t\in[0,t_n]}\la(t).
\]
From Proposition \ref{p:4.2},
 \bbea
&& (w_{0,n}(x),w_{1,n}(x)) = \biggl( \frac{1}{\la(t_n)^{\frac{(N-2)}{2}}}
u_C\left( \frac{x-x(t_n)}{\la(t_n)}, t_n\right), \\
&& \frac{1}{\la(t_n)^{\frac{N}{2}}} \dd_t
u_C\left( \frac{x-x(t_n)}{\la(t_n)}, t_n\right) \biggr) \longrightarrow
(w_0,w_1) \quad\text{in } \h\times L^2.
 \eea
Note that $E((w_0,w_1))=E_C$. Moreover, $\int \abs{\na w_0}^2 < \int
\abs{\na W}^2$, by the corresponding properties of $u_C$ and Theorem
\ref{t:3.8}. Let $w_0(x,\tau)$, $\tau\in(-T_-(w_0,w_1),0]$ be the
corresponding solution of (CP). If $T_-(w_0,w_1)<\infty$, then
Proposition \ref{p:4.2} and Proposition \ref{p:4.14} yield
$\int\na w_0.w_1=0$, and Theorem \ref{t:5.1} and Proposition \ref{p:4.2}
give a contradiction. Hence $T_-(w_0,w_1)=+\infty$. Let $w_n(x,\tau)$
be the solution of (CP), with data $\bigl(w_{0,n}(x),w_{1,n}(x)\bigr)$,
$\tau\in(-T_-(w_{0,n},w_{1,n}),0]$. Because of Remark \ref{r:2.22},
$\varliminf T_-(w_{0,n},w_{1,n})=+\infty$, and for any $\tau\in(-\infty,0]$,
\[
\bigl( w_n(x,\tau),\dd_\tau w_n(x,\tau)\bigr) \longrightarrow
\bigl( w_0(x,\tau), \dd_\tau w_0(x,\tau) \bigr)
\]
in $\h\times L^2$. Note that, by uniqueness in (CP),
for $0\leq t_n + \frac{\tau}{\la(t_n)}$,
 \bbe\label{eq:7.2}
w_n(x,\tau) = \frac{1}{\la(t_n)^{\frac{(N-2)}{2}}} u_C \left(
\frac{x-x(t_n)}{\la(t_n)},\, t_n+\frac{\tau}{\la(t_n)} \right).
 \ee
Note that,
\[
\varliminf\limits_{n} (-\tau_n) = \varliminf\limits_{n} (t_n\la(t_n))
\geq T_-((w_0,w_1)) = +\infty,
\]
so that for all $\tau\in(-\infty,0]$, for $n$ large, $0\leq
t_n + \frac{\tau}{\la(t_n)} \leq t_n$. In fact, if $-\tau_n\to-\tau_0<\infty$,
then
\[
w_n(x,-\tau_n) = \frac{1}{\la(t_n)^{\frac{(N-2)}{2}}} u_C \left(
\frac{x-x(t_n)}{\la(t_n)}, 0 \right),
\]
\[
\dd_\tau w_n(x,-\tau_n) = \frac{1}{\la(t_n)^{\frac{N}{2}}} \dd_t u_C \left(
\frac{x-x(t_n)}{\la(t_n)}, 0 \right),
\]
would converge to $(w_0(x,-\tau_0),\dd_\tau w_0(x,-\tau_0))$ in $\h\times L^2$,
with $\la(t_n)\to 0$, which is a contradiction from $(u_{0,C},u_{1,C})\not\equiv
(0,0)$, $(w_0,w_1)\not\equiv(0,0)$.

\mbr

Next, note that we must have $\norm{w_0}_{S(-\infty,0)}=+\infty$.
Otherwise, by Theorem \ref{t:2.21}, for $n$ large, $T_-(w_{0,n},w_{1,n})=\infty$
and $\norm{w_n}_{S(-\infty,0)}\leq M$, uniformly in $n$, which, in
view of \eq{7.2}, contradicts $\norm{u_C}_{S(0,+\infty)}=+\infty$.

\mbr

Fix now $\tau\in(-\infty,0]$. For $n$ sufficiently large,
$t_n+\frac{\tau}{\la(t_n)}\geq 0$ and $\la(t_n+\frac{\tau}{\la(t_n)})$ is defined. Let
 \bbea
&& \biggl( \frac{1}{\la(t_n+\frac{\tau}{\la(t_n)})^{\frac{(N-2)}{2}}}
 u_C \left( \frac{x-x(t_n+\frac{\tau}{\la(t_n)})}{\la(t_n+\frac{\tau}{\la(t_n)})},
t_n+\frac{\tau}{\la(t_n)} \right), \\
&&  \frac{1}{\la(t_n+\frac{\tau}{\la(t_n)})^{\frac{N}{2}}}
 \dd_t u_C \left( \frac{x-x(t_n+\frac{\tau}{\la(t_n)})}{\la(t_n+\frac{\tau}{\la(t_n)})},
t_n+\frac{\tau}{\la(t_n)} \right) \biggr) \\
&& = \Biggl( \frac{1}{\widetilde{\la}_n(\tau)^{\frac{(N-2)}{2}}}
 w_n \left( \frac{x-\widetilde{x}_n(\tau)}{\widetilde{\la}_n(\tau)},
\tau \right),   \frac{1}{\widetilde{\la}_n(\tau)^{\frac{N}{2}}}
 \dd_\tau w_n \left( \frac{x-\widetilde{x}_n(\tau)}{\widetilde{\la}_n(\tau)},
\tau \right) \Biggr)\in K,
 \eea
with
 \bbe\label{eq:7.3}
\widetilde{\la}_n(\tau) = \frac{\la(t_n+\frac{\tau}{\la(t_n)})}{\la(t_n)} \geq 1, \
 %
\widetilde{x}_n(\tau) = x(t_n+\frac{\tau}{\la(t_n)}) -
\frac{x(t_n)}{\widetilde{\la}_n(\tau)}.
 \ee
Now, since
$
\frac{1}{\la_n^{\frac{N}{2}}} \vec{v}\left( \frac{x-x_n}{\la_n} \right)
\underset{n\to\infty}{\xrightarrow{\ \ \ \ \ }}\vec{\widetilde{v}} \quad\text{in } L^2,
$
with either $\la_n\to 0$ or $+\infty$, or $\abs{x_n}\to\infty$
implies that $\overrightarrow{\widetilde{v}}\equiv 0$, we see that (since
$E_C>0$) we can assume, after passing to a subsequence, that
$\widetilde{\la}_n(\tau)\to\widetilde{\la}(\tau)$,
$1\leq\widetilde{\la}(\tau)<\infty$ and $\widetilde{x}_n(\tau)
\to\widetilde{x}(\tau)\in\bR^N$. But then
\[
\Biggl( \frac{1}{\widetilde{\la}_n(\tau)^{\frac{(N-2)}{2}}}
 w_0 \left( \frac{x-\widetilde{x}_n(\tau)}{\widetilde{\la}_n(\tau)},
\tau \right),   \frac{1}{\widetilde{\la}_n(\tau)^{\frac{N}{2}}}
 \dd_\tau w_0 \left( \frac{x-\widetilde{x}_n(\tau)}{\widetilde{\la}_n(\tau)},
\tau \right) \Biggr)\in \overline{K}.
\]
But then, by Proposition \ref{p:4.21} and Theorem \ref{t:5.1}, $(w_0,w_1)=(0,0)$,
contradicting $E_C=E((w_0,w_1))$. This proves i).

\mbr

For ii) note that if $u_0\in L^2$, this is the result in Theorem \ref{t:3.10}.
The proof of the general case is a modification of the one of Theorem \ref{t:3.10}.
Let $A=\norm{(u_0,u_1)}_{\h\times L^2}>0$. Recall that (from Lemma \ref{l:2.18}
and its proof) there exists $\eps_0>0$ such that, for $0<\eps<\eps_0$, there
exists $M_0=M_0(\eps)$, with
\[
\int_{\abs{x}\geq M_0+t} \abs{\na_x u(x,t)}^2 + \abs{\dd_t u(x,t)}^2
+ \abs{u(x,t)}^{2^*} + \frac{\abs{u(x,t)}^2}{\abs{x}^2} dx \leq\eps,
\]
for $t\in[0,T_+(u_0,u_1))$. Assume that $T_+(u_0,u_1)=+\infty$ to
reach a contradiction. \\
Let $f(\tau)$ be a solution to the differential
inequality ($f\geq 0$)
\bbe\label{eq:7.5}
 f'(\tau) \geq B f(\tau)^{\frac{N-1}{N-2}}, \ \
 f(0) = 1.
\ee
Then, the time of blow-up for $f$ is $\tau_*$, with $\tau_*\leq K_N
B^{-1}$. \\
Consider now, for $R$ large, $\phi\in C_0^\infty(B_2)$,
$\phi\equiv 1$ on $\abs{x}<1$, $0\leq\phi\leq 1$, 
$$y_R(t)=\int
u^2(x,t)\phi(x/R)dx.$$
Then, $y'_R(t) = 2\int u\, \dd_t u\,
\phi(x/R)dx$, and, using the notation in Lemma \ref{l:5.6}, we have
that
\[
y''_R(t) = 2 \left[ \int (\dd_t u)^2 - \abs{\na_x u}^2
+ \abs{u}^{2^*} dx \right] + O(r(R)).
\]
Arguing as in the proof of Theorem \ref{t:3.10}, we find that
\[
y''_R(t) \geq 2 \left[ 1 + \frac{N}{N-2} \right]
\int (\dd_t u)^2 \phi(x/R) + \widetilde{\widetilde{\de}}_0 + O(r(R)).
\]
Choose now $\eps_1$ so small, and $M_0=M_0(\eps_1)$, as above, so that,
for $R>2M_0$, $O(r(R))\leq\eps_1$, $\eps_1\leq\widetilde{\widetilde{\de}}_0/2$.
We then have, for $0<t<R/2$,
 \bbe\label{eq:7.6}
\begin{split}
 y''_R(t) &\geq \widetilde{\widetilde{\de}}_0/2, \\
 y''_R(t) &\geq 2 \left[ 1 + \frac{N}{N-2} \right]
\int (\dd_t u)^2 \phi(x/R).
\end{split}
 \ee
Note also that
 \bbe\label{eq:7.7}
 y_R(0) \leq C M_0^2 A^2 + \eps_1 R^2, \ \
 \abs{y_R'(0)} \leq C M_0 A^2 + \eps_1 R.
 \ee
Let
$
T = \frac{4CM_0 A^2+2\eps_1 R+2R\sqrt{\eps_1}}{\widetilde{\widetilde{\de}}_0}.
$
Then, (if $T<R/2$)
 \bbea
y'_R(T) &\geq& T \frac{\widetilde{\widetilde{\de}}_0}{2} + y'_R(0)
 \geq 2CM_0 A^2 +\eps_1 R + R \sqrt{\eps_1} - C M_0 A^2 -\eps_1 R  \\
&=& CM_0 A^2 + R \sqrt{\eps_1}.
 \eea
Thus, there exists $0<t_1<T$ such that $y_R'(t_1)=CM_0
A^2+R\sqrt{\eps_1}$, and for $0<t<t_1$, we have $y'_R(t)<CM_0 A^2 +
R\sqrt{\eps_1}$. Note that, in light of \eq{7.6},
$y'_R(t)>y'_R(t_1)>0$, $t>t_1$ ($t<R/2$) and also
 $$
y_R(t_1) \leq y_R(0) + \int_0^{t_1} y'_R \leq y_R(0) + t_1
(CM_0 A^2+R\sqrt{\eps_1}) = y_R(0) + t_1 y'_R(t_1).
$$
We next estimate $T=\frac{4CM_0A^2}{\widetilde{\widetilde{\de}}_0} +
\frac{2\eps_1 R}{\widetilde{\widetilde{\de}}_0} +
\frac{2\sqrt{\eps_1}R}{\widetilde{\widetilde{\de}}_0}$. We first
choose $\eps_1$ so small that
$\frac{2\eps_1}{\widetilde{\widetilde{\de}}_0}
+\frac{2\sqrt{\eps_1}}{\widetilde{\widetilde{\de}}_0}\leq
\frac{1}{32K_N}$, where $K_N$ is the constant defined at the
beginning of the proof, and $R$ so large that
$\frac{4CM_0A^2}{\widetilde{\widetilde{\de}}_0}\leq\frac{1}{16
K_N}R$. We then have $T\leq \frac{1}{8K_N} R$. We can also ensure
$T\leq\frac{R}{8}$. Thus,
\[
y_R(t_1) \leq CM_0^2A^2 +\eps_1R^2 + \frac{R}{8K_N}y_R'(t_1).
\]
If we now use the argument in the proof of Theorem \ref{t:3.10},
for the function $\widetilde{y}_R(\tau)=y_R(t_1+\tau)$,
$0\leq\tau\leq R/4$, in light of \eq{7.6}, we see that, for
$0<\tau<R/4$, we have $\log\bigl(\widetilde{y}'_R(\tau)\bigr)'
\geq \frac{(N-1)}{(N-2)} \log\bigl(\widetilde{y}_R(\tau)\bigr)'$,
so that, by integration,
\[
\frac{\widetilde{y}'_R(\tau)}{\widetilde{y}'_R(0)}
\geq \left[ \frac{\widetilde{y}_R(\tau)}{\widetilde{y}_R(0)}
\right] ^{\frac{(N-1)}{(N-2)}} \quad \text{for } 0\leq\tau\leq R/4.
\]
Thus, if $f(\tau)=\frac{\widetilde{y}_R(\tau)}{\widetilde{y}_R(0)}$ and
$B=\frac{\widetilde{y}'_R(0)}{\widetilde{y}_R(0)}=
\frac{y'_R(t_1)}{y_R(t_1)}$, we have that $f$ is a solution of \eq{7.5} for $0\leq\tau\leq R/4.$ Thus, we must have
\[
\frac{R}{4} \leq \frac{y_R(t_1)}{y'_R(t_1)} K_N
\leq \frac{K_N(CM_0^2 A^2 + \eps_1 R^2)}{y_R'(t_1)} +\frac{R}{8},
\]
or
 \bbea
\frac{1}{8} &\leq& \frac{K_N(CM_0^2 A^2 + \eps_1 R^2)}{CM_0A^2 R+\sqrt{\eps_1}R^2}
 = \frac{K_N[CM_0^2 A^2/R^2 + \eps_1]}{[CM_0A^2/R+\sqrt{\eps_1}]}  \\
&\leq& K_N M_0/R + K_N \sqrt{\eps_1}.
 \eea
By taking $K_N\sqrt{\eps_1}<\frac{1}{32}$, and $\frac{K_N M_0}{R}<\frac{1}{32}$ we
reach a contradiction, which gives the proof of ii).
\end{proof}

To conclude, let us give some Corollaries of our main results similarly to the NLS case (We will refer to \cite{15} for the proofs, which are identical).

\begin{cor}\label{5.15}
Let $(u_0,u_1)\in\h\times L^2$, $3\leq N\leq 5$. Assume that
$E\bigl((u_0,u_1)\bigr)<E\bigl((W,0)\bigr)$ and $\int\abs{\na u_0}^2 < \int \abs{\na W}^2$.  Then the solution $u$ of the Cauchy problem (CP) with data $(u_0,u_1)$ at $t=0$ has time interval of existence $I = (-\infty,+\infty)$, and there exists $(u_{0,\pm},u_{1,\pm})$ in $\dot H^1\times L^2$ such that if we denote by $v_\pm(t)$  the solutions of $(LCP)$ corresponding to these initial data, we have
\[
\lim_{t\to \pm\infty}\norm{(u(t),\dd_tu(t)) - (v_{\pm}(t),\dd_tv_{\pm}(t))}_{\dot H^1\times L^2} = 0.
\]
Moreover, if we define $\delta_0$ so that $E(u_0,u_1) \le (1-\delta_0)E\bigl((W,0)\bigr)$, there exists a function $M(\delta_0)$ so that $\norm{u}_{S((-\infty,+\infty))}\leq M(\delta_0).$
\end{cor}

Let us give now a different version of the main result.

\begin{cor}\label{5.16}
Let $(u_{0},u_{1})$ in $\dot H^1\times L^2$  and assume that for all
$t \in (-T_-(u_0),T_+(u_0))$  we have $\int\abs{\na u(t)}^2+\abs{\dd_tu(t)}^2 \le \int\abs{\na W}^2 - \delta_0$, for $\delta_0>0$.  Then the solution $u$ of the Cauchy problem (CP) with data $u_0$ at $t=0$ has time interval of existence $I = (-\infty,+\infty)$, $\norm{u}_{S((-\infty,+\infty))}<+\infty$.
\end{cor}

\begin{cor}\label{5.18}
Let  $3\leq N\leq 5$, $(u_0,u_1)\in\h\times L^2$ (no size restrictions) be such that $T_+((u_0,u_1)) < + \infty$ and \\ $\forall t \in[0,T_+((u_0,u_1))), \int\abs{\na u(t)}^2 +\abs{\dd_tu(t)}^2\le C_0$. Then, we have for $x_1(t)$ and $x_2(t)$, and for all $R>0$,
$$\varliminf_{t\to T_+(u_0)} \int_{|x-x_1(t)| \le R}  \abs{\na u(t)}^2 +\abs{\dd_tu(t)}^2 \ge \frac 2 N \int  \abs{\na W}^2$$
$$\varlimsup_{t\to T_+(u_0)} \int_{|x-x_2(t)| \le R}  \abs{\na u(t)}^2 +\abs{\dd_tu(t)}^2 \ge \int  \abs{\na W}^2.$$
\end{cor}


\begin{thebibliography}{99}

\bibitem{1} C.~Antonini and F.~Merle, {\em Optimal bounds on positive
blow-up solutions for a semilinear wave equation}, Internat. Math.
Res. Notices \textbf{21} (2001), 1141--1167.

\bibitem{2} N.~Aronszajn, A.~Krzywicki and J.~Szarski, {\em A unique continuation
theorem for exterior differential forms on Riemannian manifolds},
Ark. Mat. \textbf{4} (1962), 417--453.

\bibitem{3} T.~Aubin, {\em \'Equations diff\'erentielles non lin\'eaires et
probl\`eme de Yamabe concernant le courbure scalaire}, J. Math.
Pures Appl. (9), \textbf{55}, 1976, 3, 269--296.

\bibitem{4} H.~Bahouri and P.~G\'erard, {\em High frequency approximation of
solutions to critical nonlinear wave equations}, Amer. J. Math
\textbf{121} (1999), 131--175.

\bibitem{5} H.~Br\'ezis and M.~Marcus, {\em Hardy's inequalities revisited},
Ann. Scuola Norm. Piza \textbf{25} (1997), 217--237.

\bibitem{6} M.~Christ and M.~Weinstein, {\em Dispersion of small amplitude
solutions of the generalized Korteweg--de Vries equation}, J. Funct.
Anal. \textbf{100} (1991), 87--109.

\bibitem{7} P.~G\'erard, {\em Description du d\'efaut
de compacit\'e de l'injection de Sobolev}, ESAIM Control Optim.
Calc. Var. \textbf{3} (1998), 213--233.

\bibitem{8} Y.~Giga and R.~Kohn, {\em Nondegeneracy of blowup for
semilinear heat equations}, Comm. Pure Appl. Math. \textbf{42}
(1989), 223--241.

\bibitem{9} J.~Ginibre, A.~Soffer and G.~Velo, {\em The global Cauchy problem for the
critical nonlinear wave equation}, J. Funct. Anal.  \textbf{110}
(1992), 96--130.

\bibitem{10} J.~Ginibre and G.~Velo, {\em Generalized Strichartz inequalities
for the wave equation}, J. Funct. Anal.  \textbf{133} (1995), 50--68.

\bibitem{11} M.~Grillakis, {\em Regularity and asymptotic behaviour of the
wave equation with a critical nonlinearity}, Ann. of Math.
\textbf{132} (1990), 485--509.

\bibitem{12} M.~Grillakis, {\em Regularity for the wave equation
with a critical nonlinearity},  Comm. Pure Appl. Math. \textbf{45}
(1992), 749--774.

\bibitem{13} L.~H\"ormander, ``The analysis of linear partial differential operators III'',
Springer-Verlag, Berlin-Heidelberg-New York-Tokyo, 1984.

\bibitem{14} D.~Jerison and C.~E.~Kenig, {\em Unique continuation and absence of
positive eigenvalues for Schr\"odinger operators}, Ann. of Math.
\textbf{121} (1985), 463--494.

\bibitem{K} L.~Kapitanski, {\em Global and unique weak solutions of nonlinear wave equations}, Math. Res. Lett., \textbf{1} (1994), no. 2, 211--223.


\bibitem{15} C.~E.~Kenig and F.~Merle, {\em Global well-posedness, scattering and blow-up for the
energy critical, focusing, non-linear Schr\"odinger equation in the
radical case}, to appear, Invent. Math.

\bibitem{16} C.~E.~Kenig, G.~Ponce and L.~Vega, {\em Well-posedness and scattering
results for the generalized Korteweg-de Vries equation via the
contraction principle},  Comm. Pure Appl. Math.  \textbf{46} (1993),
527--620.

\bibitem{17} S.~Keraani, {\em On the defect of compactness for the Strichartz
estimates of the Schr\"odinger equations}, J. Differential Equations
\textbf{175} (2001), 353--392.

\bibitem{18} J.~Krieger and W.~Schlag, {\em On the focusing critical semi-linear wave equation},
to appear, Amer. J. of Math.

\bibitem{19} H.~Levine, {\em Instability and nonexistence of global solutions
to nonlinear wave equations of the form
$Pu_{tt}=-Au+\mathcal{F}(u)$}, Trans. Amer. Math. Soc. \textbf{192}
(1974), 1--21.

\bibitem{20} H.~Lindblad and C.~Sogge, {\em
On existence and scattering with minimal regularity for semilinear
wave equations}, J. Funct. Anal. \textbf{130} (1995), 357--426.

\bibitem{21} F.~Merle, {\em Existence of blow-up solutions in the
energy space for the critical generalized KdV equation},
J. Amer. Math. Soc. \textbf{14} (2001), 555--578.

\bibitem{22} F.~Merle and H.~Zaag, {\em Determination of the blow-up
rate for the semilinear wave equation},
Amer. J. of Math. \textbf{125} (2003), 1147--1164.

\bibitem{MZ} F.~Merle and H.~Zaag, {\em A Liouville theorem for vector-valued nonlinear heat equations and  applications}, Math. Ann. \textbf{316} (2000), no. 1, 103--137.

\bibitem{PS} L.E.~Payne and D.H.~Sattinger, {\em Saddle points 
and instability of nonlinear 
hyperbolic equations}, Israel J. Math., \textbf{22}, 
(1975), 273--303.

\bibitem{23} H.~Pecher, {\em Nonlinear small data scattering for the
wave and Klein-Gordon equation},
Math. Z. \textbf{185} (1984), 261--270.

\bibitem{S} D.H.~Sattinger, {\em On global solutions of nonlinear 
hyperbolic equations}, Arch. Rational Mech. Anal., 
\textbf{30}, (1968), 148--172.

\bibitem{24} J.~Shatah and M.~Struwe, {\em Regularity results for
nonlinear wave equations}, Ann. of Math.
\textbf{138} (1993), 503--518.

\bibitem{25} J.~Shatah and M.~Struwe, {\em Well-posedness in the
energy space for semilinear wave equations with critical growth},
Internat. Math. Res. Notices \textbf{7} (1994), 303--309.

\bibitem{26} J.~Shatah and M.~Struwe, ``Geometric wave equations,''
Courant Lecture Notes in Mathematics,
\textbf{2} (1998).

\bibitem{27} C.~Sogge, ``Lectures on nonlinear wave equations,''
Monographs in Analysis II, International Press, 1995.

\bibitem{28} G.~Staffilani, {\em
On the generalized Korteweg-de Vries-type equations}, Differential
Integral Equations \textbf{10} (1997), 777--796.

\bibitem{29} W.~Strauss, ``Nonlinear wave equations,'' CBMS Regional
Conference Series in
Mathematics, \textbf{73}, American Math. Soc., Providence, RI, 1989.

\bibitem{30} M.~Struwe, {\em Globally regular solutions to the $u\sp 5$
Klein-Gordon equation},
Ann. Scuola Norm. Sup. Pisa Cl. Sci. \textbf{15} (1988), 495--513.

\bibitem{31} G.~Talenti, {\em Best constant in Sobolev inequality}, Ann.
Mat. Pura Appl.
(4) \textbf{110} (1976), 353--372.

\bibitem{32} T.~Tao, {\em Spacetime bounds for the energy-critical
nonlinear wave equation in three spatial dimensions}, preprint, {\tt
http://arxiv.org/abs/math.AP/0601164}.

\bibitem{33} T.~Tao and M.~Visan, {\em
Stability of energy-critical nonlinear Schr\"odinger equations in
high dimensions}, Electron. J. Differential Equations \textbf{118}
(2005), 28 pp. (electronic).

\bibitem{34} M.~Taylor, ``Tools for PDE. Pseudodifferential operators,
paradifferential operators and layer
potentials,'' Math. Surveys and Monographs \textbf{81}, AMS,
Providence RI 2000.

\bibitem{35} N.~Trudinger, {\em
Remarks concerning the conformal deformation of Riemannian
structures on compact manifolds}, Ann. Scuola Norm. Sup. Pisa
\textbf{22} (1968), 265--274.

\bibitem{36} T.~Wolff, {\em Recent work on sharp estimates in
second-order elliptic unique continuation problems},
J. Geom. Anal. \textbf{3}  (1993),  621--650.












\end{thebibliography}
\end{document}